\documentclass[12pt,a4paper,openany]{article}

\usepackage{amsmath,amssymb,amsfonts,amsthm}
\usepackage{graphics}
\usepackage{mathrsfs}
\usepackage[all]{xy}
\usepackage{color}
\usepackage{bm}
\usepackage{setspace}
\usepackage{ulem}
\usepackage{amscd} 
\usepackage{hyperref}

\numberwithin{equation}{section}

\newtheorem{thm}{Theorem}[section]
\newtheorem{prop}[thm]{Proposition}
\newtheorem{lemma}[thm]{Lemma}

\newtheorem{cor}[thm]{Corollary}

\newtheorem{conj}[thm]{Conjecture}

\theoremstyle{definition}
\newtheorem{definition}[thm]{Definition}
\newtheorem{rem}[thm]{Remark}
\newtheorem{exam}[thm]{Example}

\newtheorem{hypo}[thm]{Hypothesis}

\newtheorem*{ack}{Acknowledgment}
\newtheorem*{ai}{Disclosure on the use of AI}

\newcommand{\sq}{\hfill $\square$}
\newcommand{\kah}{K\"{a}hler}
\newcommand{\dol}{\sqrt{-1}\partial \overline{\partial}}

\newcommand{\hol}{H\"{o}lder}
\newcommand{\ma}{Monge--Amp$\grave{{\rm e}}$re}

\newcommand{\Vol}[0]{\operatorname{Vol}}

\newcommand{\e}{\varepsilon}

\newcommand{\Ric}{\mathrm{Ric}}
\newcommand{\MA}{\mathrm{MA}}

\renewcommand{\b}{\beta}

\renewcommand{\phi}{\varphi}

\renewcommand{\leq}{\leqslant}
\renewcommand{\geq}{\geqslant}


\makeatletter

\@addtoreset{equation}{section}
\makeatother

\title{Coupled free energy and a diagonal approach to coupled {\kah}--Einstein metrics
}
\author{Takahiro Aoi\thanks{Kobe City College of Technology, Gakuen-Higashimachi, Nishi-ku, Kobe, Hyogo 651-2194, Japan\\  
\href{mailto:aoi@kobe-kosen.ac.jp}{aoi@kobe-kosen.ac.jp},
\href{mailto:takahiro.aoi.math@gmail.com}{takahiro.aoi.math@gmail.com}}
\and Satoshi Nakamura\thanks{Department of Mathematics, Nagoya Institute of Technology, Nagoya 466-8555, Japan\\
\href{mailto:nakamura.satoshi@nitech.ac.jp}{nakamura.satoshi@nitech.ac.jp}
}}
\date{
}

\begin{document}
\maketitle

\footnote{ 
2020 \textit{Mathematics Subject Classification}.
Primary 32Q20, 32Q26; Secondary 32U15, 32P05, 60F10.
}
\footnote{ 
\textit{Key words and phrases}.
Coupled K\"ahler--Einstein metric, Free energy, Large deviation principle,
Gibbs stability, Non-Archimedean pluripotential theory.
}

\begin{abstract}
In this paper, we introduce the coupled free energy $F_{-\gamma}^{\rm cp}$ which is a natural generalization of the free energy functional defined by Berman.
In the Fano case, our coupled free energy provides a diagonal approach to the coercivity of the coupled Ding functional.
We also define the non-Archimedean coupled free energy $F_{-\gamma}^{\rm cp,NA}$ on the space of finite energy measures on the Berkovich analytification.
As a consequence of the regularization process of the non-Archimedean entropy established by Boucksom--Jonsson, Chi Li and Trusiani, we prove that the 
existence of a unique coupled {\kah}--Einstein metric is equivalent to the stability condition that there exists $\epsilon >0$ such that $F_{-1-\epsilon}^{\rm cp,NA} \geq 0$ on the space of measures associated with a test configuration for the anticanonical line bundle.

As a coupled version of the work by Berman, we show that the law of the empirical measure satisfies a large deviation principle whose rate function is the coupled free energy under the assumption that a conjecture by Hultgren is true.
We also introduce a coupled version of the microscopic stability thresholds and prove that a uniformly coupled Gibbs stable Fano manifold has a unique coupled {\kah}--Einstein metric.
This is a generalization of the results by Fujita--Odaka and Berman.

\end{abstract}


\tableofcontents

\section{Introduction}
Let $X$ be a projective algebraic manifold.
The study of canonical {\kah} metrics is a central problem in complex geometry.
When the canonical line bundle $K_X$ is ample, i.e., $X$ is canonically polarized, 
Aubin and Yau proved the existence of a unique {\kah}--Einstein metric with negative Ricci curvature by solving the complex {\ma} equation
\cite{Au,Yau}.
Berman \cite{Ber13} developed a thermodynamical approach to K\"ahler–Einstein metrics 
based on the free energy functional 
$$
F_\beta = {\rm Ent} + \beta E
$$
on the space of probability measures on $X$.
Here, ${\rm Ent}$ is the relative entropy (with respect to some reference measure), $E$ is the pluricomplex energy and $\beta \in \mathbb{R}$ (see \S \ref{functionals} in this paper).
If $K_X$ is ample, there exists a unique minimizer of the free energy functional $F_\beta$ for $\b>0$
and it is the measure defined by the (twisted) {\kah}--Einstein metric.
Berman also studied the negatively curved {\kah}--Einstein problem
from the viewpoint of statistical mechanics \cite{Ber17,Ber18}.
He proved that the law of the empirical measure satisfies a large deviation principle whose rate function is $F_\beta - \inf F_\beta$ and there exists a ``canonical" sequence of measures which converges to the measure defined by the {\kah}--Einstein metric.

When $-K_X$ is ample, i.e., $X$ is a Fano manifold, 
the existence of a K\"ahler–Einstein metric with positive Ricci curvature is subject to stability obstructions.
It is well known that the uniqueness and existence of a {\kah}--Einstein metric is equivalent to the (uniform) $K$-stability as a solution of (the uniform version of) the Yau--Tian--Donaldson conjecture \cite{BBJ,CDS1,CDS2,CDS3,Ti}.
On the other hand, Berman also studied the positively curved {\kah}--Einstein problem from the viewpoints of thermodynamics and statistical mechanics.
He introduced the microscopic stability thresholds and the notion of the (uniform) Gibbs stability by considering the log canonical threshold for the pair of the product $X^N$ and the divisor defined by the Slater determinant which is a ``canonical" holomorphic section in $H^0 (X^N , -k K_{X^N})$.
Berman gave the inequality between the free energy functional $F_{-\gamma}$ for $\gamma>0$
and the partition function $\mathcal{Z}_N (-\gamma)$,
and proved that a uniformly Gibbs stable Fano manifold admits a unique {\kah}--Einstein metric (\cite{BerQP}).
This is an analytic counterpart of the previous work by Fujita--Odaka \cite{FO}.

Hultgren-Witt Nystr\"om \cite{HW} introduced coupled K\"ahler–Einstein metrics
as a notion of canonical metrics on a compact K\"ahler manifold with definite canonical bundle, including manifolds that admit no K\"ahler–Einstein metric.
For a decomposition of $\lambda K_X,  \,\,(\lambda \neq 0)$, i.e., an $m$-tuple of ample $\mathbb{Q}$-line bundles $\{L_i \}_{i=1}^m$;
$$
\lambda K_X = L_1 + ... +L_m,
$$
we say that an $m$-tuple of {\kah} metrics $(\omega_i)_i \in \prod_{i=1}^m c_1 (L_i)$ is a coupled {\kah}--Einstein metric if
$$
{\rm Ric}\, \omega_1 =\cdots={\rm Ric }\, \omega_m = -\lambda^{-1} \sum_{i=1}^{m} \omega_i.
$$
Hultgren--Witt Nystr\"{o}m defined the coupled Ding functional (which is called the associated Ding functional in their paper \cite{HW}) on the product of the spaces of (singular) {\kah} metrics and proved the regularity of its critical point.
Moreover, in the Fano case, they defined the coupled $K$-polystability for a tuple of test configurations with the same base space and proved that it is a necessary condition for the existence of coupled {\kah}--Einstein metrics.
After their fundamental work, Fujita--Hashimoto \cite{FH} introduced the notion of the sum of test configurations and proved that the existence of coupled {\kah}--Einstein metrics is equivalent to the uniform coupled Ding stability.
Fujita--Hashimoto also proved that the uniform coupled Ding stability is equivalent to the condition that the corresponding stability threshold which is called the coupled delta-invariant \cite{RTZ} is greater than 1.

In this paper, we extend Berman's thermodynamical formalism \cite{Ber13} to the coupled setting by introducing {\it the coupled free energy}:
$$
F^{\rm cp}_{\beta} := {\rm Ent} + \beta E^{\rm cp}.
$$
Here, the symbol $E^{\rm cp}$ denotes {\it the coupled pluricomplex energy} (see \S \ref{functionals} in this paper).
Our arguments build on
\cite{Ber13,Ber18,BBJ,BHJ,BJKst,BJnaCY,Chi Li,Tru}.
The coupled free energy $F^{\rm cp}_\b$ provides a common framework for the following results.

\begin{itemize}
\item We generalize the thermodynamical formalism of Berman to our coupled free energy. As a consequence, we prove that it suffices to check the coercivity of the coupled Ding functional $D^{\rm cp}_{-\gamma}$ on the ``diagonal" in order to prove the coercivity of $D^{\rm cp}_{-\gamma}$ on the product of the finite energy spaces in the Fano case.
\item When $K_X$ is ample, we prove that the law of the empirical measures satisfies a large deviation principle with rate function $F_\beta^{\rm cp} - \inf F_\beta^{\rm cp}$ under the assumption that a conjecture of Hultgren holds.
As a consequence, we construct a canonically determined measure converging to
the Monge--Amp\`ere measure of the unique (twisted) coupled K\"ahler--Einstein metric.
\item In the Fano case, we introduce the non-Archimedean coupled free energy $F_{-\gamma}^{\rm cp,NA}$ on the space of finite energy measures on the Berkovich analytification $X^{\rm NA}$.
This is a natural generalization of the non-Archimedean free energy introduced by Boucksom--Jonsson \cite[\S 2.5]{BJKst} and also provides a diagonal approach to the non-Archimedean coupled Ding functional $D_{-\gamma}^{\rm cp,NA}$.
Extending the stability condition of Boucksom--Jonsson \cite[\S 2.8-10]{BJKst}
to the coupled setting,
we introduce thermodynamical $K$-stability and prove that it is equivalent to the existence of a unique coupled K\"ahler–Einstein metric.
This result recovers the result of Fujita--Hashimoto \cite{FH} when the holomorphic automorphism group is finite.
Moreover, we show that it suffices to test the non-negativity of the non-Archimedean coupled free energy of a measure defined by a test configuration of $(X,-K_X)$ in order to characterize the existence of a coupled {\kah}--Einstein metric.
\item In the Fano case, we define analogs of the partition function, the microscopic stability threshold and the uniform Gibbs stability
in the coupled setting under the assumption that $L_1^n =\cdots=L_m^n$.
We give an inequality between the coupled partition function and $F_{-\gamma}^{\rm cp}$
which generalizes the result of Berman \cite[Theorem 1.2]{BerQP}.
As a consequence, we obtain the coupled version of the results by Fujita--Odaka \cite[Corollary 2.6]{FO} and Berman \cite[Corollary 1.3]{BerQP}.
\end{itemize}

\subsection{Canonically polarized case}

We assume that $X$ is a projective manifold such that the canonical line bundle is ample and there exists an $m$-tuple of ample $\mathbb{Q}$-line bundles $L_i$ such that
$$
K_X = L_1 +\cdots+ L_m .
$$
In this case, Hultgren--Witt Nystr\"{o}m proved the existence and uniqueness of the coupled {\kah}--Einstein metric (\cite{HW} and see also Pingali's work \cite{Pi}).
Their proof uses the variational characterization of coupled K\"ahler–Einstein metrics as critical points of the coupled Ding functional. 
In the coupled analogue of Berman's framework \cite{Ber13},
we prove that $F_\beta^{\rm cp}$ has a unique minimizer for every $\b>0$
(Proposition \ref{uniqueness and existence}).
In addition, the minimizer corresponds to the measure defined by the twisted coupled {\kah}--Einstein metric.

A natural question is whether Berman's work on a large deviation principle for 
the K\"ahler--Einstein metric with negative Ricci curvature
extends to the coupled setting.
Berman's proof uses a Fekete configuration which maximizes the norm of the Slater determinant $\det S^{(k)} $ for an ample line bundle.
He also uses the convergence of the logarithm of the $L^\infty$-norm of $\det S^{(k)}$ to the Aubin--Mabuchi energy in some appropriate sense (\cite{BB10,BBWN}).
In our case, the existence of a Fekete configuration is a nontrivial problem because we consider a tuple of ample $\mathbb{Q}$-line bundles.
Hultgren proposed the notion of a {\it mutually asymptotic Fekete sequence} and conjectured that there exists such a sequence if and only if weight functions have the same equilibrium measure \cite[Conjecture 1]{Hul}.
Under the assumption that this conjecture holds, we can show that our coupled free energy naturally arises as a rate function in a large deviation principle for the twisted coupled {\kah}--Einstein metric.

We assume the equality of volumes
$$
L_1^n=\cdots=L_m^n.
$$
For each $i=1,2,\dots, m$, we define a holomorphic section (which is called the Slater determinant) on the product $X^{N_i}$ by
$$
{\rm det} S_i^{(k)} (x_1, ... , x_{N_i} ) :=  \det (s^i_\alpha (x_\beta))_{\alpha \beta} \in H^{0} (X^{N_i}, \otimes_j p_j^* L_i^k),
$$
where $(s^i_\alpha)_\alpha$
is an orthogonal basis of $H^0 (X, L_i^k)$ chosen as in
Section \ref{LDP-Gibbs},
and $N_i := {\rm dim} H^0 (X, L_i^k)$.
By fixing a Hermitian metric $H(i,k,0)$ on $L_i^k$, we define a Gibbs measure $\mu^{(N)}_\b$ on the product $X^N$ by
$$
\mu_\b^{(N)}:=\frac{1}{Z_{N,\b}}
\prod_i  \left| \pi_{i}^* \det S_i^{(k)} \right|^{2\b /k}_{H(i,k,0 ) }dV^{\otimes N},
$$
where $N := \max_i N_i$ and $Z_{N,\b}$ is the normalizing constant, called the partition function, chosen so that $\mu_\b^{(N)}$ is a (non-symmetric) probability measure. 
By a slight abuse of terminology, this defines the random point process $(X^N,\mu^{(N)}_\b)$ (see Remark \ref{RPP}).
Let $\delta_N:=\frac{1}{N}\sum_{j=1}^N\delta_{x_j}$ be the
empirical measure which is a map from $X^N$ to the space of 
probability measures on $X$,
where $\delta_x$ is the Dirac measure for $x\in X$.

\begin{thm}[Theorem \ref{cp-LDP}]
Assume that Conjecture 1 in \cite{Hul} (see Conjecture \ref{Hul-conj} in this paper) is true.
Fix $\beta >0$.
Then, the law of the empirical measure
$\Gamma_{N,\b}:=(\delta_N)_* \mu^{(N)}_\beta $ satisfies 
a large deviation principle with speed $N$ and 
a rate function $F_\beta^{\rm cp} - \inf_{\mathcal{P}(X)} F_\beta^{\rm cp}.$
In particular, as $N\to\infty$, the law
$\Gamma_{N,\b}$ converges weakly to the Dirac measure
$\delta_{\mu_{\rm tcKE}}$,
where $\mu_{\rm tcKE}$ is the Monge–Amp\`ere measure of the unique twisted coupled K\"ahler–Einstein metric.
\end{thm}

As a consequence of the large deviation principle,
we construct a canonical sequence of measures
converging weakly to
$\mu_{\rm tcKE}$.
Let us consider the averaged $1$-point correlation measure
$$\rho_{\b,1}^{(N)}:=\frac{1}{N}\sum_{i=1}^N(\mathrm{pr_i})_*
\mu^{(N)}_\b,$$ where $\mathrm{pr}_i:X^N\to X$ is the projection
onto the $i$-th factor.
We emphasize that $\rho_{\b,1}^{(N)}$ depends only on
the decomposition ${L_i}$ of $K_X$,
the Hermitian metric $H(i,k,0)$ on $L_i^k$, and the volume form $dV$ on $X$,
which constitute the data defining a twisted coupled K\"ahler--Einstein metric.

\begin{cor}[Corollary \ref{r-correlation}]
If $\Gamma_{N,\b}\to\delta_{\mu_{\rm tcKE}}$,
then the averaged $1$-point correlation measure
$\rho_{\b,1}^{(N)}$ converges weakly to $\mu_{\rm tcKE}$.
\end{cor}

\subsection{Fano case}

\subsubsection{Coupled version of the thermodynamical formalism by Berman}

Assume that $X$ is a Fano manifold, i.e., the anticanonical line bundle $-K_X$ is ample, and $-K_X$ is decomposed into the sum of ample $\mathbb{Q}$-line bundles:
$$
-K_X = L_1 +\cdots+ L_m.
$$
In this case, the uniqueness of coupled {\kah}--Einstein metrics up to ${\rm Aut}_0(X)$, i.e., the identity component of the holomorphic automorphism group, is proved by Hultgren--Witt Nystr\"{o}m \cite{HW}.
In addition, Fujita--Hashimoto proved that the existence of a coupled {\kah}--Einstein metric is equivalent to the uniform coupled Ding stability \cite{FH}.
(Note that 
their result also covers Fano manifolds with a non-discrete holomorphic automorphism group.)
One of the goals of this paper is to study coupled {\kah}--Einstein metrics from the viewpoint of Berman's thermodynamical formalism \cite{Ber13}.
By replacing $\beta$ with $- \gamma < 0$, the corresponding coupled free energy is given by
$$
F_{-\gamma}^{\rm cp} (\mu) := {\rm Ent }_{dV} (\mu) - \gamma E^{\rm cp} (\mu), \,\,\, \mu \in \mathcal{M}^1(X).
$$
Here, the symbol $\mathcal{M}^1(X) $ denotes the space of measures with finite energy.
This is a natural generalization of the free energy functional introduced by Berman \cite{Ber13}.
If there exists a minimizer of $F_{-\gamma}^{\rm cp}$, it is the measure defined by a twisted coupled {\kah}--Einstein metric.
It is known that there exists a unique coupled {\kah}--Einstein metric if and only if the coupled Ding functional is coercive on the product of the finite energy spaces $\prod_{i=1}^{m} \mathcal{E}^1 (L_i)$ (\cite{HW,FH}).

The coupled free energy allows us to reduce the coercivity of 
the coupled Ding functional
$D^{\mathrm{cp}}_{-\gamma}$ on $\prod_{i=1}^{m} \mathcal{E}^1 (L_i)$ to the coercivity on the subset where all normalized Monge–Amp\`ere measures coincide. We call this subset the diagonal and denote it by $\Delta$;
$$
\Delta := \Delta_{\{L_i\}_i } := \left\{  ( \phi_i)_i \in \prod_i \mathcal{E}^1 (L_i) \, \middle|  \,\, {\rm MA_{\omega_i} }(\phi_i)= {\rm MA_{\omega_j} }(\phi_j) , \,\forall i,j=1,...,m \right\}.
$$
\begin{thm}[Theorem \ref{coercivity}]
\label{diagonal result}
Fix $\gamma >0$.
The following are equivalent.

\medskip

{\rm (i)} $F_{-\gamma}^{\rm cp}$ is coercive on $\mathcal{M}^1 (X)$;

{\rm (ii)} there exists $\epsilon >0$ such that $F_{-\gamma-\epsilon}^{\rm cp}$ is bounded from below on $\mathcal{M}^1 (X)$;

{\rm (iii)} $D_{-\gamma}^{\rm cp}$ is coercive on $\prod_i \mathcal{E}^1 (L_i)$;

{\rm (iv)} there exists $\epsilon >0$ such that $D_{-\gamma-\epsilon}^{\rm cp}$ is bounded from below on $\prod_i \mathcal{E}^1 (L_i)$;

{\rm (v)} $D_{-\gamma}^{\rm cp}$ is coercive on $\Delta $;

{\rm (vi)} there exists $\epsilon >0$ such that $D_{-\gamma-\epsilon}^{\rm cp}$ is bounded from below on $\Delta $.

\end{thm}

\subsubsection{Coupled version of the thermodynamical formalism in the non-Archimedean setting by Boucksom--Jonsson}

As a non-Archimedean analog of the thermodynamical formalism by Berman \cite{Ber13}, Boucksom--Jonsson introduced the non-Archimedean free energy which is equivalent to the non-Archimedean Mabuchi K-energy by the non-Archimedean Calabi--Yau theorem \cite{BJnaCY}.
In the non-Archimedean setting, we can define the non-Archimedean coupled free energy as follows (see Definition \ref{def of NA cp free}):
$$
F_{-\gamma}^{\rm cp,NA} (\mu) := {\rm Ent }^{\rm NA} (\mu) - \gamma E^{\rm cp,NA} (\mu), \,\,\, \mu \in \mathcal{M}^{\rm 1,NA}(X).
$$

The following non-Archimedean analogue shows that the uniform lower bounds for the coupled Ding functional can likewise be tested on the non-Archimedean diagonal
$\Delta^{\mathrm{NA}}$ for $\{L_i \}_i$;
$$
\Delta^{\rm NA} := \Delta_{\{L_i\}_i }^{\rm NA} := \left\{  ( \phi_i)_i \in \prod_i \mathcal{E}^{1,\rm NA} (L_i) \, \middle|  \, \, \ {\rm MA}_i^{\rm NA} (\phi_i) =  {\rm MA}_j^{\rm NA}(\phi_j) , \,\forall i,j=1,...,m \right\}.
$$

\begin{thm}[Theorem \ref{NAcoercivity}]
\label{NA coercivity}
Fix $\gamma >0$.
The following are equivalent.

\medskip

{\rm (i)} there exists $\epsilon >0$ such that $F_{-\gamma}^{\rm cp,NA} \geq \epsilon E^{\rm cp,NA}$ on $\mathcal{M}^{\rm 1,NA} (X)$;

{\rm (ii)} there exists $\epsilon >0$ such that $F_{-\gamma-\epsilon}^{\rm cp,NA} \geq 0$ on $\mathcal{M}^{\rm 1,NA} (X)$;

{\rm (iii)} there exists $\epsilon >0$ such that $D_{-\gamma}^{\rm cp,NA} \geq \epsilon J^{\rm cp,NA}$ on $\prod_i \mathcal{E}^{\rm 1,NA} (L_i)$;

{\rm (iv)} there exists $\epsilon >0$ such that $D_{-\gamma-\epsilon}^{\rm cp,NA} \geq 0$ on $\prod_i \mathcal{E}^{\rm 1,NA} (L_i)$;

{\rm (v)} there exists $\epsilon >0$ such that $D_{-\gamma}^{\rm cp,NA} \geq \epsilon J^{\rm cp,NA}$ on $\Delta^{\rm NA}$;

{\rm (vi)} there exists $\epsilon >0$ such that $D_{-\gamma-\epsilon}^{\rm cp,NA} \geq 0$ on $\Delta^{\rm NA}$.

\end{thm}
As a coupled version of the stability condition introduced by Boucksom--Jonsson \cite[\S 2.8-2.10]{BJKst}, we introduce the following.

\begin{definition}[Definition \ref{TK}]
We say that $(X, \{ L_i \}_i)$ is thermodynamically $K$-stable if there exists $\epsilon >0$ such that  $F^{\rm cp, NA}_{-1-\epsilon} \geq 0$ on $\mathcal{M}^{\rm 1,NA}(X)$, equivalently, 
\begin{equation}
    \label{thermo}
F^{\rm cp, NA}_{-1} (\mu) \geq \epsilon E^{\rm cp,NA} (\mu), \,\,\,\, \forall \mu \in \mathcal{M}^{\rm 1,NA}(X).
\end{equation}
\end{definition}
Note that the notion of the coupled $K$-(poly)stability was introduced by Hultgren--Witt Nystr\"{o}m \cite{HW}.
As a consequence of the slope formulae by Berman--Boucksom--Jonsson \cite{BBJ} and Chi Li \cite{Chi Li} for maximal geodesic rays, we can characterize the existence of a coupled {\kah}--Einstein metric by the thermodynamical $K$-stability.
Moreover, by the regularization process of the non-Archimedean entropy ${\rm Ent}^{\rm NA}$ with respect to the strong topology established by Boucksom--Jonsson \cite{BJKst}, Chi Li \cite{Chi Li} and Trusiani \cite{Tru}, we obtain the following result which contains a modified version of the conjecture of Hultgren--Witt Nystr\"{o}m \cite{HW} when the holomorphic automorphism group of $X$ is finite.
\begin{thm}[Theorem \ref{prob} and Corollaries \ref{analytic correspondence}, \ref{alg}]
\label{Fano}
Assume that ${\rm Aut} (X)$ is finite.
The following are equivalent.
\begin{enumerate}

\item There exists a unique coupled {\kah}--Einstein metric;
\item $F_{-1}^{\rm cp}$ is coercive on $\mathcal{M}^1 (X)$;
\item $F^{\rm cp}_{-1}$ is lower semicontinuous on the space of probability measures $\mathcal{P}(X)$ with respect to the weak topology;
\item $(X,\{L_i\}_i)$ is thermodynamically K-stable;
\item $\exists \epsilon > 0 $ s.t. $F^{\rm cp, NA}_{-1-\epsilon}  (\mu)\geq 0$ for any $\mu = {\rm MA}^{\rm NA} (\phi_{(\mathcal{X}, \mathcal
L)})$ with $\phi_{(\mathcal{X}, \mathcal
L)} \in \mathcal{H}^{\rm NA} (-K_X)$, i.e., $\mu$ is a measure associated with a test configuration for $(X,-K_X)$. 
\end{enumerate}
\end{thm}

Nonnegativity of $F^{\mathrm{cp,NA}}_{-1}$ also yields a coupled analog of Fujita's upper bound for the anticanonical volume \cite{Fj}.
\begin{thm}[Theorem \ref{cp Fj 1} and Corollary \ref{cp Fj 2}]
   Assume that $F^{\rm cp ,NA}_{-1} \geq 0$ on $\mathcal{M}^{1, \rm NA}(X)$.    Then, we have
    \begin{equation}
    \label{sum of volume}
            (L_1^n)^{\frac{1}{n}} + ... + (L_m^n)^{\frac{1}{n}} \leq n+1.
    \end{equation}
Moreover, 
if $F^{\rm cp ,NA}_{-1} \geq 0$ on $\mathcal{M}^{1, \rm NA}(X)$ and $(-K_X)^n \leq (n+1)^n$, then the equality (\ref{sum of volume}) holds iff $X \simeq \mathbb{P}^n$ and $L_i = \mathcal{O}(k_i)$ for some $k_i \in \mathbb{Q}_{>0}$.
\end{thm}
Note that the projective space $\mathbb{P}^n$ achieves the equality in (\ref{sum of volume}) for any decomposition.
The proof follows Fujita's argument \cite{Fj}.
We can also show the inequality (\ref{sum of volume}) by assuming that $(X,\{L_i\}_i)$ is coupled Ding semistable.
The characterization of the equality case also uses
the Brunn--Minkowski inequality (see \cite[Corollary 1.6.3]{La2} and \cite[Theorem 2.6]{LX}).
\begin{rem}
Note that any Fano manifold $X$ satisfies $(-K_X)^n \leq (n+1)^n$ if $n={\rm dim} \,  X \leq 3$ \cite{Isk,MM}.
Thus, the equality (\ref{sum of volume}) holds iff $X \simeq \mathbb{P}^n$ for any coupled {\kah}--Einstein Fano manifold with ${\rm dim} \,  X \leq 3$.
\end{rem}
\begin{rem}
Fujita proved that the inequality $(-K_X)^n \leq (n+1)^n$ holds if $X$ is Ding semi-stable \cite{Fj}.
On the other hand, there exists a Fano manifold which is not Ding semi-stable but admits a 2-coupled {\kah}--Einstein metric by Yotsutani \cite{Yo}.
\end{rem}

\subsubsection{Microscopic stability threshold and Gibbs stability in the coupled setting}

Berman also studies {\kah}--Einstein metrics on Fano manifolds from the viewpoint of statistical mechanics \cite{Ber17,Ber18}.
He introduced the microscopic stability thresholds and (uniform) Gibbs stability by considering a certain integrability condition of the Slater determinant for the anticanonical line bundle.
Fujita--Odaka \cite{FO} proved that a uniformly Gibbs stable Fano manifold is uniformly $K$-stable by establishing the inequality between the microscopic stability threshold and the delta-invariant introduced by themselves (see also Fujita's result \cite{Fj16}).
As an analytic counterpart of the work of Fujita--Odaka, Berman \cite{BerQP} proved that a uniformly Gibbs stable Fano manifold has a unique {\kah}--Einstein metric. 
His proof is based on the fact that the logarithm of the partition function gives a lower bound of the free energy.

As an application of the coupled free energy $F_{-\gamma}^{\rm cp}$, we can prove analogs of the works of Fujita--Odaka and Berman \cite{FO,BerQP} under the assumption that ample line bundles $L_i$ have the same volume, i.e.,
$$
V_1 =\cdots = V_m.
$$
Here, we write $V_i = L_i^n$.
For each $i$, we define a holomorphic section on the product $X^{N_i}$ by
$$
{\rm det} \, S_i^{(k)} (x_1, ... , x_{N_i} ) :=  \det (s^i_\alpha (x_\beta))_{\alpha \beta} \in H^{0} (X^{N_i}, \otimes_j p_j^* L_i^k),
$$
where $(s^i_\alpha)_\alpha $ denotes a basis of $H^0 (X, L_i^k)$ and $N_i := {\rm dim} H^0 (X, L_i^k)$.
By relabeling $i$, we assume that $N_1 \leq N_2 \leq ... \leq N_m$.
For $\boldsymbol{\gamma} = (\gamma_1 ,..., \gamma_m) \in \mathbb{R}^m$, we define the coupled partition function by
$$
\mathcal{Z}_{N}^{\rm cp} (- \boldsymbol{\gamma}) := \int_{X^N} \prod_i  \left| \pi_{i}^* \det S_i^{(k)} \right|^{-2\gamma_i /k}_{H(i,k,0 ) }  dV^{\otimes N}  \in \mathbb{R}_{>0} \cup\{+\infty\},
$$
where $N := \max_i N_i$ and $H(i,k,0)$ is a fixed Hermitian metric on $L_i^k$.
We introduce the notion of microscopic stability thresholds and the uniform Gibbs stability in the coupled setting.
\begin{definition}
The coupled microscopic stability threshold at level $k$ is defined by
$$
\gamma_{k}^{\rm cp}\left( X, \{ L_i \}_{i=1}^{m} \right) := \sup \left\{ \gamma \geq 0 \, \middle| \,  \mathcal{Z}_{N}^{\rm cp}(-(\gamma,..., \gamma)) < + \infty \right\} .
$$
We say that $(X, \{L_i\}_i)$ is uniformly coupled Gibbs stable if $$\gamma^{\rm cp} \left( X, \{ L_i \}_{i=1}^{m} \right) := \liminf_k \gamma_k^{\rm cp} \left( X, \{ L_i \}_{i=1}^{m} \right) >1.$$
\end{definition}
Note that $\gamma^{\rm cp} (X, \{L_i \}_i)$ and $\gamma_k^{\rm cp} (X, \{L_i \}_i)$ do not depend on the choice of reference data, i.e., the choice of bases and smooth Hermitian metrics.
As a coupled version of \cite{BerQP}, we prove that the logarithm of the coupled partition function gives a lower bound of the coupled free energy in the following sense.

\begin{thm}[Theorem \ref{partition and free energy}]
\label{partition and free energy}
Fix $\gamma_0 > 0$.
We set $\boldsymbol{\gamma} = (\gamma_1,...,\gamma_m)$, where $$
\gamma_i := \frac{\gamma_0}{1-\epsilon_i} \frac{N}{N_1}
$$
for sufficiently small $\epsilon_i>0 ,\,\,i=1 ,2, ....,m$.
Then, there exists a positive constant $C_0 >  0$ such that for any sufficiently integer $k \in \mathbb{Z}_{>0}$ and small $\tau >0$, we have
\begin{eqnarray*}
\hspace{-20pt}-\frac{1}{ N  (1+\tau)} \log \mathcal{Z}^{\rm cp}_{N} (-(1+\tau)\boldsymbol{\gamma})\leq  F_{-(1 -k^{-1}C_0)\gamma_0}^{\rm cp}\left( \mu \right)  +O \left( \frac{\log k}{k} \right) , \,\,\, \forall \mu \in \mathcal{M}_1(X).
\end{eqnarray*}
\end{thm}

As a consequence, we obtain the coupled version of the results by Fujita--Odaka \cite{FO} and Berman \cite{BerQP} by using Kewei Zhang's work for coupled delta-invariant \cite{KZ}.
\begin{cor}
\label{ext of BFO}
The inequality $ \delta^{\rm cp} \left( X, \{ L_i \}_{i=1}^{m} \right) \geq \gamma^{\rm cp} \left( X, \{ L_i \}_{i=1}^{m} \right)$ holds.
In particular, if $(X , \{L_i \}_i)$ is uniformly coupled Gibbs stable, then there exists a unique coupled {\kah}--Einstein metric.
\end{cor}
The following diagram summarizes the implications proved here and in earlier work in the Fano case.

\begin{flalign*}
\xymatrix{
\mbox{$F_{-1-\epsilon}^{\rm cp,NA}\circ {\rm MA}^{\rm NA}_{-K_X} \geq 0$ on $\mathcal{H}^{\rm NA}(-K_X)$}  &\mbox{$D_{-1}^{\rm cp,NA} \geq \epsilon J^{\rm cp,NA}$ on $\Delta^{\rm NA}$} \ar@{->}[dl] \ar@{<-}[r]^{\rm }  &\mbox{$D_{-1}^{\rm cp,NA} \geq \epsilon J^{\rm cp,NA}$}  \ar@{<->}[d]^{\rm \cite{FH}} \\
\mbox{$F_{-1}^{\rm cp,NA} \geq \epsilon E^{\rm cp,NA}$ on $\mathcal{M}^{\rm 1,NA}(X)$}  \ar@{<->}[u]^{\rm Thm \, \footnotesize{\ref{one para}}} \ar@{->}[urr]^{\rm Thm \, \footnotesize{\ref{NAcoercivity}}}  \ar@{->}[r]^{\rm } \ar@{<-}[rrd]^{\rm see \, proof \,  of \, Thm \, \footnotesize{\ref{TK and coer}}}  &   \mbox{\mbox{$\delta^{\rm cp} > 1$} }\ar@{->}[dr]^{\rm \cite{KZ}} \ar@{<->}[r]^{\rm \cite{FH}} & \mbox{unif cDing stable } \ar@{<->}[d]^{\rm \cite{FH}} \\
\mbox{$\gamma^{\rm cp} > 1$ and $V_1 = ... = V_m$}  & \mbox{$D_{-1}^{\rm cp}$ is coercive on $\Delta$}  \ar@{<-}[r]^{\rm  } & \mbox{$D_{-1}^{\rm cp}$ is coercive} \ar@{<->}[d]^{\rm \cite{HW,Tak,FH}} \\
\mbox{$F_{-1}^{\rm cp}$ is coercive on $\mathcal{M}^{\rm 1}(X)$} \ar@{<-}[ru] \ar@{->}[rru]^{\rm Thm \, \footnotesize{\ref{coercivity}}}\ar@{<-}[u]^{\rm Thm \, \footnotesize{\ref{partition and free energy}}}& \mbox{  $F_{-1}^{\rm cp} \, {\rm is \, lsc}$ on $\mathcal{P}(X)$}  \ar@{<-}[l]^{\rm Thm \, \footnotesize{\ref{prob}}} &\mbox{$\exists$cKE metric}  \ar@{<-}[l]^{\rm Thm \, \footnotesize{\ref{prob}} }  \\
}
\end{flalign*}

This paper is organized as follows.
Section 2 recalls coupled K\"ahler–Einstein metrics and introduces the coupled free energy. Section 3 studies the canonically polarized case and proves a large deviation principle under Hultgren's conjecture. Section 4 treats the Fano case: it establishes the diagonal criterion for the coercivity, introduces thermodynamical $K$-stability, proves a coupled volume bound, and relates uniform coupled Gibbs stability to the existence of a unique coupled K\"ahler–Einstein metric.

\begin{ack}
The authors would like to thank
Rolf Andreasson, Yoshinori Hashimoto, Eiji Inoue, Yuji Odaka,
and Shunsuke Saito
for helpful discussions and comments.
This collaboration began during the workshop
``Probabilistic methods in complex geometry" held in Como, Italy, in June 2025,
and the authors would like to thank the organizers for providing this opportunity.
The first author is supported by JSPS Grant-in-Aid for Early-Career Scientists No. 25K17260.
The second author is supported by JSPS Grant-in-Aid for Early-Career Scientists
No. 24K16917.
\end{ack}

\begin{ai}
ChatGPT-6.1 Sol was used to explore the proof of Proposition \ref{submean}.
It was also used for editorial assistance in preparing the manuscript.
The authors take full responsibility for 
all mathematical arguments and statements in this paper.
\end{ai}

\section{Coupled {\kah}--Einstein metrics and functionals}

In this section, we introduce the coupled free energy which is a natural generalization of the free energy functional in the thermodynamical formalism of {\kah}--Einstein metrics by Berman \cite{Ber13}.
We also recall the definition of coupled {\kah}--Einstein metrics introduced by Hultgren--Witt Nystr\"{o}m.

\subsection{Coupled {\kah}--Einstein metric}

Throughout this paper, we consider an $n$-dimensional projective manifold $X$ such that the canonical line bundle $K_X$ is definite, i.e., $K_X > 0$ or $K_X <0$.
We assume that there exist ample $\mathbb{Q}$-line bundles $L_i$, $1 \leq i \leq m$ such that
$$
\pm K_X = L_1 + \cdot \cdot \cdot + L_m.
$$

Since all line bundles $L_i$ are ample, we can find smooth Hermitian metrics $h_i$ on $L_i$ such that their curvature forms $\Theta(L_i, h_i)$ multiplied by $\sqrt{-1}$ are {\kah} forms $\omega_i \in c_1 (L_i)$.

Hultgren--Witt Nystr\"{o}m introduced the following as a natural generalization of the notion of {\kah}--Einstein metrics.

\begin{definition}[\cite{HW}]

An $m$-tuple of {\kah} metrics $(\omega_i)_{i=1}^m$ is said to be a coupled {\kah}--Einstein metric if the equalities
$$
 \mp  {\rm Ric}\, \omega_{j} = \sum_{i=1}^{m}\omega_{i} 
$$
hold for any $j=1, ... ,m$.

\end{definition}

\begin{exam}

\begin{enumerate}
\item Any {\kah}--Einstein metric $\omega_{\rm KE}$ on $X$ with definite canonical line bundle $K_X$ is a coupled {\kah}--Einstein metric for $m=1$. More generally, the pair $(X, \{ \lambda k_i K_X \}_{i=1}^{m})$ for $k_i \in \mathbb{Q}_{>0}$ with $k_1 +... + k_m =1$ admits a coupled {\kah}--Einstein metric $ (k_i \omega_{\rm KE})_i$ since Ricci forms are scale invariant.
\item Assume that Fano manifolds $X$ and $Y$ admit {\kah}--Einstein metrics $\omega_X$ and $\omega_Y$ respectively. 
$Z : = X \times Y$. We have $-K_Z = -p_X^* K_X - p_Y^* K_Y$ for the projections $p_X :Z \to X, \, p_Y :Z \to Y$. For $a, b \in (0,1)$, 
we set K\"aler classes
$L_1 := -ap_X^* K_X - bp_Y^* K_Y$, $L_2 := -(1-a)p_X^* K_X - (1-b)p_Y^* K_Y$. Then, the pair $(Z, \{L_1, L_2\})$ has a coupled {\kah}--Einstein metric $(a p_X^* \omega_X+ bp_Y^* \omega_Y, (1-a) p_X^* \omega_X+ (1-b)p_Y^* \omega_Y)$.
\end{enumerate}

\end{exam}

By the $\partial\bar\partial$-lemma, every K\"ahler metric in $c_1(L_j)$
can be written as
$\omega_{j, \varphi_j} := \omega_j + \dol \varphi_j$.
The tuple $(\omega_{j,\varphi_j})_j$ is coupled K\"ahler--Einstein
if and only if its potentials satisfy the following systems of
complex Monge--Amp\`ere equations:
$$
\displaystyle \frac{\omega_{j, \varphi_j}^n}{V_j}  = \dfrac{e^{ \pm \sum_i \varphi_i} dV  }{\int_X  e^{ \pm  \sum_i \varphi_i} dV},\,\,\,\,  j = 1 , ..., m.
$$
Here, $\sqrt{-1} \Theta (K_X ,dV^*) = \pm \sum_i \omega_i$ and $V_i = L_i^n$.

\begin{rem}
In the paper of Hultgren--Witt Nystr\"{o}m \cite{HW}, they introduced the notion of coupled {\kah}--Einstein metrics for the decomposition by {\kah} classes:
$$
\mp c_1 (X) = \alpha_1 + ... + \alpha_m. 
$$
Several of our results, including the variational properties of the coupled free energy, hold for a decomposition by {\kah} classes.
However, the definition of the Slater determinant and the partition function
require a decomposition by ample $\mathbb{Q}$-line bundles.
\end{rem}

\subsection{Functionals on the space of probability measures}
\label{functionals}

In this subsection, we define functionals associated with the coupled {\kah}--Einstein metric.

Recall that we fix a {\kah} metric $\omega_i$ which is the curvature form (multiplied by $\sqrt{-1}$) for a Hermitian metric $h_i$ on $L_i$ with $dV = \otimes h_i$. 
We write the spaces of $\omega_i$-plurisubharmonic functions as
$$
{\rm PSH} (X , L_i) := \left\{ u : X \to [ - \infty, \infty ) \,\, {\rm upper\,semicontinuous}  \,\,\middle| \,\, \omega_i + \dol u \geq 0 \,\, {\rm as \,\, a \,\,  current }\right\}
$$
for each $i$.
We also define the spaces of (smooth) {\kah} potentials by
$$
\mathcal{H}(L_i) := \left\{ \phi \in C^\infty(X, \mathbb{R})  \,\,\middle| \,\, \omega_i + \dol \phi > 0 \right\}
$$
for each $i$.
Note that we have $ \mathcal{H}(L_i)  \subset {\rm PSH} (X , L_i)$.

We write the space of Borel probability measures on $X$ as $\mathcal{P}(X)$.
The complex {\ma} operator on $\mathcal{H}(L_i)$ is defined by
$$
\mathcal{H}(L_i) \ni \phi_i \mapsto {\rm MA}_{\omega_i} (\phi_i) : = \frac{\omega_{i,\phi_i}^n}{V_i} \in \mathcal{P}(X) .
$$
By the Calabi--Yau theorem \cite{Yau}, 
a volume form $\mu \in \mathcal{P} (X)$ has a potential
$\phi \in \mathcal{H}(L_i)$, unique up to an additive constant,
such that
${\rm MA }_{\omega_i}(\phi) = \mu$.

Although wedge products of arbitrary currents are not generally defined, non-pluripolar products in the sense of Bedford-Taylor \cite{BD}
provide a Monge–Amp\`ere operator for $\omega_i$-plurisubharmonic potentials.
By considering the nonpluripolar product, we define the Aubin--Mabuchi energy by
$$
\mathcal{E}_i (\phi_i ) := \frac{1}{(n+1) V_i } \sum_{j=0}^{n} \int_X \phi_i \, \omega_{i , \phi_i}^j \wedge \omega_{i}^{n-j}  , \,\,\,\,\, \phi_i \in {\rm PSH}(X, L_i).
$$
By following \cite{Ber13,BBGZ}, we recall the definition of the finite energy space.

\begin{definition}
    The finite energy space $\mathcal{E}^1(X,L_i)$  is defined by
$$
\mathcal{E}^1 (L_i) := \left\{ 
u \in  {\rm PSH}(X, L_i) \, \middle| \, \mathcal{E}_i (u) > - \infty \,\, {\rm and} \, \int_X \omega_{i ,u}^n =V_i \right\} .
$$
\end{definition}

The following functionals play the role of a distance of the finite energy space $\mathcal{E}^1 (L_i)$.
\begin{definition}
For $\phi_i \in \mathcal{E}^1 (L_i)$, we set
$$
I_i(\phi_i) := \int_X\phi_i ({\rm MA}_{\omega_i}(0)- {\rm MA}_{\omega_i}(\phi_i)) , \,\,\, J_i(\phi_i) := \int_X \phi _i{\rm MA}_{\omega_i}(0)- \mathcal{E}_i(\phi_i). 
$$
\end{definition}

It is well-known that $I_i$ and $J_i$ are finite on $\mathcal{E}^1 (L_i)$ and the following inequalities hold (see \cite[Section 2]{BBGZ}):
\begin{equation}
\label{IandJ}
\frac{1}{n+1} I_i \leq I_i-J_i \leq \frac{n}{n+1} I_i .
\end{equation}
By using the Aubin--Mabuchi energy, we recall the following functional on the space of probability measures.
\begin{definition}
   We define the pluricomplex energy for $(X, L_i)$ by
    $$
    E_{i} (\mu) := \sup_{u \in {\rm  PSH}(X,L_i)}\left\{  \mathcal{E}_i (u) - \int_X u \mu \right\}, \,\,\, \mu \in \mathcal{P} (X) .
    $$
   The space of finite energy probability measures is defined by
$$
\mathcal{M}^1 ( X) := \left\{ \, \mu \in \mathcal{P} (X) \,\middle| \,  E_i (\mu) < \infty  \, \right\}.
$$
    
\end{definition}
\begin{rem}
\label{indep of class}
Di Nezza proved that the space of finite energy measures $\mathcal{M}^1 ( X)$ does not depend on the choice of a {\kah} class (\cite[Proposition 3.1]{DN}).
Moreover, the strong topology of $\mathcal{M}^1(X)$ in the sense of \cite{BBEGZ} is also independent of the choice of {\kah} class by Boucksom--Jonsson \cite{BJfin}.
\end{rem}

\begin{rem}\label{projection}
In the above definition, the pluricomplex energy $E_i$ is not a Legendre transform of the Aubin--Mabuchi energy $\mathcal{E}_i$.
But, as explained in \S 2.5 of \cite{Ber13}, $E_i$ can be written as a Legendre transform of the composition $\mathcal{E}_i \circ P_{\omega_i}$, where 
$P_{\omega_i} : C^0 (X) \to C^0 (X) \cap {\rm PSH} (X, \omega_i)$
is a projection defined by
$$ u \mapsto P_{\omega_i} (u):=\sup \{ v(x) \, | \,  v \in  {\rm PSH }(X,\omega_i) , \, v \leq u \}.$$
In particular, $E_i$ is convex and lower semicontinuous on $\mathcal{P}(X)$.    
\end{rem}

A probability measure of finite energy is characterized as follows.

\begin{thm}[\cite{BBGZ}]
\label{finite energy}
A probability measure $\mu$ has finite energy if and only if  
there exists 
$\phi_{i, \mu} \in \mathcal{E}^1 (L_i)$
unique
up to an additive constant such that $\mu ={\rm MA}_{\omega_i} (\phi_{i,\mu})$. Moreover, we have $E_i (\mu) = \mathcal{E}_i (\phi_{i,\mu}) - \int_X \phi_{i, \mu} \mu  = (I_i-J_i) (\phi_{i,\mu})$.

\end{thm}

Note that the pluricomplex energy and the functional $J$ are compatible by (\ref{IandJ}):
$$
\frac{1}{n} J_i (\phi) \leq E_i ({\rm MA}_{\omega_i}(\phi)) \leq n J_i (\phi).
$$

We define the following functional on the space of probability measures $\mathcal{P}(X)$.

\begin{definition}[coupled pluricomplex energy]
\label{coupled pc energy}
We define the coupled pluricomplex energy by
$$
E^{\rm cp} (\mu) := \sum_{i=1}^m E_i (\mu), \,\,\,\, \mu \in \mathcal{P} (X).
$$
We also write the space of finite-energy probability measures by
$$
\mathcal{M}^1 ( X)=   \left\{ \, \mu \in \mathcal{P} (X) \,\middle| \,  E^{\rm cp} (\mu) < \infty  \, \right\}.
$$
\end{definition}
As we mentioned in Remark \ref{indep of class}, finiteness of pluricomplex energy is independent of the choice of {\kah} classes.
By definition and Theorem \ref{finite energy}, we immediately obtain the following.

\begin{prop}
\label{pluricomplex}
    A probability measure $\mu$ has finite energy if and only if we can find $\phi_{j, \mu} \in \mathcal{E}^1 (L_j)$ such that $\mu = {\rm MA}_{\omega_j} (\phi_{j,\mu})$ for any $j=1, ...,m$. Moreover, we have 
    $$
    E^{\rm cp} (\mu) = \sum_{i=1}^{m} \left(  \mathcal{E}_i (\phi_{i,\mu}) - \int_X \phi_{i,\mu} \mu \right).
    $$
\end{prop}

We recall the definition of the relative entropy.

\begin{definition}
The relative entropy with respect to the smooth probability measure $dV$ is defined by
$$
{\rm Ent}_{dV}(\mu) :=\sup_{\phi \in C^0 (X)} \left( - \log \int_X e^{\phi} dV  + \int_X \phi \mu \right), \,\,\,\, \mu \in \mathcal{P} (X).
$$
    
\end{definition}
Note that the relative entropy is the Legendre transform of the functional $C^0(X) \ni \phi \mapsto \log \int_X e^{\phi} \in \mathbb{R}$.
It is well-known that ${\rm Ent}_{dV} (\mu) = \int_X \log \left( \frac{d\mu}{dV} \right) \mu$ when $\mu$ is absolutely continuous with respect to $dV$, otherwise ${\rm Ent}_{dV} (\mu) = + \infty$.

We introduce the coupled free energy which generalizes the free energy functional in \cite{Ber13}.
 
\begin{definition}[coupled free energy]
\label{coupled free energy}
For $\beta \in \mathbb{R}$, we define the coupled free energy by
$$
F^{\rm cp}_{\beta} (\mu ) := {\rm Ent}_{dV} (\mu) + \beta E^{\rm cp} (\mu), \,\,\,\, \mu \in \mathcal{P} (X).
$$
\end{definition}
For a fixed reference measure $dV$,
the entropy term of $F_\beta^{\rm cp}$ is independent of the decomposition, 
whereas $E^{\mathrm{cp}}$ depends on the line bundles $\{L_i\}_{i=1}^m$. 
We study positive inverse temperature $\beta$ in the canonically polarized case 
and write $\beta=-\gamma<0$ in the Fano case. 
We show below that minimizers correspond to (twisted) coupled K\"ahler–Einstein metrics.

\section{Coupled free energy and Large Deviation Principle when $K_X >0$}

In this section, we assume that $K_X$ is ample.

\subsection{Relationship with the coupled Ding functional}

When $K_X > 0$, Hultgren-Witt Nystr\"{o}m proved the Aubin-Yau type existence and uniqueness result of coupled {\kah}--Einstein metrics.

\begin{thm}[\cite{HW}]
Assume that the canonical line bundle $K_X$ is ample and has a decomposition $K_X = L_1 + ... + L_m$ for some ample $\mathbb{Q}$-line bundles $L_i$.
Then, there exists a unique coupled {\kah}--Einstein metric.
\end{thm}

Hultgren-Witt Nystr\"om \cite{HW} proved this theorem using the coupled Ding functional. 
Pingali \cite{Pi} subsequently gave a proof by the continuity method.
We now introduce the corresponding mean field equation and establish the coupled analogue of \cite[Theorem 3.2]{Ber13}.

\begin{definition}
We define the normalized coupled {\ma} mean field equation associated with the data $( \{\omega_i\}_{i=1}^m,dV, \beta ) $ by
\begin{equation}
\label{cME}
\displaystyle \frac{\omega_{j, \varphi_j}^n}{V_j}  = \dfrac{e^{ \beta \sum_i \varphi_i} dV  }{\int_X  e^{ \beta  \sum_i \varphi_i} dV},\,\,\,\,  j = 1 , ..., m.
\end{equation}
The solution of the normalized coupled {\ma} mean field equation satisfies the twisted coupled {\kah}--Einstein equation:
$$
{\rm Ric}\, \omega_{j, \phi_j} = -\beta \sum_{i} \omega_{i ,\phi_i} + (1-\beta) {\rm Ric }\, dV ,\,\,\,\,  j = 1 , ..., m.
$$
We call an $m$-tuple $(\omega_{j,\phi_j})_j$ which satisfies the above equation a twisted coupled {\kah}--Einstein metric.
\end{definition}

Hultgren--Witt Nystr\"{o}m proved the regularity of the weak solution of (\ref{cME}) (\cite[Theorem 2.4]{HW} and see Theorem \ref{regularity} in this paper).

We recall the following functional on the product of $\mathcal{E}^1 (L_i)$ which is introduced by  Hultgren--Witt Nystr\"{o}m as a generalization of the Ding functional in the study of {\kah}--Einstein metrics.

\begin{definition}[coupled Ding functional]
Fix $\beta > 0$.
For $(\varphi_i)_{i=1}^{m} \in \prod_{i=1}^{m}\mathcal{E}^1 ( L_i)$, we define the coupled Ding functional by
$$
D^{\rm cp}_{\beta} ((\varphi_i )_i ) := -\frac{1}{\beta} \log \int_X e^{\beta \sum \varphi_i } dV   + \sum_i \mathcal{E}_i (\varphi_i).
$$
\end{definition}
Note that the coupled Ding functional is called the associated Ding functional in \cite{HW}.
We establish the coupled analogue of \cite[Theorem 3.2]{Ber13}.
\begin{prop}
\label{uniqueness and existence}
    For any $\beta >0 $, there exists a unique solution $( \phi_i)_{i=1}^m \in \prod_{i} \mathcal{E}^1(L_i)$ (up to additive constants) 
    of the equation (\ref{cME}).
Moreover,  $( \phi_i)_{i=1}^m$ is the unique maximizer of the coupled Ding functional $D^{\rm cp}_\beta$ and $\mu_{ME} := {\rm MA}_{\omega_i} (\phi_i)$ is the unique minimizer of the coupled free energy $F^{\rm cp}_{\beta}$.
\end{prop}

\begin{proof}
The existence and uniqueness of the solution follow from
the same argument as in
\cite[Theorem A]{HW} (see also \cite{Pi}).
The regularity of weak solutions of (\ref{cME}) by \cite{HW} implies that the solution  $( \phi_i)_{i=1}^m$ is a unique maximizer of $D^{\rm cp}_{\beta}$ \cite[\S 2.3]{HW}.

Since the free energy $F^{\rm cp}_\beta$ is strictly convex on $\{ {\rm Ent}_{dV} (\mu) < \infty \} \subset \mathcal{P}(X)$ (along affine segments), if a minimizer exists, it is unique.
So, it suffices to show that $\mu_{ME} := {\rm MA}_{\omega_i} (\phi_i)$ minimizes $F_\beta^{\rm cp}$.
For $\mu $ with $F_\beta^{\rm cp} (\mu) < \infty$, we set $\mu_t := (1-t)\mu_{ME} + t \mu$.
By \cite[Proposition 2.11, Proposition 2.7]{Ber13}, we have
$$
\frac{d {\rm Ent}_{dV}}{dt}{}_{t = 0+} = - \int_X \log \frac{\mu_{ME}}{dV} (\mu_{ME} - \mu) .
$$
and
$$
\frac{d E_i(\mu_t)}{dt} {}_{t=0+} = - \int_X  \phi_i (\mu-\mu_{ME}).
$$
for each $i$.
Since $(\phi_i)_i$ solves the equation (\ref{cME}),
the right derivative of $F_\beta^{\rm cp}(\mu_t)$ at $t=0$ vanishes.
Convexity therefore gives
$$
F_\beta^{\rm cp} (\mu) \geq F_\beta^{\rm cp} (\mu_{ME}) +     \frac{dF_\beta^{\rm cp} (\mu_t) }{ dt} {}_{t=0+} =  F_\beta^{\rm cp} (\mu_{ME}) 
$$ for any $\mu$ with $F_\beta^{\rm cp}(\mu) < \infty$.
Thus, $\mu_{ME}$ minimizes $ F_\beta^{\rm cp} $.
\end{proof}

By following \cite{Ber13}, we obtain the following inequality between $F_\beta^{\rm cp}$ and $D_\beta^{\rm cp}$.
\begin{prop}
For $\beta > 0$, we have
$$
\inf_{\mu \in \mathcal{M}^1(X)} F^{\rm cp}_{\beta} (\mu) \geq \sup_{(\varphi_i)_i \in \prod_{i} \mathcal{E}^1(X,\omega_i) } \beta D^{\rm cp}_{\beta}((\varphi_i)_i) .
$$
Moreover, the equality holds if and only if
$(\varphi_i)_i$ solves \eqref{cME} and $\mu=\operatorname{MA}_{\omega_i}(\varphi_i)$ for every $i$.
\end{prop}

\begin{proof}
From the definition of $ \mathcal{M}^1 ( X) \subset \mathcal{P} (X)$, it suffices to consider finite energy probability measures.
For any $\mu \in \mathcal{M}^1 ( X)$ and $(\varphi_i)_i \in \prod_{i} \mathcal{H}(L_i )$, we have
\begin{eqnarray*}
F^{\rm cp}_{\beta} (\mu )
&=& \sup_{\phi \in C^0 (X)} \left( - \log \int_X e^{\phi} dV  + \int_X \phi \mu \right) + \beta E^{\rm cp} (\mu) \\
&\geq & - \log \int_X e^{\beta \sum_i  \varphi_i } dV  +\beta \int_X \sum_i  \varphi_i\mu - \beta \int_X \sum_i  \varphi_i\mu + \beta \sum_i \mathcal{E}_i (\varphi_i) \\
&=& \beta D^{\rm cp}_{\beta}((\varphi_i)_i).
\end{eqnarray*}
So, we obtain $\inf_{\mu \in \mathcal{M}^1 ( X)} F^{\rm cp}_{\beta} (\mu) \geq \sup_{(\varphi_i)_i \in \prod_{i} \mathcal{H}_i } \beta D^{\rm cp}_{\beta}((\varphi_i)_i)$.
In addition, $D^{\rm cp}_{\beta}$ is 
upper
semicontinuous with respect to the $L^1$-topology by \cite[Lemma 3.1]{HW}, we have $\sup_{ \prod_{i} \mathcal{H}_i } \beta D^{\rm cp}_{\beta} = \sup_{ \prod_{i} \mathcal{E}^{1}(L_i) } \beta D^{\rm cp}_{\beta}$ and we have finished the proof of the first assertion.

Since $D_\beta^{\rm cp}$ and $F_\beta^{\rm cp}$ have the unique maximizer and minimizer respectively by Proposition \ref{uniqueness and existence}, the first argument in this proposition implies that if the equality $F^{\rm cp}_{\beta} (\mu ) = \beta  D^{\rm cp}_{\beta}((\varphi_i)_i) 
$ holds, then $(\varphi_i)_i$ solves (\ref{cME}).
Conversely, if $(\varphi_i)_i$ is a solution of  (\ref{cME}), by writing $\mu = \frac{ e^{\beta \sum_i \varphi_i} dV }{  \int_X e^{\beta \sum_i \varphi_i} dV }$, we have
\begin{eqnarray*}
F^{\rm cp}_{\beta} (\mu )
&=& \int_X \log \frac{\mu}{dV} \mu + \beta E^{\rm cp} (\mu) \\
&=& - \log  \int_X e^{\beta \sum_i \varphi_i} dV + \beta \int_X \left( \sum_i \varphi_i \right) \mu + \beta \left( \sum_i \mathcal{E}_i (\phi_i) - \int_X \left( \sum_i \varphi_i \right) \mu \right) \\
&=& \beta  D^{\rm cp}_{\beta}((\varphi_i)_i) .
\end{eqnarray*}
Here, we have used Proposition \ref{pluricomplex} in the second equality.
\end{proof}

\subsection{Large deviation principle and coupled K\"ahler--Einstein metrics}\label{LDP-cpKE}
Berman \cite{Ber17, Ber18} characterized the K\"ahler--Einstein measure 
on a canonically polarized manifold as the many-particle limit of canonical random point processes whose empirical measures satisfy an LDP. 
We extend this framework to the coupled setting, 
subject to Hultgren's conjecture (Conjecture \ref{Hul-conj}).
We emphasize that our coupled free energy 
$F^{\mathrm{cp}}_\b = {\rm Ent}_{dV}+\b E^{\mathrm{cp}}$
for $\b>0$
naturally appears as the rate function in the framework.

The large deviation principle (LDP for short) provides a natural framework
for describing the macroscopic behavior of a system consisting of a large number of particles. 
In our setting, a configuration of $N$ points on $X$ will be regarded as an $N$-particle system,
and its macroscopic state is described by the corresponding empirical measure.
The LDP quantifies the concentration of the law of the empirical measure
around its limiting measure and, more precisely,
describes the exponential decay of the probability of deviations
from the limiting measure in terms of a rate function.
We recall the precise definition below.

\begin{definition}
    Let $\mathcal{Z}$ be a complete separable metric space.
    \begin{enumerate}
        \item A function $I:\mathcal{Z}\to[0,\infty]$ is called a rate function
        if it is a lower semicontinuous. It is a good rate function if 
        it is also proper, i.e. every sublevel set is compact.
        \item We say that a sequence $\{\Gamma_k\}_{k=1}^\infty$
        of probability measures on $\mathcal{Z}$ satisfies a large deviation principle
        with speed $0<r_k\to\infty$ and rate function $I$
        if the inequality
        \begin{equation}\label{upbnd}
            \limsup_{k\to\infty}\frac{1}{r_k}\log\Gamma_k(\mathcal{F})
            \leq -\inf_{\mu\in\mathcal{F}}I(\mu)
        \end{equation}
        holds for any closed subset $\mathcal{F}\subset\mathcal{Z}$,
        and the inequality
        \begin{equation}\label{lowbdn}
            \liminf_{k\to\infty}\frac{1}{r_k}\log\Gamma_k(\mathcal{G})
            \geq -\inf_{\mu\in\mathcal{G}}I(\mu)
        \end{equation}
        holds for any open subset $\mathcal{G}\subset\mathcal{Z}$.
        \item The LDP is said to be weak if \eqref{upbnd} holds
        for any compact subset $\mathcal{F}$.
    \end{enumerate}
\end{definition}

In our setting, $\mathcal{Z}$ will be
the space of Borel probability measures $\mathcal{P}(X)$,
for a compact manifold $X$,
which is compact with respect to the weak convergence topology
by the Prokhorov theorem.
If $\mathcal{Z}$ is compact, then the notions of an LDP and a weak LDP coincide.
Thus, the following proposition is useful for our purpose
(see, for instance, \cite[Theorem 4.1.11 and 4.1.18]{DZbook}).

\begin{prop}\label{def-LDP'}
    Let $B_\e(\mu)$ be the ball of radius $\e>0$ centered at $\mu\in\mathcal{Z}$.
    Then a sequence $\{\Gamma_{N}\}_{N=1}^\infty$
    of probability measures on $\mathcal{Z}$
    satisfies a weak LDP with speed $r_N$ and a rate function $I$ if and only if
    the equality
    \begin{equation}
        \lim_{\e\to 0}\liminf_{N\to\infty}\frac{1}{r_N}\log\Gamma_{N}(B_\e(\mu))
        =-I(\mu)
        =\lim_{\e\to 0}\limsup_{N\to\infty}\frac{1}{r_N}\log\Gamma_{N}(B_\e(\mu))
    \end{equation}
    holds for any $\mu\in\mathcal{Z}$.
\end{prop}

Besides giving exponential estimates for rare events,
a LDP also yields a law of large numbers when the rate function admits a unique minimizer. 
More precisely, 
the well-known observation below shows that
the corresponding probability measures concentrate asymptotically
on the unique minimizer of the rate function
(see, for instance, \cite[Corollary 3.7]{RubLDP}). 

\begin{prop}\label{Gamma-conv}
    Assume $\mathcal{Z}$ is compact.
    Let $\{\Gamma_{N}\}_{N=1}^\infty$
    be a sequence of probability measures on $\mathcal{Z}$
    satisfying a LDP with rate function $I$.
    If $I$ admits a unique minimizer $z_*\in\mathcal{Z}$,
    then $\Gamma_{N}$ weakly converges to the Dirac measure supported on $z_*$. 
\end{prop}

We will use this observation to deduce the convergence of our random point processes
toward the measure associated with a unique twisted coupled K\"ahler--Einstein metric.

We now apply this general framework to the Gibbs measures naturally associated with
a decomposition of the canonical bundle. 
Our goal is to construct a sequence of random point processes
whose empirical measures satisfy an LDP with the coupled free energy
$F^{\rm cp}_\b$ as a rate function.

\subsubsection{Large deviation principle for a Gibbs measure
associated with a decomposition of the canonical bundle}
\label{LDP-Gibbs}

Let $X$ be an $n$-dimensional canonically polarized manifold,
and fix a decomposition $K_X=L_1+\cdots +L_m$ 
into ample $\mathbb{Q}$-line bundles
satisfying
\begin{equation}\label{norm_vol}
    V_1=\cdots=V_m.
\end{equation}
For sufficiently divisible $k$,
set $N_i:=\dim H^0(X,kL_i)$ and $N=\max_iN_i$.
After relabeling, assume
$N_1\leq N_2 \leq\cdots\leq N_m=N$.
The assumption \eqref{norm_vol} means that the growth of $N_1,\dots,N_m$ as $k\to\infty$
are comparable, that is, $\lim\limits_{k\to\infty}N_i/N=1$ for any $i$.
This comparability will be used to define 
a mutually asymptotic Fekete configuration
in Definition \ref{mutually Fekete}.

One may define the Slater determinant
$\det S_i^{(k)}\in H^0(X^{N_i}, (kL_i)^{{\boxtimes}N_i})$ 
as a generator of
the line $\Lambda^{N_i}H^0(X, kL_i)$.
More precisely,
for each $i$, we define a holomorphic section on the product $X^{N_i}$ by
\begin{equation}
\label{Slater}
{\rm det} S_i^{(k)} (x_1, ... , x_{N_i} ) :=  \det (s^i_\alpha (x_\beta))_{\alpha \beta} \in H^{0} (X^{N_i}, \otimes_j p_j^* L_i^k),
\end{equation}
where $(s^i_\alpha)_\alpha$
is the orthonormal basis specified below,
and $p_j:X^N\to X$ is the natural projection.
In order to consider the Slater determinant in the coupled setting,
let us consider the line bundle
\begin{equation}
{L}_i^\prime:=
\left( \otimes_{j=1}^{N_i} p_j^*  L_i^k \right) \otimes \left( \otimes_{j=N_i + 1}^{N} p_j^* \mathscr{O}_X \right)
\to X^N.
\end{equation}
We identify $\det S_i^{(k)}$ with a holomorphic section
$\pi_i^* \det S_i^{(k)}$ of $L_i^\prime$,
where
\begin{equation}
\pi_i : X^N \to X^{N_i} \,\, ; \,\, \left(x_1 , ..., x_{N_i} , x_{N_i +1} ,... , x_{N} \right) \mapsto \left(x_1 , ..., x_{N_i} \right).
\end{equation}

Retain the Hermitian metric $h_i$ and the compatible volume form $dV$
fixed in Section \ref{functionals}.
For each $i=1,\ldots,m$ and each sufficiently large and divisible
integer $k$, choose an orthonormal basis
$S_i^{(k)}=(s_\alpha^i)_{\alpha=1}^{N_i}$ of $H^0(X,kL_i)$
with respect to the inner product
\begin{equation}
 \langle s,t\rangle_{i,k,0}
 := \int_X h_i^{\otimes k}(s,t)\,dV.
\end{equation}
These bases are used throughout Section \ref{LDP-cpKE}
and are kept fixed.
For $\phi_i \in C^0(X)$, 
we can define a Hermitian metric on the line bundle $L_i^\prime$ by
\begin{equation}
H(i, k, \phi_i) (u,u) \left(x_1 , ..., x_{N_i} , x_{N_i +1} ,... , x_{N} \right) 
:= |u|^2 \left(  \prod_{j=1}^{N_i} h_i^k e^{-k\phi_i}    (x_j) \right) e^{-k (\phi_i (x_{N_i +1}) + ... + \phi_i (x_{N}))},
\end{equation}
where $u$ denotes a local holomorphic frame of ${L}_i^\prime$.
Let us consider the (non-symmetric) Hamiltonian
\begin{equation}\label{Sing-Ham}
    H^{(N)}(x_1,\dots,x_N):=
    -\log\prod_{i=1}^m\left|\pi_i^*\det S_i^{(k)}(x_1,\dots,x_N)\right|_{H(i,k,0)}^{2/k}
\end{equation}
to define the Gibbs measure
\begin{equation}\label{Gibbs-meas}
\mu_\b^{(N)}:=\frac{e^{-\b H^{(N)}}dV^{\otimes N}}{Z_{N,\b}}.
\end{equation}
Here, $\b$ is a fixed positive constant and
$Z_{N,\b}$ is the normalizing constant defined by
$$Z_{N,\b}:=\int_{X^N}e^{-\b H^{(N)}}dV^{\otimes N}$$
and is called the ($N$-particle) partition function.
The ensemble $(X^N, \mu^{(N)}_\b)$ defines a random point process
with $N$ particles on $X$.
\begin{rem}
\label{RPP}
    Although the measure $\mu_{\beta}^{(N)}$ is not invariant
    under the action of the symmetric group on $X^N$, 
    we refer to $(X^N,\mu_{\beta}^{(N)})$ as a random point process 
    by a slight abuse of terminology. 
    On the other hand, the averaged $r$-point correlation measure 
    $\rho_{\beta,r}^{(N)}$, defined below, is invariant 
    under the action of the symmetric group on $X^r$.
\end{rem}
The corresponding empirical measure is
the map $\delta_N:X^N\to\mathcal P(X)$ defined by 
\begin{equation}
\delta_N:=\frac{1}{N}\sum_{j=1}^N\delta_{x_j},
\end{equation}
where $\delta_x$ is the Dirac measure for $x\in X$.
Its law under $\mu^{(N)}_\beta$ is
\begin{equation}
\Gamma_{N,\b}:=(\delta_N)_*\mu^{(N)}_\b\in\mathcal{P}(\mathcal{P}(X)).
\end{equation}

Now we state the main theorem in Section \ref{LDP-cpKE}.

\begin{thm}\label{cp-LDP}
    Assume Conjecture \ref{Hul-conj} below holds.
    Then the law $\Gamma_{N,\b}:=(\delta_N)_*\mu^{(N)}_\b$ satisfies 
    a large deviation principle with speed $N$
    whose rate function is the (inf-normalized) coupled free energy
    $F^{\rm cp}_\b-\inf_{\mathcal{P}(X)}F^{\rm cp}_\b$.
    In particular, as $N\to\infty$, the law
$\Gamma_{N,\b}$ converges weakly to the Dirac measure
$\delta_{\mu_{\rm tcKE}}$,
where $\mu_{\rm tcKE}$ is the Monge–Amp\`ere measure of the unique twisted coupled K\"ahler–Einstein metric.

\end{thm}

Note that we have
$0\leq\inf F_\beta^{\rm cp}<+\infty$.
This is because $\beta>0$, both functionals $\mathrm{Ent}_{dV}$ and $E^{\rm cp}$ are
non-negative and $F_\beta^{\rm cp}(dV)<+\infty$.

The proof of Theorem \ref{cp-LDP} consists of proving
Theorems \ref{upperLDP} and \ref{lowerLDP} and
Corollary \ref{lim-Z_N}.
Theorems \ref{upperLDP} and \ref{lowerLDP} establish local upper and lower
bounds for the unnormalized measures $Z_{N,\b}\Gamma_{N,\b}$.
Corollary \ref{lim-Z_N} determines the asymptotic normalization
$N^{-1}\log Z_{N,\beta}\to-\inf F^{\mathrm{cp}}_\beta$.
Subtracting this limit from the local bounds and applying
\cite[Lemma 3.4]{Ber17} yield the LDP in Theorem \ref{cp-LDP}.
Proposition \ref{Gamma-conv} then gives weak convergence to
$\delta_{\mu_{\mathrm{tcKE}}}$.

As a consequence of the weak convergence of $\Gamma_{N,\b}$,
we construct a canonical sequence of measures converging weakly to the measure
defining a unique twisted coupled K\"ahler--Einstein metric.
More generally, for a fixed positive integer $r$, set
$(N)_r:=N(N-1)\cdots(N-r+1)$
and define the averaged $r$-point correlation measure by
$$
  \rho_{\b,r}^{(N)}
  :=
  \frac{1}{(N)_r}
  \sum_{\substack{1\leq i_1,\ldots,i_r\leq N\\
  i_1,\ldots,i_r\ {\rm distinct}}}
  (\operatorname{pr}_{i_1,\ldots,i_r})_*\mu^{(N)}_\b
  \in\mathcal{P}(X^r),
$$
where
$
\operatorname{pr}_{i_1,\ldots,i_r}:X^N\to X^r
$
is the projection onto the $i_1,\ldots,i_r$-th factors, i.e.
$
(x_1,\ldots,x_N)\mapsto(x_{i_1},\ldots,x_{i_r}).
$
Note that for fixed $\b$, the measures  $\mu^{(N)}_\b$ and $\rho_{\b,r}^{(N)}$
are determined by the decomposition, the reference Hermitian metrics, and the volume form $dV$. 
The normalization makes them independent of the choice of bases.

\begin{cor}\label{r-correlation}
For each fixed positive integer $r$,
if $\Gamma_{N,\b}\to\delta_{\mu_{\rm tcKE}}$,
then the averaged $r$-point correlation measure
$\rho_{\b,r}^{(N)}$ converges weakly to $\mu_{\rm tcKE}^{\otimes r}$.
\end{cor}

As in the case $m=1$, 
if $\mu^{(N)}_\beta$ is symmetric under permutations of the $N$ factors, 
all ordered $r$-point marginals coincide, 
and the conclusion follows from the standard argument \cite[Lemma 4.1]{Ber18}. 
Since our Gibbs measure need not be symmetric, 
we instead average over all ordered $r$-tuples of distinct indices.
    
\subsubsection{Proof of a local upper estimate for the unnormalized measure
$Z_{N,\b}\Gamma_{N,\b}$}

We first establish the local upper estimate for the unnormalized measure $Z_{N,\beta}\Gamma_{N,\beta}$. 
The argument follows the general strategy in Berman's proof of the LDP in the uncoupled setting \cite{Ber18}.
The main point is that the contribution of each Slater determinant can be controlled
from above by its supremum norm,
whose asymptotic behavior is governed by the corresponding Aubin--Mabuchi energy.
Combining these estimates with Sanov's theorem for the reference product measure
yields the coupled pluricomplex energy and the entropy term, respectively.
This leads naturally to the coupled free energy $F^{\mathrm{cp}}_\beta$.

\begin{thm}\label{upperLDP}
For any $\mu\in\mathcal{P}(X)$,
we have
    \begin{equation}
        \lim_{\e\to 0}\limsup_{k\to \infty}
        \frac{1}{N}\log\int_{\delta_N^{-1}(B_{\e}(\mu))}
        e^{-\b H^{(N)}}
        dV^{\otimes N}
        \leq -F^{\mathrm{cp}}_\b(\mu).  
    \end{equation}
\end{thm}
\begin{proof}
Fix $\phi_i\in\mathcal{H}(L_i)$ and $\e>0$.
We have
\begin{eqnarray*}
    &&\int_{\delta_N^{-1}(B_{\e}(\mu))}
        e^{-\b H^{(N)}}
        dV^{\otimes N} \\
    &=&\int_{\delta_N^{-1}(B_{\e}(\mu))}
        \prod_{i=1}^m\left|\pi_i^*\det S_i^{(k)}\right|_{H(i,k,\phi_i)}^{2\b/k}
        \prod_i^me^{\b\sum_{j=1}^N\phi_i(x_j)}
        dV^{\otimes N} \\
    &\leq&\sup_{\delta_N^{-1}(B_{\e}(\mu))}
        \prod_{i=1}^m\left|\pi_i^*\det S_i^{(k)}
        \right|_{H(i,k,\phi_i)}^{2\b/k}
        \int_{\delta_N^{-1}(B_{\e}(\mu))}\mu_{\b\phi}^{\otimes N},
\end{eqnarray*}
where $\mu_{\b\phi}:=\exp({\b\sum_i\phi_i})dV$.

It follows from \cite[Corollary A (i)]{BB10}
and Sanov's theorem (see \cite[Proposition 4.6]{Ber18} for example) that

\begin{eqnarray*}
    &&\lim_{\e\to 0}\limsup_{k\to\infty}\frac{1}{N}\log
       \int_{\delta_N^{-1}(B_{\e}(\mu))}
       e^{-\b H^{(N)}}
       dV^{\otimes N} \\
    &\leq& \sum_{i=1}^m\limsup_{k\to\infty}\frac{1}{N}\log\sup_{X^N}
           \left|\pi_i^*\det S_i^{(k)}\right|_{H(i,k,\phi_i)}^{2\b/k}
           +\lim_{\e\to 0}\limsup_{k\to\infty}\frac{1}{N}\log
            \int_{\delta_N^{-1}(B_{\e}(\mu))}\mu_{\b\phi}^{\otimes N} \\
    &\leq& \sum_{i=1}^m -\b\mathcal{E}_i(\phi_i)-\mathrm{Ent}_{\mu_{\b\phi}}(\mu) \\
           &=& \sum_{i=1}^m\b\left(-\mathcal{E}_i(\phi_i)+\int_X\phi_i\mu\right)
           -\mathrm{Ent}_{dV}(\mu).
\end{eqnarray*}
Therefore,
\begin{eqnarray*}
    &&\lim_{\e\to 0}\limsup_{k\to\infty}\frac{1}{N}\log
       \int_{\delta_N^{-1}(B_{\e}(\mu))}
       e^{-\b H^{(N)}}
       dV^{\otimes N} \\
    &\leq& \sum_{i=1}^m\inf_{\phi_i\in\mathcal{H}_i}
           \b\left(-\mathcal{E}_i(\phi_i)+\int_X\phi_i\mu\right)
           -\mathrm{Ent}_{dV}(\mu) \\
    &=& -F^{\mathrm{cp}}_\b(\mu).
\end{eqnarray*}
This completes the proof.
\end{proof}

\subsubsection{Proof of a local lower estimate for the unnormalized measure
$Z_{N,\b}\Gamma_{N,\b}$}

We next turn to the local lower estimate. In contrast to the upper estimate, 
one needs to find configurations whose empirical measures approximate a prescribed measure $\mu$
while simultaneously controlling $H^{(N)}$ from above.
Two ingredients are essential for this purpose. 
The first is a submean inequality, 
which allows us to compare the value of the Gibbs density at a suitable configuration with its integral over a neighborhood of that configuration. 
The second is the existence of mutually asymptotic Fekete configurations, 
which provides configurations whose empirical measures converge to the prescribed equilibrium measure 
and whose Slater determinants have the desired asymptotic behavior. 
We begin by recalling these two ingredients.

We first prove the submean inequality.
By the Calabi-Yau theorem, we choose
an auxiliary K\"ahler metric $g$ on $X$ whose Riemannian volume
form is proportional to $dV$.
In what follows, $B_r(\mu)$ denotes the open ball of radius $r$
centered at $\mu\in\mathcal P(X)$ for the $L^2$-Wasserstein distance $W_2$ induced
by $g$. Since $X$ is compact, this distance induces the weak topology
on $\mathcal P(X)$.

\begin{prop}\label{submean}
For $x\in X^N$ and $r>0$, write
\begin{equation}
 A_r(x):=\delta_N^{-1}\bigl(B_r(\delta_N(x))\bigr).
\end{equation}
There exist a constant $C_0>0$ and, for each $\e\in(0,1]$,
a constant $C_\e>0$ such that
\begin{equation}\label{submean-ineq}
 e^{-\beta H^{(N)}(x)}
 \le C_\e e^{C_0\e N}
 \frac{\int_{A_\e(x)}
           e^{-\beta H^{(N)}}\,dV^{\otimes N}}
      {\int_{A_{\e^2}(x)}dV^{\otimes N}}
\end{equation}
for $x\in X^N$.
The constants are independent of $k$ and $x$, and $C_0$ is also
independent of $\e$.
\end{prop}
\begin{proof}
Put $q:=\min_iN_i$, and let $S_q\subset S_N$ be the symmetric
group acting on $X^N$ by permuting the first $q$ coordinates. 
For $t>0$, define
\begin{equation}
 r_{i,k}:=\left|\pi_i^*\det S_i^{(k)}\right|_{H(i,k,0)}^2
 \quad\text{and}\quad
 H_{k,t}:=-\frac1k\sum_{i=1}^m\log(r_{i,k}+t),
\end{equation}
and set $f_{k,t}:=e^{-\beta H_{k,t}}$ and
$v_{k,t}:=f_{k,t}^{1/2}$. These functions are smooth and $S_q$-invariant.
Moreover, $f_{k,t}\to e^{-\beta H^{(N)}}$ pointwise as $t\downarrow0$.

The curvature of the line bundle carrying $\pi_i^*\det S_i^{(k)}$,
with its reference metric, is
$k\sum_{\ell=1}^{N_i}p_\ell^*\omega_i$. Consequently,
\begin{equation}
 \sqrt{-1}\,\partial\bar\partial\log(r_{i,k}+t)
 \ge -k\sum_{\ell=1}^{N_i}p_\ell^*\omega_i.
\end{equation}
Indeed, in a local frame with $r_{i,k}=|F|^2e^{-\Psi}$, the function
$\log(r_{i,k}+t)+\Psi=\log(|F|^2+t e^\Psi)$ is plurisubharmonic.
Taking the trace in each coordinate gives
\begin{equation}\label{coord-bound}
 \Delta_{x_j}H_{k,t}\le C
\end{equation}
for $1\le j\le N$,
where $C>0$ is independent of $k$, $t$, and $j$. Here
$\Delta_{x_j}$ is the (negative) Laplacian of $g$ in the $j$-th coordinate.
The summand indexed by $i$ is independent of $x_j$ when
$j>N_i$. It follows that
\begin{equation}
 \Delta_{x_j}v_{k,t}
 =v_{k,t}\left(-\frac\beta2\Delta_{x_j}H_{k,t}
       +\frac{\beta^2}{4}|\nabla_{x_j}H_{k,t}|^2\right)
 \ge -\frac{\beta C}{2}\,v_{k,t}.
\end{equation}

Equip $X^N$ with $g_N:=N^{-1}g^{\oplus N}$, and denote its distance
by $d_N$. Since $\Delta_{g_N}=N\sum_{j=1}^N\Delta_{x_j}$, we have
\begin{equation}\label{prod-bound}
 \Delta_{g_N}v_{k,t}\ge -\frac{\beta C}{2}N^2v_{k,t}.
\end{equation}
If $\Ric_g\ge-Kg$, then
$\Ric_{g_N}\ge-KN g_N$. As the real dimension is
$2nN$, the Ricci lower bound required in
\cite[Theorem 2.1]{Ber17} is uniform in $N$.
For $r>0$, put
\begin{equation}
 Q_r(x):=\bigcup_{\tau\in S_q}B_{d_N}(\tau x,r).
\end{equation}
This is the inverse image of the radius $r$ ball centered at $[x]$ in $X^N/S_q$
with respect to the quotient map $X^N\to Y_N:=X^N/S_q$.
Apply \cite[Theorem 2.1]{Ber17} to the function $v_{k,t}$ on $Y_N$, using
\eqref{prod-bound}, with
$\delta=\e$, $\eta=\e$, and
$\lambda=\sqrt{\beta C/2}\,N$.
We obtain
\begin{equation}\label{part-quot}
 f_{k,t}(x)
 \le C_\e e^{a\e N}
 \frac{\int_{Q_\e(x)}f_{k,t}\,dV^{\otimes N}}
      {\int_{Q_{e^2}(x)}dV^{\otimes N}},
\end{equation}
where $a>0$ is independent of $k,t$, and $\e$, and
$C_\e$ is independent of $k$ and $t$.

We next identify the Wasserstein balls 
with unions of balls for the product metric. 
As discussed in the proof of \cite[Proposition 3.9]{Ber17},
the $L^2$-Wasserstein distance between 
$\delta_N(x)$ and $\delta_N(y)$
can be expressed as
\begin{equation}
    W_2\bigl(\delta_N(x),\delta_N(y)\bigr)^2
=\min_{\sigma\in S_N}\frac{1}{N}\sum_{j=1}^N
d_g(x_{\sigma(j)},y_j)^2
=\min_{\sigma\in S_N}d_N(\sigma x,y)^2.
\end{equation}
In particular, $y\in A_r(x)$ if and only if
$d_N(\sigma x,y)<r$ for some $\sigma\in S_N$, which proves
\begin{equation}
 A_r(x)=\bigcup_{\sigma\in S_N}B_{d_N}(\sigma x,r).
\end{equation}
Choose left coset representatives
$S_N=\coprod_{\alpha=1}^{M_k}\sigma_\alpha S_q$, where
$M_k:=N!/q!$. Then
\begin{equation}
 A_r(x)=\bigcup_{\alpha=1}^{M_k}\sigma_\alpha Q_r(x).
\end{equation}
Writing $\Vol_N(E):=\int_EdV^{\otimes N}$ and using
the invariance of the product measure under coordinate permutations,
we obtain
\begin{equation}
 Q_r(x)\subset A_r(x)
\quad\text{and}\quad
 \Vol_N(A_r(x))
 \le M_k\Vol_N(Q_r(x)).
\end{equation}
Substitution into \eqref{part-quot} yields
\begin{equation}
 f_{k,t}(x)
 \le M_k C_\e e^{a\e N}
 \frac{\int_{A_\e(x)}f_{k,t}\,dV^{\otimes N}}
      {\int_{A_{\e^2}(x)}dV^{\otimes N}}.
\end{equation}
For each fixed $k$, bounded convergence as $t\downarrow0$ gives
\begin{equation}\label{index-loss}
 e^{-\beta H^{(N)}}(x)
 \le M_k C_\e e^{a\e N}
 \frac{\int_{A_\e(x)}e^{-\beta H^{(N)}}\,dV^{\otimes N}}
      {\int_{A_{\e^2}(x)}dV^{\otimes N}}.
\end{equation}

Finally, by normalization
$L_1^n=\cdots=L_m^n=:V$
in \eqref{norm_vol},
the Hilbert
polynomial gives
\begin{equation}
 N_i(k)=\frac{V}{n!}k^n+O(k^{n-1}).
\end{equation}
for any $i=1,\dots,m$.
Thus $N-q=O(k^{n-1})$, $N=O(k^n)$, and
\begin{equation}
 0\le\frac{\log M_k}{N}
 =\frac1N\sum_{\ell=q+1}^{N}\log\ell
 \le\frac{(N-q)\log N}{N}
 =O\!\left(\frac{\log k}{k}\right)\to 0.
\end{equation}
For each fixed $\e>0$, it follows that
$M_k\le D_\e e^{\e N}$ with a constant
$D_\e$ independent of $k$.
This completes the proof of the inequality \eqref{submean-ineq}.
\end{proof}

The second observation concerns the mutually Fekete configuration
in the sense of Hultgren \cite{Hul}.
Let us consider $m$ pairs $\{(L_i,\phi_i)\}_i$
where $\phi_i$ is a continuous Hermitian metric on $L_i$
relative to $h_i$.

\begin{definition}\label{mutually Fekete}
    A sequence $\{P_k\}_{k=1}^\infty$, with $P_k\in X^N$,
    is called a mutually asymptotic Fekete configuration
    for pairs $\{(L_i,\phi_i)\}_i$
    if for all $i=1,\dots,m$,
    \begin{equation}
        \lim_{k\to\infty}\frac{1}{kN}\log
        \left|\pi_i^*\det S_i^{(k)}(P_k)\right|_{H(i,k,\phi_i)}^2
        =
        \lim_{k\to\infty}\frac{1}{kN}\log\max_{X^N}
        \left|\pi_i^*\det S_i^{(k)}\right|_{H(i,k,\phi_i)}^2.
    \end{equation}
\end{definition}

Let $\mu_{(L_i,\phi_i)}$ be the equilibrium measure of $(L_i,\phi_i)$,
that is, $\mu_{(L_i,\phi_i)}=d_{\phi_i}(\mathcal{E}_i \circ P_{\omega_i})$.
One of the main roles of the mutually Fekete configuration
is to approximate each $\mu_{(L_i,\phi_i)}$.

\begin{prop}[See the proof of Lemma 1 in\cite{Hul}]\label{Fekete_conv}
    Let $\{P_k\}_{k=1}^\infty\subset X^N$ be
    a mutually asymptotic Fekete configuration
    for pairs $\{(L_i,\phi_i)\}_i$.
    Then for any $i=1,\dots,m$, we have
    \begin{equation}
        \lim_{k\to\infty}\delta_N(P_k)=\mu_{(L_i,\phi_i)}.
    \end{equation}
\end{prop}

When $m=1$, the Fekete configuration always exists \cite{BB10}.
In our general coupled case, Hultgren \cite{Hul} conjectures the following.
\begin{conj}\label{Hul-conj}
    A mutually asymptotic Fekete configuration exists
    for $\{(L_i,\phi_i)\}_i$
    if and only if
    \begin{equation}\label{mu-balance}
        \mu_{(L_1,\phi_1)}=\cdots=\mu_{(L_m,\phi_m)}.
    \end{equation}
\end{conj}

The direction ``only if" follows from Proposition \ref{Fekete_conv}.
The converse is established in the toric setting under the regularity assumptions of
\cite[Theorem 1]{Hul}, in particular when the weights are smooth.

Summarizing the above observations, we obtain the lower bound of a LDP
provided that Conjecture \ref{Hul-conj} is true.
\begin{thm}\label{lowerLDP}
    If Conjecture \ref{Hul-conj} is true, then for any $\mu\in\mathcal{P}(X)$,
    we have
    \begin{equation}
        \lim_{\e\to 0}\liminf_{k\to \infty}
        \frac{1}{N}\log\int_{\delta_N^{-1}(B_{\e}(\mu))}
        e^{-\b H^{(N)}}
        dV^{\otimes N}
        \geq -F^{\mathrm{cp}}_\b(\mu).
    \end{equation}
\end{thm}
\begin{proof}
    Fix $\mu\in\mathcal{P}(X)$ and $\e\in(0,1]$.
    It suffices to consider the case where $E^{\rm cp}(\mu)<\infty$,
    i.e. $\mu\in\mathcal{M}^1(X)$.
    We first assume that
    there exists a continuous Hermitian metric $\phi_i$ on $L_i$
    such that
    \begin{equation}\label{mu=equiv}
      \mu_{(L_i,\phi_i)}=\mu  
    \end{equation}
    for any $i=1,\dots,m$.   
    By our assumptions, the pair $\{(L_i,\phi_i)\}_i$ admits
    a mutually asymptotic Fekete configuration $\{P_k\}_{k=1}^\infty$.
    Put $P_k=(P_{k,1},P_{k,2},\dots,P_{k,N})\in X^N$.
    It follows from Proposition \ref{Fekete_conv} that
    $B_{\e}(\delta_N(P_k))\subset B_{2\e}(\mu)$
    for any large $k\gg1$.
    We can apply
    the submean inequality in
    Propositions \ref{submean} to show
    \begin{eqnarray*}
        &&\liminf_{k\to\infty}\frac{1}{N}\log
        \int_{\delta_N^{-1}(B_{2\e}(\mu))}
        e^{-\b H^{(N)}}
        dV^{\otimes N} \\
        &\geq&
        \liminf_{k\to\infty}\frac{1}{N}\log
        \left(
        \prod_{i=1}^m\left|\pi_i^*\det S_i^{(k)}(P_k)
        \right|_{H(i,k,\phi_i)}^{2\b/k}
        \right)
        +\liminf_{k\to\infty}\frac{1}{N}\log
        \prod_{i=1}^me^{\b\sum_{j=1}^N\phi_i(P_{k,j})}\\
        &&+\liminf_{k\to\infty}\frac{1}{N}\log
        \int_{\delta_N^{-1}(B_{\e^2}(\delta_N(P_k)))}dV^{\otimes N}
        +\liminf_{k\to\infty}-\frac{1}{N}\log C_\e-C_0\e.
    \end{eqnarray*}
    The first term in the right hand side equals
     $\beta\sum_{i=1}^m-(\mathcal{E}_i\circ P_{\omega_i})(\phi_i)$
    by definition of a mutually asymptotic Fekete configuration
    and \cite[Corollary A (i)]{BB10}.
    The Second term in the right hand side equals
    $\beta\sum_{i=1}^m\int_X\phi_i\mu$
    by Proposition \ref{Fekete_conv} and the assumption \eqref{mu=equiv}.
    Now we can assume $B_{\e^2/2}(\mu)\subset B_{\e^2}((\delta_N(P_k)))$
    for any large $k\gg 1$ by Proposition \ref{Fekete_conv}.
    Thus, by Sanov's theorem, we have
    \begin{eqnarray*}
        &&\lim_{\e\to 0}\liminf_{k\to\infty}\frac{1}{N}\log
        \int_{\delta_N^{-1}(B_{2\e}(\mu))}
        e^{-\b H^{(N)}}
        dV^{\otimes N} \\
        &\geq&
        \sum_{i=1}^m\b\left(-(\mathcal{E}_i\circ P_{\omega_i})(\phi_i)
        +\int_X\phi_i\mu\right)
        -\lim_{\e\to 0}\inf_{B_{\e^2/2}(\mu)}\mathrm{Ent}_{dV} \\
        &\geq&
        \sum_{i=1}^m\b\left(-(\mathcal{E}_i\circ P_{\omega_i})(\phi_i)
        +\int_X\phi_i\mu\right)-\mathrm{Ent}_{dV}(\mu).
    \end{eqnarray*}
    By definition of $\phi_i$ (recall Remark \ref{projection}),
    \begin{equation}
    -(\mathcal{E}_i\circ P_{\omega_i})(\phi_i)+\int_X\phi_i\mu=-E_i(\mu).
    \end{equation}
    This completes the proof for the case when we assume the condition \eqref{mu=equiv}.

    To treat a general finite-energy measure, 
    we approximate its normalized potential in the $d_1$-topology 
    and transfer the resulting strong convergence to each polarization.
    Let $\phi_\mu\in\mathcal{E}^1(X,\omega)$
    be a potential of $\mu$
    for a K\"ahler metric $\omega$, that is, $\mu=\MA_\omega(\phi_\mu)$.
    Take a sequence of K\"ahler potentials
    $\{\phi_\alpha\}_\alpha\subset\mathcal{H}(X,\omega)$
    such that $d_{1,\omega}(\phi_\alpha,\phi_\mu)\to0$.
    In particular, 
    $
        \mathcal{E}_\omega(\phi_\alpha)\to\mathcal{E}_\omega(\phi_\mu).
    $
    Since $\MA_\omega(\cdot)$ gives a homeomorphism from $\mathcal{E}^1(X,\omega) / \mathbb{R}$
    to $\mathcal{M}^1(X)$ (see \cite[Proposition 2.6]{BBEGZ} and Remark \ref{indep of class} in this paper),
    \begin{equation}\label{MA-approx}
     \mu_\alpha:=\MA_\omega(\phi_\alpha)\to\MA_\omega(\phi_\mu).
    \end{equation}
    Thus we also have
    $
        E_\omega(\mu_\alpha)\to E_\omega(\mu).
    $
    For any $i=1,2,\dots,m$, one can define potentials
    $\phi_{\mu,i}\in\mathcal{E}^1(L_i)$
    and $\phi_{i,\alpha}\in\mathcal{H}(L_i)$
    by the conditions
    $
     \mu=\MA_{\omega_i}(\phi_{\mu,i}), 
    $
    $\mu_{\alpha}=\MA_{\omega_i}(\phi_{i,\alpha})
    $
    and $\sup_X \phi_{i, \alpha} = \sup_X \phi_{\mu, i} =0$.
    Then for any $i$, we have
    \begin{equation}\label{E_i-approx}
        E_i(\mu_\alpha)\to
        E_i(\mu)
    \end{equation}
    and $d_{1,\omega_i}(\phi_{i,\alpha},\phi_{\mu,i})\to 0$ as $\alpha\to\infty$.
    By repeating the above argument showing the lower bound
    with $\mu$ replaced by $\MA_\omega(\phi_\alpha)$
    and $\phi_i$ by $\phi_{i,\alpha}$,
    and using the approximation formulas \eqref{MA-approx} and \eqref{E_i-approx},
    we obtain the desired lower bound estimate in the LDP.
    Indeed, since there exists $\alpha_\e\in\mathbb{N}$ such that for any $\alpha\geq\alpha_\e$
    we have $B_{2\e}(\MA_\omega(\phi_\alpha))\subset B_{3\e}(\mu)$,
    the following estimate holds for any $\alpha\geq\alpha_\e$.
    \begin{eqnarray*}
        &&
        \liminf_{k\to\infty}\frac{1}{N}\log
        \int_{\delta_N^{-1}(B_{3\e}(\mu))}
        e^{-\b H^{(N)}}
        dV^{\otimes N} \\
        &\geq&
        \liminf_{k\to\infty}\frac{1}{N}\log
        \int_{\delta_N^{-1}(B_{2\e}(\mu_\alpha))}
        e^{-\b H^{(N)}}
        dV^{\otimes N} \\
        &\geq&
        -C_0\e
        +\sum_{i=1}^m-\b\mathcal{E}_i(\phi_{i,\alpha})
        -\inf_{B_{\e^2/2}(\mu_\alpha)}\mathrm{Ent}_{\mu_{\b\phi}} \\
        &\geq&
        -C_0\e+
        \sum_{i=1}^m\b\left(-\mathcal{E}_i(\phi_{i,\alpha})
        +\int_X\phi_{i,\alpha}\mu\right)-\mathrm{Ent}_{dV}(\mu).
    \end{eqnarray*}
    The third inequality follows from the fact that
    $\mu_\alpha$ is in $B_{\e^2/2}(\mu)$ for sufficiently large $\alpha$.
    Letting $\varepsilon \to 0$ (then $\alpha_\e\to\infty$) and
    applying Lemma \ref{phi_a-approx} below with
    $\omega=\omega_i$, $\phi=\phi_{\mu,i}$, and $\phi_\alpha=\phi_{i,\alpha}$,
    we obtain the following estimate.
    \begin{eqnarray*}
        &&
        \lim_{\e\to 0}\liminf_{k\to\infty}\frac{1}{N}\log
        \int_{\delta_N^{-1}(B_{3\e}(\mu))}
        e^{-\b H^{(N)}}
        dV^{\otimes N} \\
        &\geq&
        \lim_{\alpha\to\infty}
        \sum_{i=1}^m\b\left(-\mathcal{E}_i(\phi_{i,\alpha})
        +\int_X\phi_{i,\alpha}\mu_\alpha\right)-\mathrm{Ent}_{dV}(\mu)\\
        &=&
        \lim_{\alpha\to\infty}\sum_{i=1}^m-\b E_i(\mu_\alpha)
        -\mathrm{Ent}_{dV}(\mu)\\
        &=&-F^{\rm cp}_\b(\mu).
    \end{eqnarray*}
    This completes the proof.
\end{proof}

\begin{lemma}\label{phi_a-approx}
Let $\phi, \phi_\alpha\in\mathcal{E}^1(X,\omega)$
satisfy $\lim\limits_{\alpha\to\infty}d_{1,\omega}(\phi,\phi_\alpha)=0$.
Then
\begin{equation}
\lim_{\alpha\to\infty}\int_X\phi_\alpha
\left(\MA_\omega(\phi)-\MA_\omega(\phi_\alpha)\right)=0.
\end{equation}
\end{lemma}

\begin{proof}
We apply \cite[Propositions 3.4 and 3.7]{Tru23} with the
reference potential $\psi=0$. Since
$d_{1,\omega}(\phi_\alpha,\phi)\to0$, the triangle inequality gives
$
\sup_{\alpha}d_{1,\omega}(\phi_\alpha,0)<+\infty.
$
Hence there exists $C>0$ such that
$
d_{1,\omega}(\phi_\alpha,0)+d_{1,\omega}(\phi,0)\le C
$
for all $\alpha$.

Set
\begin{equation}
I_\omega(u,v)
:=
\int_X (u-v)
\left(\MA_\omega(v)-\MA_\omega(u)\right).
\end{equation}
By \cite[Proposition 3.4]{Tru23}, after enlarging $C$ if
necessary, there exists a continuous increasing function
$f_C\colon\mathbb R_{\geq0}\to\mathbb R_{\geq0}$ with
$f_C(0)=0$ such that
\begin{equation}
I_\omega(\phi,\phi_\alpha)
\le
f_C\bigl(d_{1,\omega}(\phi,\phi_\alpha)\bigr).
\end{equation}
Thus
$
I_\omega(\phi,\phi_\alpha)\to 0.
$

We now apply \cite[Proposition 3.7]{Tru23} with
$
u=\phi_\alpha, v=0,
\varphi_1=\phi,$ and $\varphi_2=\phi_\alpha.
$
There exists a constant $\widetilde C>0$, independent of
$\alpha$, such that
\begin{equation}
\left|
\int_X\phi_\alpha
\left(\MA_\omega(\phi)-\MA_\omega(\phi_\alpha)\right)
\right|
\le
\widetilde C\,
I_\omega(\phi,\phi_\alpha)^{1/2}.
\end{equation}
This completes the proof.
\end{proof}
\subsubsection{Asymptotics for the partition function $Z_{N,\b}$
and the proof of \hbox{Theorem \ref{cp-LDP}}}


We now determine the exponential growth rate of the partition function. 
In view of \cite[Lemma 3.4]{Ber17},
this converts the local bounds for $Z_{N,\beta}\Gamma_{N,\beta}$ 
into an LDP for the probability measures $\Gamma_{N,\beta}$,
and completes the proof of Theorem \ref{cp-LDP}.

\begin{cor}\label{lim-Z_N}
Assume that the estimates in Theorems \ref{upperLDP} and \ref{lowerLDP}
hold for any $\mu\in\mathcal{P}(X)$. Then
\begin{equation}
\lim_{k\to\infty}
\frac{1}{N}\log Z_{N,\beta}
=
-\inf_{\mu\in\mathcal P(X)}
F_\beta^{\rm cp}(\mu).
\end{equation}
\end{cor}


\begin{proof}
Set
$
\widetilde{\Gamma}_{N,\beta}
:=
Z_{N,\beta}\Gamma_{N,\beta}.
$
Then
$
Z_{N,\beta}
=
\widetilde{\Gamma}_{N,\beta}(\mathcal P(X)).
$
Set
$
A_\beta
:=
\inf_{\mu\in\mathcal P(X)}
F_\beta^{\rm cp}(\mu).
$

We first prove the lower bound estimate for the left hand side.
Fix $\eta>0$. By the definition of
$A_\beta$, there exists $\mu_\eta\in\mathcal P(X)$ such that
\begin{equation}
F_\beta^{\rm cp}(\mu_\eta)
\leq A_\beta+\eta.
\end{equation}
For every $\e>0$, we have
\begin{equation}
Z_{N,\beta}
\geq
\widetilde{\Gamma}_{N,\beta}
(B_\varepsilon(\mu_\eta)).
\end{equation}
Hence
\begin{equation}
\liminf_{k\to\infty}
\frac{1}{N}\log Z_{N,\beta}
\geq
\liminf_{k\to\infty}
\frac{1}{N}
\log
\widetilde{\Gamma}_{N,\beta}
(B_\varepsilon(\mu_\eta)).
\end{equation}
Letting $\e\to0$ and applying the local lower bound
in Theorem \ref{lowerLDP}, we obtain
\begin{equation}
\liminf_{k\to\infty}
\frac{1}{N}\log Z_{N,\beta}
\geq
-F_\beta^{\rm cp}(\mu_\eta)
\geq
-A_\beta-\eta.
\end{equation}
Since $\eta>0$ is arbitrary,
\begin{equation}
\liminf_{k\to\infty}
\frac{1}{N}\log Z_{N,\beta}
\geq
-A_\beta.
\end{equation}

We next prove the upper bound estimate for the left hand side.
Fix $\eta>0$.
By the local upper bound in Theorem \ref{upperLDP}, for every
$\mu\in\mathcal P(X)$ there exists $r_\mu>0$, independent of $k$,
such that
\begin{equation}
\limsup_{k\to\infty}
\frac{1}{N}
\log
\widetilde{\Gamma}_{N,\beta}
(B_{r_\mu}(\mu))
\leq
-A_\beta+\eta.
\end{equation}
Since $X$ is compact, $\mathcal P(X)$ is compact with respect to the
weak topology. Hence there exist
$\mu_1,\ldots,\mu_M\in\mathcal P(X)$ such that
\begin{equation}
\mathcal P(X)
\subset
\bigcup_{j=1}^M
B_{r_{\mu_j}}(\mu_j).
\end{equation}
Notice that $M$ is independent of $k$.

For each $j$, by the definition of the limsup, there exists
$k_j\in\mathbb N$ such that
\begin{equation}
\widetilde{\Gamma}_{N,\beta}
(B_{r_{\mu_j}}(\mu_j))
\leq
\exp\left\{
N(-A_\beta+2\eta)
\right\}
\end{equation}
for every $k\geq k_j$. Thus, setting
$
k_0:=\max_{1\leq j\leq M}k_j,
$
we obtain, for every $k\geq k_0$,
\begin{eqnarray*}
Z_{N,\beta}
&=&
\widetilde{\Gamma}_{N,\beta}(\mathcal P(X))
\\
&\leq&
\sum_{j=1}^M
\widetilde{\Gamma}_{N,\beta}
(B_{r_{\mu_j}}(\mu_j))
\\
&\leq&
M
\exp\left\{
N(-A_\beta+2\eta)
\right\}.
\end{eqnarray*}
Therefore,
\begin{equation}
\frac{1}{N}\log Z_{N,\beta}
\leq
\frac{\log M}{N}
-A_\beta+2\eta.
\end{equation}
Since $M$ is independent of $k$,
we have
$
\frac{1}{N}\log M\to0.
$
Consequently,
\begin{equation}
\limsup_{k\to\infty}
\frac{1}{N}\log Z_{N,\beta}
\leq
-A_\beta+2\eta.
\end{equation}
Letting $\eta\to0$, we obtain
\begin{equation}
\limsup_{k\to\infty}
\frac{1}{N}\log Z_{N,\beta}
\leq
-A_\beta.
\end{equation}
This completes the proof.
\end{proof}

\subsubsection{Convergence of the averaged $r$-point correlation measures}

For fixed $r$, sampling $r$ distinct indices is asymptotically equivalent to
sampling with replacement. 
This allows us to deduce convergence of 
the averaged $r$-point correlation measures $\rho^{(N)}_{\b,r}$
from convergence of the law of the empirical measure.

\begin{proof}[Proof of Corollary \ref{r-correlation}]
For $x^{(N)}:=(x_1,\ldots,x_N)\in X^N$, define
\begin{equation}
  \delta_N^{[r]}(x^{(N)})
  :=
  \frac{1}{(N)_r}
  \sum_{\substack{1\leq i_1,\ldots,i_r\leq N\\
  i_1,\ldots,i_r\ {\rm distinct}}}
  \delta_{(x_{i_1},\ldots,x_{i_r})}
  \in\mathcal{P}(X^r),
\end{equation}
where $(N)_r:=N(N-1)\cdots(N-r+1)$.
By definition,
\begin{equation}
  \rho_{\b,r}^{(N)}
  =
  \int_{X^N}
  \delta_N^{[r]}(x^{(N)})\,\mu^{(N)}_\b(x^{(N)}).
\end{equation}
On the other hand,
\begin{equation}
  \delta_N(x^{(N)})^{\otimes r}
  =
  \frac{1}{N^r}
  \sum_{i_1,\ldots,i_r=1}^N
  \delta_{(x_{i_1},\ldots,x_{i_r})}.
\end{equation}

Let $f\in C^0 (X^r)$. Since the only difference between the two sums
above comes from tuples with repeated indices, we have
\begin{equation}
  \left|
  \int_{X^r}f\,\delta_N^{[r]}(x^{(N)})
  -
  \int_{X^r}f\,\delta_N(x^{(N)})^{\otimes r}
  \right|
  \leq
  2\|f\|_\infty
  \left(1-\frac{(N)_r}{N^r}\right).
\end{equation}
Moreover,
$1-(N)_r/N^r$
converges to zero, uniformly in
$x^{(N)}\in X^N$.
It follows that
\begin{equation}
  \int_{X^r}f\,\rho_{\b,r}^{(N)}=
  \int_{X^N}
  \left(
  \int_{X^r}f\,\delta_N(x^{(N)})^{\otimes r}
  \right)
  \mu^{(N)}_\b(x^{(N)})
  +o(1).
\end{equation}
Define
$
\Phi_f:\mathcal{P}(X)\to\mathbb{R}
$
by
\begin{equation}
  \Phi_f(\nu):=
  \int_{X^r}f\,\nu^{\otimes r}.
\end{equation}
The map $\Phi_f$ is bounded and continuous with respect to the weak
topology on $\mathcal{P}(X)$. 
By the definition of $\Gamma_{N,\b}$, we thus obtain
\begin{equation}
  \int_{X^N}
  \left(
  \int_{X^r}f\,\delta_N(x^{(N)})^{\otimes r}
  \right)
  \mu^{(N)}_\b(x^{(N)})
  =
  \int_{\mathcal{P}(X)}
  \Phi_f(\nu)\,\Gamma_{N,\b}(\nu).
\end{equation}
Since
$
  \Gamma_{N,\b}\to\delta_{\mu_{\rm tcKE}},
$
we have
\begin{equation}
  \int_{\mathcal{P}(X)}
  \Phi_f(\nu)\,\Gamma_{N,\b}(\nu)
  \to
  \Phi_f(\mu_{\rm tcKE})
  =
  \int_{X^r}f\,\mu_{\rm tcKE}^{\otimes r}.
\end{equation}
Therefore,
\begin{equation}
  \int_{X^r}f\,\rho_{\b,r}^{(N)}
  \to
  \int_{X^r}f\,\mu_{\rm tcKE}^{\otimes r}
\end{equation}
for every $f\in C^0 (X^r)$.
This completes the proof.
\end{proof}
\section{Coupled free energy and stability when $K_X < 0$}

We now turn to the Fano case. 
We first relate the coercivity of the coupled free energy to 
the coercivity of the coupled Ding functional 
and establish the diagonal criterion. 
We then introduce the non-Archimedean coupled free energy 
and thermodynamical $K$-stability. 
Finally, we derive a coupled volume bound 
and compare the coupled microscopic stability threshold with the coupled delta-invariant.

\subsection{Coercivity and the existence of coupled {\kah}--Einstein metrics}

In this section, we assume that $X$ is a Fano manifold, i.e., the anti-canonical line bundle $-K_X$ is ample.
We assume that there exist ample $\mathbb{Q}$-line bundles $L_i$, $1 \leq i \leq m$ such that
$$
 - K_X = L_1 + \cdot \cdot \cdot + L_m.
$$
We fix {\kah} metrics $\omega_i \in c_1 (L_i)$,
each of which is the curvature form (multiplied by $\sqrt{-1}$) of a fixed Hermitian metric $h_i$ on $L_i$.

Recall that we  say that an $m$-tuple $(\omega_{j, \varphi_j})_j , \, \varphi_j \in \mathcal{H}_i$ is a coupled {\kah}--Einstein metric if the system
\begin{equation}
\label{cKE}
    {\rm Ric}\, \omega_{j, \varphi_j} = \sum_{i=1}^{m}\omega_{i, \varphi_i} 
\end{equation}
holds for any $j=1,...,m$.
This  is equivalent to the following:
$$
\displaystyle \frac{\omega_{j, \varphi_j}^n}{V_j}  = \dfrac{e^{ - \sum_i \varphi_i} dV  }{\int_X  e^{ - \sum_i \varphi_i} dV},\,\,\,\, \forall j =1,...,m.
$$
Here, we fix a smooth volume form $dV$ on $X$ such that $\sqrt{-1} \Theta (-K_X ,dV) = \sum_i \omega_i$ and $\int_X dV =1$.
We can easily check the following.
\begin{prop}
\label{proportional}
Assume that $L_i$ are proportional to $-K_X$, i.e., there exist $k_i \in \mathbb{Q}_{>0}$ such that $\sum_ik_i = 1$ and $L_i = -k_i K_X$.
If $(\omega_{i, \phi_i})_i$ is a coupled {\kah}--Einstein metric, then $\omega_{i,\phi_i}$ are {\kah}--Einstein metrics.
\end{prop}

\begin{proof}
We may assume that $\omega_i = k_i \omega_0$ for some {\kah} form $\omega_0 \in c_1 (X)$.
If we can find a solution $(\omega_{i,\phi_i})_i = (k_i \omega_0 + \dol \phi_i)_i$ of (\ref{cKE}), the uniqueness of the solution of the Calabi--Yau theorem implies that we have
$$
\frac{\phi_1}{k_1} = ... = \frac{\phi_m}{k_m}
$$
up to additive constants.
Then, we obtain ${\rm Ric} \, \omega_{i, \phi_i} = k_i^{-1} \omega_{i, \phi_i}$.
\end{proof}

For $\gamma>0$, we write $\gamma = -\beta$.
The coupled free energy is written as
$$
F^{\rm cp}_{- \gamma } (\mu ) = {\rm Ent}_{dV} (\mu) - \gamma E^{\rm cp} (\mu) , \,\,\,\, \mu \in \mathcal{M}^1 ( X).
$$
If $\mu \in \mathcal{P}(X) \setminus \mathcal{M}^1 (X)$, we define $F_{-\gamma}^{\rm cp}
 (\mu) = +\infty $.
Here, ${\rm Ent}_{dV} $ denotes the relative entropy with respect to the (smooth) probability measure $dV$.

As we have seen in Section 3, the coupled Ding functional plays a very important role in the study of the coupled free energy.

\begin{definition}[coupled Ding functional]
For $\varphi_i \in \mathcal{E}^1 (L_i)$, we define the coupled Ding functional by
$$
D^{\rm cp}_{ - \gamma} ((\varphi_i )_i ) := -\frac{1}{\gamma } \log \int_X e^{ - \gamma \sum \varphi_i } dV  - \sum_i \mathcal{E}_i (\varphi_i).
$$
\end{definition}

A twisted coupled {\kah}--Einstein metric is defined by the following system of equations:
$$
{\rm Ric} \, \omega_{j, \phi_j} = \gamma \sum_i \omega_{i, \phi_i} + (1-\gamma){\rm Ric}\, dV,
$$
equivalently,
\begin{equation}
    \label{tcKE}
\displaystyle \frac{\omega_{j, \varphi_j}^n}{V_j}  = \dfrac{e^{ - \gamma \sum_i \varphi_i} dV  }{\int_X  e^{ - \gamma \sum_i \varphi_i} dV},\,\,\,\, \forall j=1,...,m.
\end{equation}
We can show that a twisted coupled {\kah}--Einstein metric is a critical point of the coupled Ding functional.
Hultgren--Witt Nystr\"{o}m studied the regularity of a minimizer of $D_{-\gamma}^{\rm cp}$ as follows.

\begin{thm}{\rm (\cite[Theorem 2.4]{HW})}
\label{regularity}
If $(\phi_i)_i \in \prod_i \mathcal{E}^1 (L_i)$ is a minimizer of $D_{-\gamma}^{\rm cp}$, then $(\phi_i)_i$ is a weak solution of  (\ref{tcKE}).
Moreover, $(\phi_i)_i$ is smooth and $(\omega_i + \dol \phi_i)_i$ is a twisted coupled {\kah}--Einstein metric.
\end{thm}

The following coercivity property provides the variational criterion for
the existence used below.

\begin{definition}
We say that the coupled Ding functional $D_{- \gamma}^{\rm cp}$ is coercive on $\Pi_i \mathcal{E}^1(L_i)$ if there exist positive constants $C_1 ,C_2$ such that
$$D_{- \gamma}^{\rm cp} \geq C_1 J^{\rm cp} - C_2$$ on 
$\Pi_i \mathcal{E}^1(L_i)$.
Here, the functional $ J^{\rm cp}$ is defined by
$$
 J^{\rm cp} ((\varphi_i)_i) := \sum_i J_i (\phi_i)=  \sum_i \left( \int_X \varphi_i \rm MA_{\omega_i}(0) - \mathcal{E}_i (\varphi_i) \right).
$$
\end{definition}

Note that the functional $J^{\rm cp}$ is compatible with $E^{\rm cp}$ (see Theorem \ref{finite energy} and (\ref{IandJ})).

\begin{thm}{\rm (\cite{HW,Tak,FH})}
\label{Ding criterion}
Fix $\gamma \in (0,1]$.
If $\gamma =1$, we assume that ${\rm Aut}(X)$ is finite.
The unique existence of the twisted coupled {\kah}--Einstein metric is equivalent to the coercivity of the coupled Ding functional $D_{-\gamma}^{\rm  cp}$ on $\prod_i \mathcal{E}^1 (L_i)$.
\end{thm}

We introduce the notion of the diagonal in the product space $\prod_i \mathcal{E}^1 (L_i)$.
We will see later that it is enough to check the coercivity on the ``diagonal" in order to characterize the existence of the twisted coupled {\kah}--Einstein metric.

\begin{definition}
\label{diagonal}
We define the diagonal for $\{L_i \}_i$ by
$$
\Delta = \Delta_{\{L_i\}_i } := \left\{  ( \phi_i)_i \in \prod_i \mathcal{E}^1 (L_i) \, \middle|  \, \,  {\rm MA}_{\omega_i} (\phi_i)= {\rm MA}_{\omega_j} (\phi_j) , \,\forall i,j=1,...,m \right\}
$$
We say that the coupled Ding functional $D_{-\gamma}^{\rm cp}$ is coercive on $\Delta$ if there exist $C_1,C_2 >0$ such that $D_{-\gamma}^{\rm cp}
 \geq C_1 J^{\rm cp} -C_2$ on $\Delta$.
\end{definition}

We introduce the notion of coercivity for the coupled free energy.

\begin{definition}
We say that the coupled free energy $F_{- \gamma}^{\rm cp}$ is coercive if there exist positive constants $C_1 ,C_2>0$ such that $F_{- \gamma}^{\rm cp} \geq C_1 E^{\rm cp} - C_2$ on $\mathcal{M}^{1}(X)$.
\end{definition}

We generalize \cite[Theorem 3.4]{Ber13} to the coupled setting.
\begin{prop}
\label{F geq G}
Fix $\gamma  > 0$. For any $\mu \in \mathcal{M}^1 ( X)$, we have
\begin{equation}
\label{FgeqG}
F^{\rm cp}_{- \gamma } (\mu ) \geq \gamma   D^{\rm cp}_{- \gamma }((\varphi_i)_i) , \,\,\, \mu = \frac{ \omega_{ j , \varphi_j}^n }{ V_j }, \,\,\, \forall j.
\end{equation}
(Note that we can find $(\phi_i)_i \in \prod_i \mathcal{E}^1 (L_i)$ such that each $\phi_i$ solves the above equation by Proposition \ref{pluricomplex}.)
If the equality $\gamma   D^{\rm cp}_{- \gamma }((\varphi_i)_i) = F^{\rm cp}_{- \gamma } (\mu )$ holds, then $(\omega_{i,\phi_i})_i$ is a twisted coupled {\kah}--Einstein metric.
Moreover, if $D_{-\gamma}^{\rm cp}$ is coercive on $\Delta$, then $F_{-\gamma}^{\rm cp}$ is coercive.
\end{prop}

\begin{proof}

Pick sequences $(u_i^j) \subset C^0(X)$ such that $u_i^j \searrow \phi_i$ as $j \to \infty$ for each $i = 1, ...,m$.
Then, we have
\begin{eqnarray*}
F^{\rm cp}_{- \gamma } (\mu )
&=& \sup_{\phi \in C^0 (X)} \left( - \log \int_X e^{\phi} dV  + \int_X \phi \mu \right) - \gamma  E^{\rm cp} (\mu) \\
&\geq & - \log \int_X e^{ - \gamma \sum_i  u_i^j } dV- \gamma  \sum_i \mathcal{E}_i (\varphi_i)  - \gamma  \int_X \sum_i  u_i^j \mu + \gamma \int_X \sum_i  \varphi_i\mu .
\end{eqnarray*}
Then, the monotone convergence theorem implies that the right hand side converges to $- \log \int_X e^{ - \gamma \sum_i  \phi_i } dV- \gamma  \sum_i \mathcal{E}_i (\varphi_i) +  0 = \gamma D_{-\gamma}^{\rm cp} ((\phi_i )_i)$ and we obtain the inequality (\ref{FgeqG}).

 If the equality $\gamma   D^{\rm cp}_{- \gamma }((\varphi_i)_i) = F^{\rm cp}_{- \gamma } (\mu )$ holds, we directly obtain the following:
 $$
 {\rm Ent}_{dV} (\mu)=\int_X \log \frac{\mu}{dV} \mu = - \log \int_X e^{-\gamma \sum_i \phi_i } dV - \gamma \int_X \sum_i \phi_i \mu.
 $$
 The definition of ${\rm Ent}_{dV} $ and the fact that the function $\phi \mapsto -\log \int_X e^\phi dV$ is strictly
 concave imply that the equations
 $$
 \frac{\omega_{j,\phi_j}}{V_j} =\mu = \frac{e^{-\gamma \sum_i \phi_i} dV}{ \int_
 Xe^{-\gamma \sum_i \phi_i} dV}
 $$
 hold for any $j$.
 So, Theorem \ref{regularity} implies that $(\omega_{j,\phi_j})_j$ is a twisted coupled {\kah}--Einstein metric.
By the compatibility of $J^{\rm cp}$ and $E^{\rm cp}$, the coercivity of $D_{-\gamma}^{\rm cp}$ on the diagonal $\Delta \subset \prod_i \mathcal{E}^1(L_i)$ implies that the coupled free energy $F_{-\gamma}^{\rm cp}$ is coercive on $\mathcal{M}^{1}(X)$. 
\end{proof}

We next obtain the reverse inequality by associating a Gibbs measure with a tuple of potentials.
This is also a coupled version of \cite[Theorem 3.4]{Ber13}.
\begin{prop}
\label{converse ineq}
For $(\varphi_i)_i \in \prod_{i} \mathcal{E}^1 (L_i) $, we set
$$
\mu = \frac{ e^{ - \gamma \sum_i  \varphi_i } dV }{\int_X e^{ - \gamma \sum_i  \varphi_i } dV}.
$$
Then, we have
$$
 \gamma   D^{\rm cp}_{- \gamma }((\varphi_i)_i) \geq F^{\rm cp}_{- \gamma } (\mu ) .
$$
If the equality $\gamma   D^{\rm cp}_{- \gamma }((\varphi_i)_i) = F^{\rm cp}_{- \gamma } (\mu )$ holds, then $(\omega_{i,\phi_i})_i$ is a twisted coupled {\kah}--Einstein metric.
\end{prop}

\begin{proof}

Since the quasi-psh function $\sum_i \phi_i$ has zero Lelong number (\cite[Corollary 1.8]{GZw}), Skoda's integrability theorem (\cite{Sk}) implies that $\mu$ is absolutely continuous
with respect to the reference volume form $dV$.
So, we immediately have
$$
{\rm Ent}_{dV} (\mu) = \int_X \log \frac{\mu}{dV} \mu =  - \log \int_X e^{ - \gamma \sum_i  \varphi_i } dV- \gamma \int_X \sum_i \varphi_i \mu .
$$
Thus, we can compute as follows.
\begin{eqnarray*}
F^{\rm cp}_{- \gamma } (\mu )
&=&- \log \int_X e^{ - \gamma \sum_i  \varphi_i } dV- \gamma \int_X \sum_i \varphi_i \mu- \gamma  E^{\rm cp} (\mu) \\
&\leq & - \log \int_X e^{ - \gamma \sum_i  \varphi_i } dV - \gamma  \int_X \sum_i  \varphi_i\mu + \gamma \int_X \sum_i  \varphi_i\mu - \gamma  \sum_i \mathcal{E}_i (\varphi_i) \\
&=& \gamma D^{\rm cp}_{-\gamma}((\varphi_i)_i).
\end{eqnarray*}
Note that the equality holds iff $\mu = {\rm MA}_{\omega_i} (\phi_i)$ for all $i$, so $(\omega_{i,\phi_i})_i$ is a twisted coupled {\kah}--Einstein metric.
\end{proof}

As a consequence of the above results, we obtain the following.

\begin{cor}
\label{inf}
We have
$$
\inf_{\mu \in \mathcal{M}^1(X)} F^{\rm cp}_{ -\gamma } (\mu ) = \inf_{(\varphi_i)_i \in \prod_{i} \mathcal{E}^1(X,\omega_i) } \gamma D^{\rm cp}_{ - \gamma }((\varphi_i)_i) .
$$
Moreover, the infima of $F_{-\gamma}^{\rm cp}$ and $\gamma D_{-\gamma}^{\rm cp}$ coincide and are attained by the twisted coupled {\kah}--Einstein metric.
\end{cor}

The next lemma converts boundedness below at a larger parameter into coercivity at a smaller parameter.

\begin{lemma}
 \label{Dcoer}
    Assume that $D_{-\gamma^\prime}^{\rm cp}$ is bounded from below.
   Then, for any $0<\gamma < \gamma^\prime$ and  $\alpha \in (0 , \min { \{ \gamma ,  \alpha (X) \} }) $, there exists $C>0$ such that
   $$
   D_{-\gamma}^{\rm cp} \geq \frac{\alpha (\gamma^\prime -\gamma)}{\gamma (\gamma^\prime - \alpha)} J^{\rm cp}  -C.
   $$
In particular, $D_{-\gamma}^{\rm cp}$ is coercive on $\prod_i \mathcal{E}^1 (L_i)$.
Here, the symbol $\alpha(X)$ denotes the alpha invariant
introduced by Tian \cite{Tian}.
\end{lemma}
\begin{proof}
We can prove this lemma by adapting the proof of \cite[Proposition 4.9]{RTZ}. (Note that our hypothesis is slightly different.)
We set $p := \frac{\gamma^\prime - \alpha}{\gamma^\prime - \gamma}$ and $q =  \frac{\gamma^\prime - \alpha }{\gamma - \alpha}$.
Note that $p^{-1} + q^{-1} =1$ and $\frac{\alpha}{\gamma} p^{-1} + \frac{\gamma^\prime }{\gamma} q^{-1} =1$.
Then, we have
\begin{eqnarray*}
- \frac{1}{\gamma} \log \int_X e^{- \gamma \sum_i \varphi_i} dV
&=& - \frac{1}{\gamma} \log \int_X e^{- \alpha p^{-1} \sum_i \varphi_i} e^{-\gamma^\prime q^{-1} \sum_i \varphi_i}  dV \\
&\geq& - \frac{1}{\gamma p} \log \int_X e^{- \alpha \sum_i \varphi_i} dV - \frac{1}{\gamma q} \log \int_X e^{- \gamma^\prime \sum_i \varphi_i} dV \\
&\geq& - \log C_\alpha + \frac{\alpha}{\gamma} p^{-1} \sup \sum_i \varphi_i - \frac{1}{\gamma q} \log \int_X e^{- \gamma^\prime \sum_i \varphi_i} dV. \\
\end{eqnarray*}
Here, we have used the fact that there exists $C_\alpha > 0 $ such that $\int_X e^{-\alpha u} dV < C_\alpha$ for any $u \in {\rm PSH}(X, \omega)$ with $\sup_X u =0$. 
Thus, we have
\begin{eqnarray*}
D_{-\gamma}^{\rm cp} ((\varphi_i)_i )
&\geq& - \log C_\alpha + \frac{\alpha}{\gamma} p^{-1} \left( \sup \sum_i \varphi_i - \sum_i \mathcal{E}_i (\varphi_i) \right) \\
&& \hspace{50pt} - \frac{\gamma^\prime}{\gamma} q^{-1}\mathcal{E}_i (\varphi_i) - \frac{1}{\gamma q} \int_X e^{- \gamma^\prime \sum_i \varphi_i} dV . \\
&=& - \log C_\alpha + \frac{\alpha}{\gamma} p^{-1} \left( \sup \sum_i \varphi_i - \sum_i \mathcal{E}_i (\varphi_i) \right)  + \frac{\gamma^\prime}{\gamma q}  D^{\rm cp}_{-\gamma^\prime} ((\varphi_i)_i ) .
\end{eqnarray*}

Recall that there exists a uniform constant $C>0$ depending only on $X$ and $dV $ such that $\sum_i \sup_X  \phi_i \leq \sum_i \int_X \phi_i dV +C  \leq \sup_X  \sum_i\phi_i +C$ (see \cite[Proposition 2.7]{GZ} and \cite[Lemma 3.45]{Da}).
Since $\sum_i\sup  \varphi_i - \sum_i \mathcal{E}_i (\varphi_i)$ is compatible with the functional $J^{\rm cp}$, we obtain the proof.    
\end{proof}

In summary, we obtain the following.
\begin{thm}
\label{coercivity}
Fix $\gamma >0$.
The following are equivalent.

\medskip

{\rm (i)} $F_{-\gamma}^{\rm cp}$ is coercive on $\mathcal{M}^1 (X)$;

{\rm (ii)} there exists $\epsilon >0$ such that $F_{-\gamma-\epsilon}^{\rm cp}$ is bounded from below on $\mathcal{M}^1 (X)$;

{\rm (iii)} $D_{-\gamma}^{\rm cp}$ is coercive on $\prod_i \mathcal{E}^1 (L_i)$;

{\rm (iv)} there exists $\epsilon >0$ such that $D_{-\gamma-\epsilon}^{\rm cp}$ is bounded from below on $\prod_i \mathcal{E}^1 (L_i)$;

{\rm (v)} $D_{-\gamma}^{\rm cp}$ is coercive on $\Delta $;

{\rm (vi)} there exists $\epsilon >0$ such that $D_{-\gamma-\epsilon}^{\rm cp}$ is bounded from below on $\Delta $;

\end{thm}

\begin{proof}
The equivalence of (i) and (ii) follows from the definition of the coercivity of $F^{\rm cp}_{-\gamma}$.
The equivalence of (ii) and (iv) follows from Corollary \ref{inf}.
(iii) immediately implies (v) by definition and (v) implies (i) by Proposition \ref{F geq G}.
(iv) implies (iii) by Lemma \ref{Dcoer}.
(vi) implies (ii) by Proposition \ref{F geq G} and (iv) implies (vi) by definition.
\end{proof}
By Theorem \ref{Ding criterion}, we immediately obtain the following characterization of coupled {\kah}--Einstein metrics.
\begin{cor}
\label{analytic correspondence}
Assume that ${\rm Aut}(X)$ is finite.
There exists a unique coupled {\kah}--Einstein metric if and only if the coupled free energy $F_{-1}^{\rm cp}$ is coercive.
\end{cor}

\begin{rem}
\label{reduction to sm}
By \cite[Lemma 3.1]{BDL} and \cite[Theorem B]{BJfin}, we can show that in order to check the coercivity of $F_{-1}^{\rm cp} $ on $\mathcal{M}^{1}(X)$ (so, $\mathcal{P}(X)$ by definition), it suffices to check the coercivity of $F_{-1}^{\rm cp} $ with respect to smooth probability measures, i.e., normalized volume forms.
In the next subsection, we will see the non-Archimedean version of this fact.
\end{rem}

We can also characterize the existence of a coupled {\kah}--Einstein metric by the lower semicontinuity of the coupled free energy on the space of probability measures $\mathcal{P}(X)$ with respect to the weak topology. 
The following theorem corresponds to the coupled version of \cite[Theorem 5.9]{Berinv}.

\begin{thm}
   \label{prob}
    Fix $\gamma \in (0,1]$.
    If $\gamma =1$, we assume that ${\rm Aut}(X)$ is finite.
     The following are equivalent.
\begin{enumerate}
    \item  $F_{-\gamma}^{\rm cp}$ has a unique minimizer on $\mathcal{P}(X)$.
    \item There exists $\epsilon >0$ such that $F_{-\gamma-\epsilon}^{\rm cp}$ is bounded from below, i.e., there exists $C =C_\epsilon >0$ such that
$$
F_{-\gamma}^{\rm cp} (\mu)  \geq \epsilon E^{\rm cp} (\mu)  -C ,\,\,\,\, \mu \in \mathcal{P}(X).
$$
\item $F_{-\gamma}^{\rm cp}$ is lower semi-continuous on $\mathcal{P} (X)$ with respect to the weak topology.

\end{enumerate}

\end{thm}

\begin{proof}
((ii)$\implies$(iii))
We take a sequence $\mu_j \to \mu_\infty $ in $\mathcal{P}(X)$.
Recall that ${\rm Ent}_{dV}$ is lower semi-continuous with respect to the weak topology, so it suffices to show that there exists a subsequence whose coupled pluricomplex energy converges to $E^{\rm cp}(\mu_\infty)$.
As we want to show the lower semi-continuity of $F_{-\gamma}^{\rm cp}$, we may assume that there exists a uniform constant $F_{-\gamma}^{\rm cp} (\mu_j) \leq C < \infty$.
Since ${\rm Ent}_{dV} (\mu_j) =\epsilon^{-1} ( (\gamma + \epsilon)F_{-\gamma}^{\rm cp} -  \gamma F_{-\gamma - \epsilon}^{\rm cp}) (\mu_j) < C^\prime$ for some uniform constant $C^\prime$, we can find a subsequence $(\mu_{j_\nu})_\nu$ such that $\mu_{j_\nu} \to \mu_\infty$ weakly and $E^{\rm cp} (\mu_{j_\nu}) \to E^{\rm cp} (\mu_\infty)$ by the compactness theorem of the relative entropy \cite[Theorem 2.17]{BBEGZ}.
Thus, we obtain $F_{-\gamma}^{\rm cp}(\mu_\infty) \leq \liminf F_{-\gamma}^{\rm cp}(\mu_j)$.

((iii)$\implies$(i))
By Prokhorov's theorem $\mathcal{P} (X)$ is compact with respect to the weak topology.
So we can find a minimizer since $F_{-\gamma}^{\rm cp}$ is lsc.
The uniqueness of a minimizer of $F^{\rm cp}_{-\gamma}$ follows from Corollary \ref{inf} and the uniqueness of a minimizer of $D_{-\gamma}^{\rm cp}$.

((i)$\implies$(ii)) 
Recall that a minimizer of $F_{-\gamma}^{\rm cp}$ corresponds to a minimizer of $D_{-\gamma}^{\rm cp}$ (twisted coupled {\kah}--Einstein metric).
In addition, Theorem \ref{Ding criterion} implies that the existence of the twisted  coupled {\kah}--Einstein metric is equivalent to the coercivity of  $D_{-\gamma}^{\rm cp}$ (\cite{HW}, see also \cite[Section 8]{FH}).
Since $F_{-\gamma}^{\rm cp}$ is coercive iff  $D_{-\gamma}^{\rm cp}$ is coercive by Theorem \ref{coercivity}, we have finished the proof.
\end{proof}

\subsection{Non-Archimedean coupled free energy and thermodynamical $K$-stability}

 In this subsection, we generalize the thermodynamical formalism in the non-Archimedean setting by Boucksom--Jonsson \cite{BJKst} to the coupled case when $X$ is a Fano manifold.
 
We recall some definitions in order to discuss our applications of the non-Archimedean pluripotential theory to coupled {\kah}--Einstein metrics (for more details, see \cite{BJKst,BJnaCY, BJGlo,Chi Li}).
Let $X^{\rm val}$ be the set of valuations $v : \mathbb{C} (X)^{*} \to \mathbb{R}  $ on $X$ with $v |_{\mathbb{C}^{*}} = 0$ .
Let $X^{\rm NA}$ be the Berkovich analytification of $X$ with respect to the trivial absolute value.
It is known that $X^{\rm NA}$ is a compactification of $X^{\rm val}$ with respect to the topology of pointwise convergence.
We denote $X^{{\rm div}}$ the set of divisorial valuations which is a (dense) subset of $X^{\rm NA}$.
Any divisorial valuation $v \in X^{{\rm div}} $ is of the form 
$v=\lambda {\rm ord_F} , \, \lambda \in \mathbb{Q}_{\geq0}$ for some prime divisor $F $ over $X$.

\begin{definition}
We say that a flat family of projective varieties $\pi_\mathcal{X}: (\mathcal{X,L} )\to \mathbb{C}$ is a model if the following conditions are satisfied:
\begin{enumerate}
\item $\mathcal{X}$ is normal;
\item  there is an equivariant $\mathbb{C}^{*}$-action on $\pi_\mathcal{X}: (\mathcal{X,L} ) \to \mathbb{C}$ ;
\item  $\mathcal{X}_t := \pi_{\mathcal{X}}^{-1} (t) \simeq (X,L)$ for any $t \neq 0$.
\end{enumerate}
The trivial model of $(X,L)$ is $(X_\mathbb{C}, L_\mathbb{C}):=(X,L) \times \mathbb{C}$ with the trivial $\mathbb{C}^*$-action.
If we ignore the data $L$ and $\mathcal{L}$ , we simply say that $\mathcal{X}$ is a model of $X$.
We say that a model $\mathcal{X}$ is a simple normal crossing model (for short, a SNC model) if a pair $(\mathcal{X,X_{\rm 0}^{\rm red}})$ is a simple normal crossing pair.
We write the set of SNC models for $X$ as ${\rm SNC}(X)$.
\end{definition}

Note that the set of SNC models is a directed set by Hironaka's resolution theorem.
For any SNC model $\mathcal{X}$ of $X$, we can attach the dual complex $\Delta_\mathcal{X}$.
It is known that an element of $\Delta_\mathcal{X}$ corresponds to a quasi-monomial valuation.
There is a retraction map $r_\mathcal{X} : X^{\rm NA} \to \Delta_\mathcal{X}$ and an isomorphism $\displaystyle X^{\rm NA} \simeq \lim_{\longleftarrow} \Delta_\mathcal{X}$ induced by the direct system $(\Delta_\mathcal{X}, r_\mathcal{X})_{\mathcal{X} \in {\rm SNC}(X) }$.

\begin{definition}
    We call $(\mathcal{X,L})$ a semiample (resp. ample) test configuration of $(X,L)$ if $(\mathcal{X,L}) $ is a model of $(X,L)$ and $\mathcal{L}$ is $\pi_\mathcal{X}$-semiample (resp. $\pi_\mathcal{X}$-ample).
\end{definition}

\begin{definition}[Non-Archimedean metric]
For a semi-ample test configuration $(\mathcal{X},\mathcal{L})$ of $(X,L)$ dominating the trivial test configuration, i.e., there exists a $\mathbb{C}^*$-equivariant birational morphism $ p_{\mathcal{X}} : \mathcal{X} \to X \times \mathbb{C}$, we define a smooth and positive non-Archimedean metric $\phi = \phi_{(\mathcal{X},\mathcal{L})}$ by
$$
( \phi_{(\mathcal{X}, \mathcal{L})}   - \phi_{{\rm triv}}) (v):= \sigma(v) ( \mathcal{L} - p_{\mathcal{X}}^* L_{\mathbb{C}}) , \,\,\, v \in X^{\rm NA}.
$$
Here, $\sigma : X^{\rm NA} \to (X \times \mathbb{C})^{\rm NA}$ denotes the Gauss extension and $\phi_{{\rm triv}}$ is the trivial metric associated with the trivial test configuration (see \cite[\S 4.1]{BHJ}).
Note that by subtracting the trivial metric, we can regard $\phi_{(\mathcal{X}, \mathcal{L})}   - \phi_{{\rm triv}}$ as a function on the Berkovich analytification $X^{\rm NA}$.
We write the space of smooth and positive non-Archimedean metrics defined in this way as $\mathcal{H}^{\rm NA} (L)$.
    We define the space of non-Archimedean psh metrics by
    $$
    {\rm PSH}^{\rm NA} (L) := \left\{ \,   \phi : X^{\rm NA
    } \to \mathbb{R} \cup \{- \infty \}  \,  \middle|  \,  \phi \not\equiv -\infty  \,\, {\rm and} \,\,  \exists \{\phi_i\}_i\subset \mathcal{H}^{\rm NA} (L) \,\, {\rm s.t.} \,\,  \phi_i \searrow \phi   \,  \right\}.
    $$
\end{definition}
If $(\mathcal{X,L})$ is a model of $(X,L)$, we can define a function $f_{\mathcal{(X,L)}}$ by
$$
f_{(\mathcal{X}, \mathcal{L})}  (v):= \sigma(v) ( \mathcal{L} - p_{\mathcal{X}}^* L_{\mathbb{C}}),
$$
but in general $f_{\mathcal{(X,L)}} + \phi_{\rm triv} \notin \mathcal{H}^{\rm NA} (L)$.
In this case, we define $\phi_{(\mathcal{X,L})}: = \phi_{\rm triv} + P(f_{(\mathcal{X,L})})$, where
$$
P(f_{(\mathcal{X,L})}) (v)= \sup \left\{  ( \psi  - \phi_{{\rm triv}}) (v)  \, \middle| \,   \psi \in {\rm PSH}^{\rm NA}(L),\, \psi  - \phi_{{\rm triv}} \leq f_{(\mathcal{X,L})}\right\}  .
$$
By \cite[Theorem 8.3]{BFJ}, the function $\phi_{(\mathcal{X,L})}$ defined above is continuous and satisfies $\phi_{(\mathcal{X,L})} \in {\rm PSH}^{\rm NA} (L)$.
We can define 
functionals on the space of smooth and positive non-Archimedean metrics analogous to those on the space of K\"ahler potentials by following \cite{BHJ}.

\begin{definition}[Non-Archimedean Aubin--Mabuchi energy]
For $\phi = \phi_{(\mathcal{X},\mathcal{L})} \in \mathcal{H}^{\rm NA} (L)$, we set the non-Archimedean Aubin--Mabuchi energy by
$$    \mathcal{E}^{\rm NA} (\phi) := \frac{1}{(n+1)V} \bar{\mathcal{L}}^{n+1}.
    $$
    Here, we consider the trivial compactification $(\bar{\mathcal{X}},\bar{\mathcal{L}}) $ of $(\mathcal{X,L})$ over $\mathbb{P}^1$.
\end{definition}

We can extend the non-Archimedean {\ma} energy to ${\rm PSH}^{\rm NA} (L)$.
\begin{definition}{\rm (\cite[\S 6.1]{BJnaCY} and see \cite[\S 1.13]{BJKst})}
    For $\phi \in {\rm PSH}^{\rm NA} (L)$, we set
    $$
    \mathcal{E}^{\rm NA} (\phi) := \inf \left\{ \,    \mathcal{E}^{\rm NA} (\tilde{\phi})          \,     \middle| \,  \phi \leq \tilde{\phi} \in \mathcal{H}^{\rm NA} (L)     \,    \right\} .
    $$
We define the set of non-Archimedean metric with finite energy by
$$
\mathcal{E}^{\rm 1,NA} (L) := \left\{ \, \phi \in {\rm PSH}^{\rm NA} (L) \, \middle|  \, \mathcal{E}^{\rm NA} (\phi) > - \infty \, \right\}
$$
    
\end{definition}

For a non-Archimedean metric with finite energy, we can define the (mixed) {\ma} operator in the non-Archimedean setting.
\begin{thm}{\rm (\cite[Theorem 6.9]{BJnaCY}, see also \cite[Theorem 2.10]{Chi Li})}    There exists a unique operator ${\rm MA}^{\rm NA}$ from $\mathcal{E}^{\rm 1,NA}(L) \times ... \times \mathcal{E}^{\rm 1,NA}(L) $ to the set of Radon measures on $X^{\rm NA}$ such that

\begin{enumerate}

\item for semiample test configurations $(\mathcal{X}, \mathcal{L}_k)$ of $(X,L)$ such that the central fiber $\mathcal{X}_0 = \sum_j b_j E_j$ is Cartier, we have
$$
{\rm MA}^{\rm NA} (\phi_{(\mathcal{X}, \mathcal{L}_1)} , ... , \phi_{(\mathcal{X}, \mathcal{L}_n)} ) = V^{-1} \sum_j b_j ( \mathcal{L}_1 |_{E_j} \cdot ... \cdot \mathcal{L}_n |_{E_j} ) \delta_{x_j}
$$
where $x_j = b_{j}^{-1} r ({\rm ord}_{E_j})$ which is called the Shilov point associated with the irreducible component $E_j \subset \mathcal{X}_0$ and $r  : (X \times \mathbb{C})^{{\rm div}} \to   X^{{\rm div}}$ is the restriction map;

\item for any $\phi, \phi_1, ...., \phi_n \in \mathcal{E}^{\rm 1,NA} (L)$, we have $\int_{X^{\rm NA}} ( \phi - \phi_{{\rm triv}} ) {\rm MA}^{\rm NA}  (\phi_1, ... , \phi_n )   > - \infty$ ;

\item for any decreasing nets $\phi_k^{l} \to \phi_k$ in $\mathcal{E}^{\rm 1,NA}(L)$, we have
$$
\int_{X^{\rm NA}} ( \phi_0^j - \phi_{{\rm triv}} ) {\rm MA}^{\rm NA}  (\phi_1^j, ... , \phi_n^j ) \to \int_{X^{\rm NA}} ( \phi_0 - \phi_{{\rm triv}} ) {\rm MA}^{\rm NA}  (\phi_1, ... , \phi_n ).
$$
\end{enumerate}

\end{thm}

It is known that ${\rm MA}^{\rm NA} (\phi_{\rm triv}) $ is the Dirac mass at the trivial valuation $v_{\rm triv} \in X^{\rm NA}$.
For simplicity, we write ${\rm MA}^{\rm NA} ( \phi) = {\rm MA}^{\rm NA} (\phi,...,\phi)  $.
We will also use the non-Archimedean version of the functionals $I$ and $J$.
\begin{definition}{\rm (\cite[\S 7.2]{BHJ}, see also \cite[\S 1.13]{BJKst})}
For a test configuration $(\mathcal{X}, \mathcal{L})$, we define the following functionals:
$$
I^{\rm NA} (\phi_{(\mathcal{X}, \mathcal{L})}) := V^{-1} ( \bar{\mathcal{L}} \cdot L_{\mathbb{P}^1}^{n} - \bar{\mathcal{L}}^{n+1} + \bar{\mathcal{L}}^{n} \cdot L_{\mathbb{P}^1}) = \int_{X^{\rm NA}} (\phi_{(\mathcal{X}, \mathcal{L})} - \phi_{{\rm triv}} ) ( {\rm MA}^{\rm NA} (\phi_{{\rm triv}} ) -  {\rm MA}^{\rm NA} (\phi_{(\mathcal{X}, \mathcal{L})}))
$$
and
$$
J^{\rm NA} (\phi_{(\mathcal{X}, \mathcal{L})}) := V^{-1} \bar{\mathcal{L}} \cdot L_{\mathbb{P}^1}^{n} - E^{\rm NA} (\phi) = \int_{X^{\rm NA}} (\phi_{(\mathcal{X}, \mathcal{L})}- \phi_{{\rm triv}} ) {\rm MA}^{\rm NA} ( \phi_{{\rm triv}} ) - \mathcal{E}^{\rm NA} (\phi_{(\mathcal{X}, \mathcal{L})}).
$$
    Moreover, we can naturally extend $I^{\rm NA}$ and $J^{\rm NA}$ to $\mathcal{E}^{\rm 1,NA} (L)$ so that they are continuous under decreasing limits (see \cite[\S 7.4]{BJGlo}).
\end{definition}

The non-Archimedean version of the pluricomplex energy can be defined as follows.

\begin{definition}[non-Archimedean pluricomplex energy]
For a Radon probability measure $\mu$ on $X^{\rm NA}$, we define
$$
E^{\rm NA} (\mu) =E^{\rm NA}_{L} (\mu) := \sup_{\phi \in \mathcal{E}^{\rm 1,NA} (L)} \left(  \mathcal{E}^{\rm NA} (\phi) - \int_{X^{\rm NA} } ( \phi  - \phi_{{\rm triv}})    \mu \right) \in \mathbb{R} \cup \{+\infty \}.
$$
A measure $\mu$ is said to have finite energy if $E^{\rm NA} (\mu) < + \infty$.
We denote the set of probability
measures with finite energy as $\mathcal{M}^{\rm 1,NA}(X)$.
We say that a net $\{ \mu_j \}_j \subset \mathcal{M}^{\rm 1,NA} (X)$ converges strongly to $\mu$ if $\mu_j \to \mu$ weakly and $\lim\limits_j E^{\rm NA} (\mu_j) = E^{\rm NA} (\mu)$. 
    \end{definition}

By \cite[Lemma 9.3]{BJGlo}, we have $E^{\rm NA}_{\ell L} (\mu) = \ell E^{\rm NA}_{L} (\mu), \,\mu \in \mathcal{M}^{\rm{1,NA}} (X)$ for any $\ell \in \mathbb{Q}_{>0}.$

\begin{rem}
\label{same set}
Boucksom--Jonsson \cite[Corollary 9.21]{BJGlo} proved that the space $\mathcal{M}^{\rm 1,NA}(X)$ does not depend on the choice of ample line bundle as a non-Archimedean analog of Di Nezza's result \cite[Proposition 3.1]{DN}.
\end{rem}
Moreover, Boucksom--Jonsson proved the following.
\begin{prop}{\rm (\cite[Proposition 10.16]{BJGlo})}
\label{same strong topology}
The strong topology on $\mathcal{M}^{\rm 1,NA} (X)$ is independent of the choice of ample line bundle.
\end{prop}

In order to study the (coupled) {\kah}--Einstein metric from the non-Archimedean viewpoint, we need the following non-Archimedean Calabi--Yau theorem by Boucksom--Jonsson.
\begin{thm}{\rm (\cite[Theorem 7.3]{BJnaCY})}
\label{naCY}
The non-Archimedean {\ma} operator gives a bijection
$$
{\rm MA}^{\rm NA} : \mathcal{E}^{\rm 1,NA} (L)/\mathbb{R} \to \mathcal{M}^{\rm 1,NA} (X).
$$
Moreover, $\mu = {\rm MA}^{\rm NA} (\phi) \in \mathcal{M}^{\rm 1,NA} (X)$ if and only if $E^{\rm NA} (\mu) = \mathcal{E}^{\rm NA} (\phi) - \int_{X^{\rm NA}} ( \phi -  \phi_{{\rm triv}}) \mu$. 
In this case, the equality $E^{\rm NA} (\mu) = (I^{\rm NA} - J^{\rm NA})(\phi) $ holds.
\end{thm}

\begin{rem}
\label{compatibility}
It is well-known that the following inequalities hold on $\mathcal{E}^{\rm 1,NA}(L)$:
$$
\frac{1}{n+1} I^{\rm NA} \leq I^{\rm NA}-J^{\rm NA} \leq \frac{n}{n+1} I^{\rm NA}
$$
(see \cite[Proposition 7.8]{BHJ} and \cite[\S 1.11]{BJKst}).
So, if $\mu = {\rm MA}^{\rm NA} (\phi)$, the non-Archimedean pluricomplex energy is equivalent to the functionals $I^{\rm NA}, J^{\rm NA}$ and $I^{\rm NA} - J^{\rm NA}$.
\end{rem}

From now on, for each line bundle $L_i$, we use the symbols $\mathcal{E}^{\rm  NA}_i, I^{\rm  NA}_i, J^{\rm  NA}_i, {\rm MA}_i^{\rm NA}, E^{\rm NA}_i$.
We define the coupled pluricomplex energy in the non-Archimedean setting as follows.
\begin{definition}[non-Archimedean coupled pluricomplex energy]
The non-Archimedean coupled pluricomplex energy for a Radon probability measure $\mu $ on $X^{\rm NA}$ is defined by
$$
E^{\rm cp, NA}(\mu) := \sum_{i=1}^{m} E_i^{\rm NA}(\mu).
$$
The set of finite energy measures on $X^{\rm NA}$ can also be written as
$$
\mathcal{M}^{\rm 1,NA} (X) = \{ \mu : {\rm Radon \,\, probability \,\, measure \,\, on }\, X^{\rm NA}\,\,  {\rm s.t.} \,\, E^{\rm cp,NA} (\mu) < \infty\}.
$$
Recall that $\mathcal{M}^{\rm 1,NA}(X)$ is independent of the choice of $L_i$ (see Remark \ref{same set}).
\end{definition}

In order to define the free energy in the non-Archimedean setting, we recall the non-Archimedean entropy.
The log discrepancy function 
$$
A_X : X^{{\rm div }} \to \mathbb{R}_{\geq 0}
$$
is defined by
$$
A_X (\lambda {\rm ord }_F ) := \lambda (1 + {\rm ord}_F (K_Y - \pi^* K_X))
\quad\text{and}\quad
A_X(\phi_{\rm triv})=0
$$
where $F \subset Y \xrightarrow{\pi}  X$ is a prime divisor over $X$ and $\lambda \in \mathbb{Q}_{>0}$.
Boucksom--Jonsson proved that there exists a unique maximal lower-semicontinuous extension of $A_X$ to $X^{\rm NA}$  \cite[Theorem 2.1]{BJKst}.

\begin{definition}{\rm (\cite[Definition 7.17]{BHJ}, \cite[\S 2.2]{BJKst})}
For $\mu \in \mathcal{M}^{\rm 1,NA} (X)$, we define the non-Archimedean entropy of $\mu$ by
$$
{\rm Ent}^{\rm NA} (\mu ) := \int_{v \in X^{\rm NA}} A_X (v) \mu .
$$

\end{definition}

We introduce the following natural functional which is analogous to the non-Archimedean free energy in \cite[\S 2.5]{BJKst}.

\begin{definition}[non-Archimedean coupled free energy]
\label{def of NA cp free}
The non-Archimedean coupled free energy is defined by
$$
F^{\rm cp, NA}_{-\gamma } (\mu) := {\rm Ent}^{\rm NA} (\mu) - \gamma E^{\rm cp,NA} (\mu) , \,\,\, \mu \in \mathcal{M}^{\rm 1,NA} (X).
$$

\end{definition}

The space $\mathcal M^{1,\mathrm{NA}}(X)$ and 
the entropy $\mathrm{Ent}^{\mathrm{NA}}$ are independent of the decomposition, 
whereas $E^{\mathrm{cp,NA}}$ depends on the summands $L_i$. 

We can show the following non-Archimedean analog of Proposition \ref{proportional}.
\begin{prop}
\label{NA proportional}
Assume that $L_i$ are proportional to $-K_X$, i.e., there exist $k_i \in \mathbb{Q}_{>0}$ such that $\sum_ik_i = 1$ and $L_i = -k_i K_X$.
Then, we have $F_{-\gamma}^{\rm cp, NA} = {\rm Ent}^{\rm NA} - E^{\rm NA}_{-K_X}$.
In other word, the non-Archimedean coupled free energy corresponds to the non-Archimedean free energy of Boucksom-Jonsson \cite[\S 2.5]{BJKst}.

\end{prop}

\begin{proof}
By \cite[Lemma 9.3]{BJGlo}, we can compute as follows.
$$
F^{\rm cp, NA}_{-\gamma } = {\rm Ent}^{\rm NA} - \gamma \sum_{i=1}^m E^{\rm NA}_{-k_i K_X} = {\rm Ent}^{\rm NA} - \gamma \sum_{i=1}^m k_iE^{\rm NA}_{- K_X}= {\rm Ent}^{\rm NA} - \gamma E^{\rm NA}_{- K_X}.
$$
\end{proof}

As in the complex setting, we compare the coupled free energy with the coupled Ding functional.

\begin{definition}{\rm (\cite{BJKst,FH})} 
The non-Archimedean coupled Ding functional is defined by
$$
    D^{\rm cp,NA}_{-\gamma} ((\phi_i)_i) := \gamma^{-1} L^{\rm NA} (\gamma \sum_i  \phi_i) -\sum_{i=1}^{m} \mathcal{E}_i^{\rm NA} (\phi_i),\,\,\, (\phi_i)_i \in \prod_i \mathcal{E}^{\rm 1,NA} (L_i).
    $$
Here, the non-Archimedean $L$-functional is given by   $$
L^{\rm NA} (\phi) := \inf_{x \in X^{{\rm div}}} \left\{  A_X (x) + (\phi- \phi_{{\rm triv}} )(x) \right\} .
$$
\end{definition}

\begin{rem}
Boucksom--Jonsson proved that the non-Archimedean entropy is the Legendre transform of the non-Archimedean $L$-functional \cite[Proposition 2.8]{BJKst}:
$$
{\rm Ent}^{\rm NA} (\mu) = \sup_{f \in C^0 (X^{\rm NA})} \left( L^{\rm NA}(f) - \int_{X^{\rm NA}} f \mu \right).
$$
\end{rem}

We recall the following.
\begin{prop}{\rm (\cite[Proposition 2.6]{BJKst})}
\label{LNA}
The following equality holds:
 $$
L^{\rm NA} (\phi) = \inf_{x \in X^{{\rm val}},\, A_X(x) < \infty} \left\{  A_X (x) + (\phi- \phi_{{\rm triv}} ) (x) \right\} .
$$

\end{prop}

In the coupled setting, we can show the following non-Archimedean counterpart of Proposition \ref{F geq G}.
The proof follows the argument of \cite[Theorem 2.10]{BJKst} (see also \cite[Proposition 7.32]{BHJ}).

\begin{prop}
\label{NA ineq}
If $\mu = {\rm MA}_i^{\rm NA} (\phi_i)$ for all $i$, we have
$$
F_{-\gamma}^{\rm cp, NA} (\mu) \geq \gamma D_{-\gamma}^{\rm cp,NA} ((\phi_i)_i).
$$

\end{prop}
\begin{proof}
By the non-Archimedean Calabi--Yau theorem (Theorem \ref{naCY}), we have $E^{\rm cp,NA} (\mu) = \sum_i \left(  \mathcal{E}^{\rm NA}_i (\phi_i) - \int_{X^{\rm NA}}(\phi_i - \phi_{{\rm triv}} ) \mu\right)$.
Thus, we can compute as follows:
\begin{eqnarray}
F_{-\gamma}^{\rm cp,NA} (\mu) 
&=& \int_{X^{\rm NA}} A_X (x) \mu - \gamma \sum_i\mathcal{E}^{\rm NA}_i (\phi_i) + \gamma \int_{X^{\rm NA}} \sum_i( \phi_i - \phi_{{\rm triv}} ) \mu \\
&\geq& \inf_{ x \in X^{{\rm val}}, \, A_X (x) < \infty} \left( A_X (x)+ \gamma \sum_i (\phi_i  - \phi_{{\rm triv}} )(x) \right) - \gamma \sum_i\mathcal{E}^{\rm NA}_i (\phi_i) \\
&=& \gamma D_{-\gamma}^{\rm cp,NA} ((\phi_i)_i).
\end{eqnarray}
Here, we have used Proposition \ref{LNA} in the last inequality.
\end{proof}
We write the log canonical divisor for $\pi : \mathcal{X} \to \mathbb{C} $ as
$$
K_{\mathcal{X}/ \mathbb{C}}^{\rm log} = K_{\mathcal{X}} - \pi^* K_{\mathbb{C}} -\mathcal{X}_0 + \mathcal{X}_0^{\rm red}.
$$
The following analogue of Proposition \ref{converse ineq}
is motivated by \cite[Definition 7.31, Proposition 7.32]{BHJ}.
\begin{thm}
\label{NA converse}
Assume that 
$$
\mathcal{L}= \mathcal{L}_1 + ... + \mathcal{L}_m = -K_{\mathcal{X}/\mathbb{C}}^{\log} + c \mathcal{X}_0
$$
for some $c\in \mathbb{R}$ and $(\mathcal{X}, \mathcal{L}_i)$ are normal semiample test configurations.
We also assume that $\mu \in \mathcal{M}^{\rm 1,NA} (X)$ is written as
$$
\mu = \sum_{E} c_E \delta_{v_E}
$$
where $E$ runs over the irreducible components of $\mathcal{X}_0$.
If the pair $(\mathcal{X}, \mathcal{X}_{0}^{{\rm red}})$ is log canonical, i.e., the log discrepancy of $(\mathcal{X}, \mathcal{X}_{0}^{{\rm red}})$ satisfies $A_{(\mathcal{X}, \mathcal{X}_{0}^{{\rm red}})} \geq 0$ on $X^{\rm div}$, we have
\begin{equation}
\label{NA anticanoncial}
    F_{-1}^{\rm cp ,NA} (\mu) \leq D_{-1}^{\rm cp ,NA} ((\phi_{(\mathcal{X},\mathcal{L}_i })_i).
\end{equation}
The equality holds iff
$$
\mu = {\rm MA}_i^{\rm NA} (\phi_{(\mathcal{X},\mathcal{L}_i)}  ), \, \forall i = 1, ...,m.
$$
\end{thm}
Note that if $(\mathcal{X}, \mathcal{X}_{0})$ is log smooth, i.e., $\mathcal{X}$ is smooth and $ \mathcal{X}_{0}^{ \rm{red}}$ is simply normal crossing, then $(\mathcal{X}, \mathcal{X}_{0}^{{\rm red}})$ is automatically log canonical (see \cite[\S 1.5 and 7.7]{BHJ}).
\begin{rem}
If $m=1$, Proposition 7.32 of \cite{BHJ} implies that the equality of (\ref{NA anticanoncial}) holds by the definition of the non-Archimedean Mabuchi functional.    
\end{rem}
\begin{proof}
For simplicity, we write $
\phi_i  = \phi_{(\mathcal{X,L}_i)}.$
Then, we have $\phi_{(\mathcal{X,L})} = \sum_i \phi_i\in \mathcal{H}^{\rm NA} (-K_X)$. 
One can compute as follows:
\begin{eqnarray*}
    F_{-1}^{\rm cp ,NA} (\mu) 
    &\leq& \int_{X^{ \rm NA}} A_X \mu - \sum_i \left( \mathcal{E}_i^{ \rm NA} (\phi_i) + \int_{X^{\rm NA}} ( \phi_i - \phi_{{\rm triv}})\mu\right) \\
    &=& \int_{X^{ \rm NA}} ( A_X + \sum_i( \phi_i - \phi_{{\rm triv}}) )\mu - \sum_i \mathcal{E}_i^{ \rm NA} (\phi_i) \\
    &=&  \sum_E c_E (A_X (v_E) + ( \phi_{(\mathcal{X,L})}- \phi_{{\rm triv}})(v_E) ) - \sum_i \mathcal{E}_i^{ \rm NA} (\phi_i) \\
    &=& L^{\rm NA} (\sum_i \phi_i)  - \sum_i \mathcal{E}_i^{ \rm NA} (\phi_i)  = D_{-1}^{\rm cp, NA} ((\phi_i)_i).
\end{eqnarray*}
Here, we have used the fact that $A_X (v_E) + (\phi_{(\mathcal{X,L})} - \phi_{\rm triv})(v_E) $ is constant for any $E$ \cite[Lemma 7.30]{BHJ} in the second equality and  \cite[Proposition 7.29]{BHJ} in the third equality.
Note that the equality of the first inequality holds if and only if ${\rm MA}_i^{\rm NA} (\phi_{(\mathcal{X},\mathcal{L}_i)}  ) =  \mu $ for all $i =1,..,m$ by Theorem \ref{naCY}.
\end{proof}

\begin{rem}
By the above theorem, we can consider the coupled version of a (weak) non-Archimedean {\kah}--Einstein metric which was mentioned in \cite[Remark 7.33]{BHJ} as follows:
    $$
\mathcal{L}_1 + ... + \mathcal{L}_m = -K_{\mathcal{X}/\mathbb{C}}^{\log} + c \mathcal{X}_0, \,\,\,\, {\rm MA}_i^{\rm NA} (\phi_{(\mathcal{X},\mathcal{L}_i)}  ) =  {\rm MA}_j^{\rm NA} (\phi_{(\mathcal{X},\mathcal{L}_i)}  ), \, \forall i,j = 1,...,m.
    $$
By Theorem \ref{NA converse}, such a metric satisfies $F_{-1}^{\rm cp,NA} (\mu) = D_{-1}^{\rm cp,NA} ((\phi_i)_i)$.
\end{rem}

By following the proof of \cite[Theorem 2.10]{BJKst}, we can show the non-Archimedean version of Corollary \ref{inf}.

\begin{thm}
\label{NAbdd}
The following are equivalent.

{\rm (i)} $F_{-\gamma}^{\rm cp,NA} \geq 0$ on $\mathcal{M}^{\rm 1,NA}(X)$;

{\rm (ii)} $D_{-\gamma}^{\rm cp,NA} \geq 0$ on $\prod_{i} \mathcal{E}^{\rm 1,NA}(L_i)$.
\end{thm}

\begin{rem}
By \cite[Lemma 2.9]{BJKst} and the proof of \cite[Theorem 8.10]{FH}, we can show that $D_{-\gamma}^{\rm cp,NA} \geq 0$ on $\prod_{i} \mathcal{E}^{\rm 1,NA}(L_i)$ iff $D_{-\gamma}^{\rm cp,NA} \geq 0$ on $\prod_{i} \mathcal{H}^{\rm NA}(L_i)$.
\end{rem}

\begin{proof}

It immediately follows from Proposition \ref{NA ineq} that (ii) implies (i).

We show that (i) implies (ii) by following the proof of \cite[Theorem 2.10]{BJKst}.
For $x \in X^{{\rm div}}$, we have $A_X (x) - \gamma \sum_i E_i (\delta_x) \geq 0$  by taking $\mu = \delta_x$.
So, for any $\phi_i \in  \mathcal{E}^{\rm 1,NA}(L_i)$, we have $A_X (x) +\gamma \sum_i (\phi_i- \phi_{{\rm triv}} ) (x) \geq \gamma \sum_i \mathcal{E}_i(\phi_i)$.
By taking the infimum over $x \in X^{{\rm div}}$, we have $L^{\rm NA}(\gamma\sum_i \phi_i) \geq \gamma  \sum_i \mathcal{E}^{\rm NA}_i(\phi_i)$, so we obtain $D_{- \gamma}^{\rm cp,NA} ((\phi_i)_i) \geq 0$ for any $((\phi_i)_i) \in \prod_i \mathcal{E}^{\rm 1,NA}(L_i)$.
\end{proof}


We can show the non-Archimedean analog of Lemma \ref{Dcoer}.

\begin{lemma}
\label{bdd and coer}
Fix $\gamma < \gamma^\prime$.
Assume that the non-Archimedean coupled Ding functional satisfies $D_{-\gamma^\prime}^{\rm cp,NA} \geq 0$ on $\prod_i \mathcal{E}^{\rm 1,NA} (L_i)$.
Then, for any $\gamma < \gamma^\prime$ and  $\alpha \in (0 , \min { \{ \gamma ,  \alpha (X) \} }) $, we have
$$
D_{-\gamma}^{\rm cp,NA} \geq \frac{\alpha(\gamma^\prime -\gamma)}{\gamma (\gamma^\prime - \alpha)} J^{\rm cp,NA}
$$
on $\prod_i \mathcal{E}^{\rm 1,NA} (L_i)$.
Here, we set $J^{\rm cp, NA}((\phi_i)_i) := \sum_i J_i^{\rm cp, NA} (\phi_i)$.
\end{lemma}

In order to prove this proposition, we need the following.

\begin{definition}{\rm (\cite[\S 5.1]{BJKst})} The Alexander--Taylor capacity of a point $x \in X^{\rm NA}$ with respect to $L$ is defined by
$$
T_L(x) := \sup \{  \sup_{X^{\rm NA}}  (\phi - \phi_{{\rm triv}} ) - (\phi- \phi_{{\rm triv}} ) (x) \, | \, \phi \in {\rm PSH}^{\rm NA} (L)\}.
$$
\end{definition}
Blum--Jonsson \cite[Theorem E]{BlJ} (see also \cite[Theorem 5.15]{BJKst}) proved that the alpha invariant of $L$ is written as
$$
\alpha (L) = \inf_{x \in X^{\rm div}} \frac{A_X (x) }{T_L (x)}.
$$

\begin{proof}
We prove Lemma \ref{bdd and coer} by applying the proof of \cite[Proposition 4.9]{RTZ}.
We take $\alpha \in (0 ,  \min { \{ \gamma ,  \alpha (X) \} }) $.
We set $p := \frac{\gamma^\prime - \alpha}{\gamma^\prime - \gamma}$ and $q =  \frac{\gamma^\prime - \alpha }{\gamma - \alpha}$.
Note that $p^{-1} + q^{-1} =1$ and $\frac{\alpha}{\gamma} p^{-1} + \frac{\gamma^\prime }{\gamma} q^{-1} =1$.
For $x \in X^{\rm div}$, we have
\begin{eqnarray*}
&&\gamma^{-1}A_X (x) + \sum_i (\phi_i- \phi_{{\rm triv}} )  (x) - \sum_i \mathcal{E}^{\rm NA}_i (\phi_i) \\
&=& \frac{\gamma^\prime}{q\gamma} \left( (\gamma^{\prime})^{-1} A_X(x)  + \sum_i (\phi_i- \phi_{{\rm triv}} ) (x)  - \sum_i \mathcal{E}_i^{\rm NA} (\phi_i)\right) \\
&&  \,\,\,\,\,\,  + \frac{1}{\gamma p} \left( A_X (x) + \alpha (\sum_i   (\phi_i- \phi_{{\rm triv}} ) (x) - \mathcal{E}_i^{\rm NA} (\phi_i) ) \right) \\
&\geq& \frac{\gamma^\prime}{q\gamma} D_{-\gamma^\prime}^{\rm  cp,NA} ((\phi_i)_i)  + \frac{1}{\gamma p} \left( \alpha T_L (x)  + \alpha (\sum_i   (\phi_i- \phi_{{\rm triv}} )(x)  - \mathcal{E}_i^{\rm NA} (\phi_i) )\right) \\
&\geq& \frac{\alpha}{\gamma p} \sum_i \left( 
(\sup_{X^{\rm NA}}( \phi_i- \phi_{{\rm triv}} )  -\mathcal{E}_i^{\rm NA} (\phi_i) )\right)\\
&\geq&  \frac{\alpha }{\gamma p} J^{\rm cp,NA}((\phi_i)_i).
\end{eqnarray*}
Note that we have used the maximum principle (see \cite[\S 6.2]{BBJ}), i.e., $\sup_{X^{\rm NA} } \sum_i \phi_i  = \sum_i \phi_i (v_{{\rm triv}})= \sum_i \sup_{X^{\rm NA}} \phi_i$.
By taking the infimum over $x \in X^{\rm div}$, we obtain the proof.
\end{proof}

As an analog of Definition \ref{diagonal}, we introduce the notion of the diagonal in the product space $\prod_i \mathcal{E}^{\rm 1,NA} (L_i)$.

\begin{definition}
\label{diagonal^NA}
We define the non-Archimedean diagonal for $\{L_i \}_i$ by
$$
\Delta^{\rm NA} = \Delta_{\{L_i\}_i }^{\rm NA} := \left\{  ( \phi_i)_i \in \prod_i \mathcal{E}^{\rm 1,NA} (L_i) \, \middle|  \, \,  {\rm MA}_i^{\rm NA} (\phi_i)= {\rm MA}_j^{\rm NA}(\phi_j) , \,\forall i,j=1,...,m \right\}
$$
\end{definition}

By the above results, we can prove the non-Archimedean version of Theorem \ref{coercivity}.

\begin{thm}
\label{NAcoercivity}
Fix $\gamma >0$.
The following are equivalent.

\medskip

{\rm (i)} there exists $\epsilon >0$ such that $F_{-\gamma}^{\rm cp,NA} \geq \epsilon E^{\rm cp,NA}$ on $\mathcal{M}^{\rm 1,NA} (X)$;

{\rm (ii)} there exists $\epsilon >0$ such that $F_{-\gamma-\epsilon}^{\rm cp,NA} \geq 0$ on $\mathcal{M}^{\rm 1,NA} (X)$;

{\rm (iii)} there exists $\epsilon >0$ such that $D_{-\gamma}^{\rm cp,NA} \geq \epsilon J^{\rm cp,NA}$ on $\prod_i \mathcal{E}^{\rm 1,NA} (L_i)$;

{\rm (iv)} there exists $\epsilon >0$ such that $D_{-\gamma-\epsilon}^{\rm cp,NA} \geq 0$ on $\prod_i \mathcal{E}^{\rm 1,NA} (L_i)$;

{\rm (v)} there exists $\epsilon >0$ such that $D_{-\gamma}^{\rm cp,NA} \geq \epsilon J^{\rm cp,NA}$ on $\Delta^{\rm NA}$;

{\rm (vi)} there exists $\epsilon >0$ such that $D_{-\gamma-\epsilon}^{\rm cp,NA} \geq 0$ on $\Delta^{\rm NA}$.

\end{thm}

\begin{proof}
(i) is equivalent to (ii) by definition.
(ii) is equivalent to (iv) by Theorem \ref{NAbdd}.
(iii) implies (v) by definition and (v) implies (i) by the compatibility of $E^{\rm cp,NA}$ and $J^{\rm cp,NA}$ and Proposition \ref{NA ineq}.
(iv) implies (iii) by Lemma \ref{bdd and coer}.
(iv) implies (vi) by definition and (vi) implies (ii) by Proposition \ref{NA ineq}.
\end{proof}



    

In the {\kah}--Einstein case, i.e., $m=1$, the non-Archimedean free energy and the non-Archimedean Mabuchi K-energy play a very important role in the study of (uniform) $K$-stability.
In the coupled case, we can naturally define the following stability notion by using the non-Archimedean coupled free energy.
This is a generalization of the stability condition of \cite{BJKst}.

\begin{definition}
\label{TK}
We say that $(X, \{ L_i \}_i)$ is thermodynamically $K$-stable if there exists $\epsilon >0$ such that  $F^{\rm cp, NA}_{-1-\epsilon} \geq 0$ on $\mathcal{M}^{\rm 1,NA}(X)$, equivalently, 
\begin{equation}
    \label{thermo}
F^{\rm cp, NA}_{-1} (\mu) \geq \epsilon E^{\rm cp,NA} (\mu), \,\,\,\, \forall \mu \in \mathcal{M}^{\rm 1,NA}(X).
\end{equation}
\end{definition}

\begin{rem}
If $m=1$, i.e., the {\kah}--Einstein case, thermodynamical $K$-stability is exactly equivalent to uniform $K$-stability by the thermodynamical formalism in the non-Archimedean setting by Boucksom--Jonsson \cite[Theorem 2.13]{BJKst}.
We use the term thermodynamical K-stability to distinguish condition
\eqref{thermo}, formulated in terms of finite-energy measures, 
from the notion of coupled K-(poly)stability defined using tuples of test configurations (with the same base space) in Hultgren--Witt Nystr\"{o}m's
paper \cite{HW}.

\end{rem}

We can characterize the existence of a coupled {\kah}--Einstein metric by thermodynamical $K$-stability.

\begin{thm}
\label{TK and coer}
The coupled free energy $F_{-1}^{\rm cp}$ is coercive if and only if $(X, \{ L_i \}_i)$ is thermodynamically $K$-stable.
\end{thm}

\begin{rem}
We can prove Theorem \ref{TK and coer} by combining \cite{FH} and our results (see the diagram in Introduction).
Here, we give a proof that does not use the uniform coupled Ding stability.

\end{rem}

\begin{proof}

Assume that $(X, \{L_i\}_i)$ is thermodynamically $K$-stable.
We take sufficiently divisible $r \in \mathbb{Z}_{>0}$ such that $r L_i$ are line bundles for any $i$.
By \cite[Lemma 7.4]{BBJ} and \cite[Lemma 9.3]{BJGlo}, for any divisorial valuation $v_E$ of a prime divisor $E$ over $X$, we have
$$
A_X (E) = {\rm Ent}^{\rm NA} (\delta_{v_E}) \geq \frac{1+\epsilon}{r}\sum_i E^{\rm NA}_{rL_i} (\delta_{v_E}) = \frac{1+\epsilon}{r} \sum_i S_{rL_i} (E)= (1+\epsilon) \sum_i S_{L_i} (E),
$$
where the symbol $S_{L_i}$ denotes the expected vanishing order of $L_i$ (for instance, see \cite[Section 7]{BBJ} or \S 4.3 in this paper).
It follows that the coupled version of the delta-invariant introduced by 
Rubinstein-Tian-Zhang \cite{RTZ} is greater than 1,
and, by Kewei Zhang's result \cite{KZ}, the coupled Ding functional $D^{\rm cp}_{-1}$ is coercive.
Theorem \ref{coercivity} then shows the coercivity of
the coupled free energy $F^{\rm cp}_{-1}$.

We show that the coercivity of  $F_{-1}^{\rm cp}$ implies the existence of  $\epsilon >0$ such that $F_{-1}^{\rm cp,NA} \geq \epsilon E^{\rm cp,NA}$ on $\mathcal{M}^{\rm 1,NA} (X )$.
Note that the coercivity of $F_{-1}^{\rm cp}$ implies that $D_{-1}^{\rm cp}$ is also coercive on $\prod_i \mathcal{E}^1(L_i)$ by Theorem \ref{coercivity}.
For any $\mu \in \mathcal{M}^{\rm 1,NA} (X)$, Theorem \ref{naCY} implies that there exist $\phi_i \in \mathcal{E}^{\rm 1,NA} (L_i)$ such that $\mu = {\rm MA}^{\rm NA}_i (\phi_i)$ for any $i$.
(We note that $\sum_i \phi_i \in \mathcal{E}^{\rm 1,NA} (-K_X)$ by \cite[Proposition 5.3]{FH}.)
By \cite[Theorem 6.6]{BBJ}, there exists a unique maximal geodesic ray $\Phi_i (t) \in \mathcal{E}^1(L_i)$ (emanating from $0$) such that the corresponding non-Archimedean metric is $\phi_i$.
Then, we have $\mathcal{E}_i^{\rm NA} (\phi_i) = \lim\limits_{t \to \infty} t^{-1} \mathcal{E}_i (\Phi_i(t))$ by \cite[Corollary 6.7]{BBJ}.
Moreover, we have $J^{\rm NA}_i (\phi_i) = \lim\limits_{t \to \infty} t^{-1} J_i (\Phi_i(t))$ for each $i$ by \cite[Lemma 2.41]{Chi Li}.
The sum of maximal geodesic rays $\sum_i \Phi_i (t)$ defines a psh ray in $\mathcal{E}^1 (-K_X)$ with linear growth (see Section 8 in \cite{FH}) since any psh geodesic ray has linear growth, i.e., $ \sup_X \Phi_i(t) = O(t)$ for each $i$ by \cite[Proposition 4.1]{BBJ}. 
Thus, we obtain $L^{\rm NA} (\sum_i \phi_i ) = \lim\limits_{t \to \infty} t^{-1} L(\sum_i \Phi_i (t))$ from \cite[Theorem 5.4]{BBJ}, where we write $$L(\phi) = - \log \int_X e^{-\phi} dV.$$
Then, by Proposition \ref{NA ineq}, we can compute as follows.
\begin{eqnarray*}
F_{-1}^{\rm cp,NA}(\mu)
&\geq & D_{-1}^{\rm cp,NA} ( (\phi_i)_i) \\
&=& L^{\rm NA} (\sum_i \phi_i) - \sum_i \mathcal{E}^{\rm NA}_i (\phi_i)  \\
&=& \lim_{t \to \infty} \frac{1}{t} \left( L \left( \sum_{i} \Phi_i(t) \right)  - \sum_i \mathcal{E}_i (\Phi_i (t)) \right) \\
&=& \lim_{t \to \infty} \frac{1}{t} \left(  D_{-1}^{\rm cp} ((\Phi_i(t))_i )\right) \\
&\geq &  \lim_{t \to \infty }  \frac{1}{t}\epsilon J^{\rm cp}((\Phi_i(t))_i ) \\
&=&  \epsilon \sum_i J_i^{\rm NA} (\phi_i) \geq  \frac{\epsilon}{n} \sum_i E_i^{\rm NA} ({\rm MA}_i^{\rm NA}(\phi_i)) =  \frac{\epsilon}{n} E^{\rm cp,NA} (\mu) .
\end{eqnarray*}
Here, we have used the compatibility of the pluricomplex energy and the $J$-functional (see Remark \ref{compatibility}).
This inequality implies that $(X, \{L_i \}_i) $ is thermodynamically $K$-stable.
\end{proof}

We next show that thermodynamical K-stability can be tested on 
Monge–Amp\`ere measures arising from test configurations for
a single polarization.
The proof is based on the 
regularization process of ${\rm Ent}^{\rm NA}$ by Chi Li, Boucksom--Jonsson and Trusiani.

\begin{thm}
\label{one para}
Fix $\e>0$.
In order to check that $(X,\{L_i\}_i )$ is thermodynamically $K$-stable, it suffices to show that the inequality $F^{\rm cp, NA}_{-1-\epsilon}  (\mu)\geq 0$ holds for any $\mu = {\rm MA}^{\rm NA}_{-K_X} (\phi_{(\mathcal{X,L})}) \in \mathcal{M}^{\rm 1,NA} (X)$, 
where $(\mathcal{X,L})$ is a test configuration of $(X,-K_X)$.
In other words, $F^{\rm cp, NA}_{-1-\epsilon} \geq 0$ on $\mathcal{M}^{\rm 1,NA} (X)$ if and only if  $F^{\rm cp, NA}_{-1-\epsilon} \circ {\rm MA}_{-K_X}^{\rm NA} \geq 0$ on $\mathcal{H}^{\rm NA}(-K_X)$. 
\end{thm}

\begin{rem}
Actually, we can show the above theorem for any ample line bundle which is not necessarily the anticanonical line bundle $-K_X$.
\end{rem}


The proof of Theorem \ref{one para} follows from the regularization process of \cite{BJKst,BJnaCY,Chi Li,Tru} and Proposition \ref{same strong topology}.
Firstly, we extend Chi Li's work \cite[Step 1 in Proposition 6.3]{Chi Li} to the coupled case.
\begin{prop}
\label{dual complex}
In order to check that $(X,\{Li\}_i )$ is thermodynamically $K$-stable, it suffices to show that the inequality $F^{\rm cp, NA}_{-1-\epsilon}  (\mu)\geq 0$ holds for any $\mu \in \mathcal{M}^{\rm 1,NA} (X)$ supported on the dual complex $\Delta_{\mathcal{X}}$ for some SNC model $\mathcal{X}$.
\end{prop}
\begin{rem}
We can prove the above proposition by combining Step 1 of \cite[Proposition 6.3]{Chi Li} and Proposition \ref{same strong topology}.
In this paper, we give a (direct) proof by modifying the proof of \cite{Chi Li} for our coupled setting.
\end{rem}
\begin{proof}
It is enough to show that for any $\mu \in \mathcal{M}^{\rm 1,NA} (X )$, we can always find a sequence
$\{\mu_j\}_j \subset \mathcal{M}^{\rm 1,NA} (X)$, each supported on the dual complex $\Delta_{\mathcal{X}_j}$ for some SNC model $\mathcal{X}_j$ such that $F_{-1-\epsilon}^{\rm cp,NA}(\mu_j) \to F_{-1-\epsilon}^{\rm  cp,NA}(\mu)  $ as $j \to \infty$.
We follow the proof of Step 1 in Proposition 6.3 of \cite{Chi Li}.
By \cite[Corollary 7.20]{BJnaCY}, the net $\{ E_i^{\rm NA} ( (r_{\mathcal{X}})_* \mu) \}_{\mathcal{X} \in {\rm SNC(X)}}$ is increasing and $E_i^{\rm NA} (\mu) = \lim_\mathcal{X} E_i^{\rm NA} ( (r_{\mathcal{X}})_* \mu) $ for each $i=1,...,m$.
Here, $r_\mathcal{X} : X^{\rm NA} \to \Delta_\mathcal{X}$ denotes the retraction map associated with $\mathcal{X}$.
For each $i$, we pick sequences $\{ \mathcal{X}^i_j\}_j$ such that $  E_i^{\rm NA} ( (r_{\mathcal{X}_i^j})_* \mu) \nearrow E_i^{\rm NA} (\mu)$.
We can take a SNC model $\mathcal{X}_j$ such that $\mathcal{X}_j \to \mathcal{X}_j^i$ for all $i =1 ,..., m$.
Since $E^{\rm NA}_i ((r_{\mathcal{X}^j})_* \mu) \geq  E^{\rm NA}_i ((r_{\mathcal{X}_i^j})_* \mu) $, we have $E^{\rm NA}_i (\mu) = \lim\limits_j E^{\rm NA}_i ((r_{\mathcal{X}^j})_* \mu) $ for all $i$.
In particular, we have
$$
E^{\rm cp,NA} (\mu) =\lim_j E^{\rm cp,NA} ((r_{\mathcal{X}^j})_* \mu) .
$$
On the other hand, by \cite[Lemma 2.3]{BJKst}, the net $\{ {\rm Ent}^{\rm NA} ( (r_{\mathcal{Y}})_* \mu) \}_{\mathcal{Y} \in {\rm SNC}(X)}$ is increasing and ${\rm Ent}^{\rm NA} (\mu) = \lim_\mathcal{Y} {\rm Ent}^{\rm NA} ( (r_{\mathcal{Y}})_* \mu) $.
We take a sequence $\mathcal{Y}_j$ such that ${\rm Ent}^{\rm NA} (\mu) = \lim\limits_j {\rm Ent}^{\rm NA} ( (r_{\mathcal{Y}_j})_* \mu) $.
We take a SNC model $\mathcal{Z}_j$ such that $\mathcal{Z}_j \to \mathcal{X}_j$ and $\mathcal{Z}_j \to \mathcal{Y}_j$.
Since both functionals $E^{\rm cp,NA}$ and ${\rm Ent}^{\rm NA}$ are increasing, we can easily show that
$$
E^{\rm cp,NA} (\mu) =\lim_j E^{\rm cp,NA} ((r_{\mathcal{Z}_j})_* \mu) \, \, \,{\rm and} \,\,\, {\rm  Ent}^{\rm NA} (\mu) =\lim_j {\rm Ent}^{\rm NA} ((r_{\mathcal{Z}_j})_* \mu) .
$$
Thus, we obtain $F_{-1-\epsilon}^{\rm cp,NA} (\mu) = \lim\limits_j F_{-1-\epsilon}^{\rm cp,NA} ((r_{\mathcal{Z}_j})_* \mu)$.
\end{proof}

We recall the following in order to improve Proposition \ref{dual complex}.
\begin{prop}{\rm (\cite[Step 2 in Proposition 6.3]{Chi Li})}
\label{model}
For any finite energy measure $\mu ={\rm MA}^{\rm NA}_{-K_X} (\phi)\in \mathcal{M}^{\rm 1,NA} (X)$ supported on the dual complex $\Delta_\mathcal{X}$ for some SNC model $\mathcal{X}$ for $(X,-K_X)$, we can find a sequence of models $(\mathcal{X}_j, \mathcal{L}_j)$ such that the associated non-Archimedean metrics $\phi_j = \phi_{(\mathcal{X}_j, \mathcal{L}_j)}$ satisfy $E_i^{\rm NA} ({\rm MA}_{-K_X}^{\rm NA}(\phi_j)) \to E_{i}^{\rm NA} ({\rm MA}_{-K_X}^{\rm NA}(\phi))$ for any $i$ and ${\rm Ent}^{\rm NA} ({\rm MA}^{\rm NA}_{-K_X} (\phi_j)) \to {\rm Ent}^{\rm NA} ({\rm MA}^{\rm NA}_{-K_X} (\phi)) $.
\end{prop}
\begin{rem}We can immediately prove this proposition by combining Step 2 of \cite[Proposition 6.3]{Chi Li} and Proposition \ref{same strong topology}.
Note that we can directly prove this proposition by using only the proof of Step 2 of \cite[Proposition 6.3]{Chi Li} because its proof of the strong convergence of the approximation measures can be applied to the pluricomplex energy for any ample line bundle.
\end{rem}

Recently, Trusiani gave a proof of  the uniform Yau--Tian--Donaldson conjecture \cite{Tru} for (weighted) constant scalar curvature {\kah} metrics.
The important ingredient of his result is that he established the special Fujita approximation and solved Boucksom--Jonsson's regularization conjecture of the non-Archimedean entropy.
By using Proposition \ref{same strong topology}, we have
\begin{thm}{\rm (\cite{Tru}, see also \cite{Chi LiFj})}\label{test} 
For any finite energy measure $\mu = {\rm MA}^{\rm NA}_{-K_X} (\phi)\in \mathcal{M}^{1,NA} (X)$ where $\phi = \phi_{\mathcal{(X,L)}}$ is defined by a model $(\mathcal{X}, \mathcal{L})$, we can find a sequence of test configurations $(\mathcal{X}_j, \mathcal{L}_j)$ for $(X,-K_X)$ such that the associated non-Archimedean metrics $\phi_j = \phi_{(\mathcal{X}_j, \mathcal{L}_j)}$ satisfy $E_i^{\rm NA} ({\rm MA}_{-K_X}^{\rm NA}(\phi_j)) \to E_{i}^{\rm NA} ({\rm MA}_{-K_X}^{\rm NA}(\phi))$ for any $i$ and ${\rm Ent}^{\rm NA} ({\rm MA}^{\rm NA}_{-K_X} (\phi_j)) \to {\rm Ent}^{\rm NA} ({\rm MA}^{\rm NA}_{-K_X} (\phi)) $ as $ j \to \infty$.
\end{thm}

{\it Proof of Theorem \ref{one para}.}
Thanks to Proposition \ref{dual complex}, Proposition \ref{model} and Theorem \ref{test}, 
for any $\mu \in \mathcal{M}^{\rm 1,NA} (X)$, we can find a sequence of measures $\mu_j$ associated with a test configuration for $(X,-K_X)$ such that
$$
 \lim_j E_{i}^{\rm NA} (\mu_j) = E_{i}^{\rm NA} (\mu) \,\,{\rm and} \,\, \lim_j {\rm Ent}^{\rm NA} (\mu_j) = {\rm Ent}^{\rm NA} (\mu) .
$$
As a consequence, we obtain 
$$
  F_{-1-\epsilon}^{\rm cp,NA}  (\mu) = \lim_j F_{-1-\epsilon}^{\rm cp,NA} (\mu_j) \geq 0.
$$\sq

By Theorems \ref{TK and coer} and \ref{one para} 
and Corollary \ref{analytic correspondence}, we can show the following.
In particular, we can characterize the existence of a coupled {\kah}--Einstein metric by test configurations for a single polarization.

\begin{cor}
\label{alg}
Assume that ${\rm Aut}(X)$ is finite.
The following are equivalent.
\begin{itemize}

\item There exists a unique coupled {\kah}--Einstein metric;
\item $(X,\{L_i\}_i)$ is thermodynamically K-stable;
\item there exists $\epsilon>0$ such that $F^{\rm cp, NA}_{-1-\epsilon}  (\mu)\geq  0$ for any $\mu \in \mathcal{M}^{\rm 1,NA} (X)$ with $\mu = {\rm MA}^{\rm NA}_{-K_X}(\phi_{(\mathcal{X}, \mathcal{L})})$ for a test configuration $(\mathcal{X}, \mathcal{L})$ of $(X,-K_X)$.
\end{itemize}
\end{cor}

\subsection{A coupled analogue of Fujita's anticanonical volume bound}

Let $F$ be a prime divisor over $X$.
For an ample line bundle $L$, the expected vanishing order of $L$ for $F$ is defined by
$$
S_L (F) := \frac{1}{L^n} \int_0^{\infty} {\rm vol } (\pi^* L - x F) dx,
$$
where $F \subset Y \xrightarrow{\pi} X$ (\cite{Fj} and see also \cite[\S 7.1 and 7.2]{BBJ}).
If $L$ is a ample $\mathbb{Q}$-line bundle, we set $S_L (F) := m^{-1} S_{mL} (F)$ where $m$ is a sufficiently divisible integer such that $mL_i$ is a line bundle.

\begin{thm}
\label{cp Fj 1}    Assume that $F^{\rm cp ,NA}_{-1} (\delta_{v_F}) \geq 0$ for any Dirac mass at a divisorial point $v_F \in X^{\rm div}$.
    Then, we have
    $$
    (L_1^n)^{\frac{1}{n}} + ... + (L_m^n)^{\frac{1}{n}} \leq n+1.
    $$
    In particular, if $L_i$ are ample line bundles, then $m$ must be less than or equal to $n+1$.
\end{thm}

\begin{proof}
The proof of this theorem follows the proof by Fujita \cite{Fj}.
Let $F \subset Y \xrightarrow{\pi} X$ be a blowing up at a point.
We take a sufficiently divisible integer $m $ such that $m L_i$ are line bundles for any $i$.
Then, by \cite[Lemma 9.4]{BBJ}, we have
$$
n = A_X (F) \geq \frac{1}{m}\sum_i E^{\rm NA}_{mL_i} (\delta_{v_F}) = \frac{1}{m} \sum_i S_{mL_i} (F)=  \sum_i S_{L_i} (F),
$$
By \cite[Theorem 2.3]{Fj} for each $L_i$, we can compute as follows:
\begin{eqnarray*}
    n &\geq&  \sum_{i=1}^{m}   \frac{1}{L_i^n} \int_0^{\infty} {\rm vol } (\pi^* L_i - x F) dx \\
    &\geq&   \sum_{i=1}^{m}   \frac{1}{L_i^n} \int_0^{(L_i^n)^{\frac{1}{n}}}{\rm vol } (\pi^* L_i - x F)dx \\
    &\geq&   \sum_{i=1}^{m}   \frac{1}{L_i^n} \int_0^{(L_i^n)^{\frac{1}{n}}} (L_i^n  - x^n )dx \\
     &=&   \sum_{i=1}^{m}   \frac{n}{n+1} (L_i^n)^{\frac{1}{n}}. 
\end{eqnarray*}
Thus, we obtain the desired inequality.    
\end{proof}

\begin{rem}
If $(X,\{L_i\}_i) = (\mathbb{P}^n , \{ \mathcal{O}(k_i)\}_i)$ where $k_1 +... + k_m = n+1$, there trivially exists a coupled {\kah}--Einstein metric and the equality $(L_1^n)^{\frac{1}{n}} + ...+ (L_m^n)^{\frac{1}{n}}=n+1$ holds.
\end{rem}

\begin{cor}
\label{cp Fj 2}
Assume that $F^{\rm cp ,NA}_{-1} \geq 0$ on $\mathcal{M}^{\rm 1,NA} (X)$.
If $(L_1^n)^{\frac{1}{n} } + ... + (L_m^n)^{\frac{1}{n}} =n+1$ and $(-K_X)^{n} \leq (n+1)^n$, then $X$ is isomorphic to $\mathbb{P}^n$ and $L_i = \mathcal{O}(k_i)$ for some $k_i \in \mathbb{Q}_{>0}$.
\end{cor}

\begin{proof}
By the Brunn--Minkowski inequality for ample line bundles, we
obtain
$$
n+1 = (L_1^n)^{\frac{1}{n} } + ... + (L_m^n)^{\frac{1}{n}}  \leq ( (-K_X)^n)^{\frac{1}{n} } \leq n+1.
$$
Hence equality holds throughout. 
The equality case of the Brunn--Minkowski inequality implies that
all line bundles are proportional, that is,
there exist $k_i$ such that $L_i = - k_i K_X$ (\cite[Theorem 2.6]{LX}).
By Proposition \ref{NA proportional}, a Fano manifold $X$ is Ding semistable (see \cite[Lemma 2.9 and Theorem 2.10]{BJKst}).
By Fujita's result \cite{Fj}, the equality $(-K_X)^n = (n+1)^n$ implies that $X$ must be isomorphic to $\mathbb{P}^n$.
\end{proof}

\begin{rem}
For sufficiently large $n$, the following construction gives 
a $2n$-dimensional toric Fano manifold not isomorphic to $\mathbb P^{2n}$, together with a decomposition $-K_X = L_1 + L_2$
into ample $\mathbb R$-line bundles satisfying 
$(L_1^{2n})^{\frac{1}{2n}} + (L_2^{2n})^{\frac{1}{2n}} = 2n+1$.

In Debarre's book \cite[Proposition 5.2]{Deb}, we can find an $n$-dimensional Fano manifold $Y$ whose automorphism group is not reductive such that
$$
(-K_Y)^n \geq \left( \frac{3n^2}{10 \log n} \right)^n .
$$
More precisely, $Y$ is given by the projective bundle over some projective space.
In particular, we can find sufficiently large $n$ such that $(-K_Y)^n \geq (2n+1)^n$.
For $X := Y \times Y$ and $a,b \in (0,1)$, we set $L_1 = - a p_1^* K_Y - b p_2^* K_Y $ and $L_2 = - (1-a) p_1^* K_Y - (1-b) p_2^* K_Y $.
(Note that $-K_X = L_1 +L_2$.)
By the direct computation, we have
\begin{eqnarray*}
    (L_1^{2n})^{\frac{1}{2n}} + (L_2^{2n})^{\frac{1}{2n}} 
    &=& \left( \sqrt{ab}  + \sqrt{(1-a)(1-b)} \right)  \left(  \frac{\sqrt{(2n)!}}{n!} (-K_Y)^n\right)^{\frac{1}{n}} .  
\end{eqnarray*}
If $a = b$, we have $((-K_X)^{2n})^{\frac{1}{2n}}=(L_1^{2n})^{\frac{1}{2n}} + (L_2^{2n})^{\frac{1}{2n}}   \geq 2n+1$.
On the other hand, we can make the value of $(L_1^{2n})^{\frac{1}{2n}} + (L_2^{2n})^{\frac{1}{2n}}$ arbitrarily small as $a \to 0$ and $b \to 1$.
Thus, the intermediate value theorem implies that we can find a decomposition $\{ L_1, L_2 \}$ of $-K_X$ such that $(L_1^{2n})^{\frac{1}{2n}} + (L_2^{2n})^{\frac{1}{2n}}   = 2n+1$ by choosing appropriate $a,b \in (0,1)$.
(Note that $a$ and $b$ are not rational numbers in general.)

\end{rem}

\subsection{Uniform coupled Gibbs stability}

In this subsection, we introduce uniform coupled Gibbs stability,
extending Berman's uniform version of the Gibbs stability \cite{Ber13,Berinv,BerQP} to the coupled setting,
and prove that it implies the existence of a unique coupled K\"ahler–Einstein metric under the equal-volume assumption.

For each sufficiently large and divisible integer $k$, 
we write $N_j = \dim H^0 (X, kL_j)$.
We write $V_j = L_j^n$.
As assumed in \S 3.2, we need the following hypothesis in this subsection.
\begin{hypo}
\label{volume hypothesis}
$$
V_1 = V_2 = \cdots = V_m.
$$
\end{hypo}

It is well-known that $ N_j =\frac{ V_j}{n!} k^n + \frac{S_j}{(n-1)!} k^{n-1} + O(k^{n-2})$, where $S_j$ is the total scalar curvature of $X$ with respect to $L_j$.
We order $N_j$ so that
$$
N_1 \leq N_2 \leq \cdots \leq N_m 
$$
at level $k$.
For simplicity, we write $N := N_m$.
\begin{rem}
\label{volume comparison}
By Hypothesis \ref{volume hypothesis}, 
$N_i /N \to 1$ for every $i$ as $k\to \infty$.
This fact is important for the perturbation of the parameter $\gamma_0$ and our proof of Theorem \ref{partition and free energy} (see (\ref{purturb})).
\end{rem}
We write the trivial line bundle over $X$ as $\mathscr{O}_X$ and define the line bundle over $X^{N}$ by 
$\left( \otimes_{j=1}^{N_i} p_j^*  L_i^k \right) \otimes \left( \otimes_{j=N_i + 1}^{N} p_j^* \mathscr{O}_X \right)$.
For each $i = 1,2,...,m$, we regard each Slater determinant $\det S_i^{(k)}$ (see (\ref{Slater}) in \S 3.2.1) as a holomorphic section in $H^0 \left( X^{N} , \left( \otimes_{j=1}^{N_i} p_j^*  L_i^k \right) \otimes \left( \otimes_{j=N_i + 1}^{N} p_j^* \mathscr{O}_X \right) \right)$, i.e.,
$$
\pi_i^* \det S_i^{(k)} \left(x_1 , ..., x_{N_i} , x_{N_i +1} ,... , x_{N} \right) \in H^0 \left( X^{N} , \left( \otimes_{j=1}^{N_i} p_j^*  L_i^k \right) \otimes \left( \otimes_{j=N_i + 1}^{N} p_j^* \mathscr{O}_X \right) \right),
$$
where $\pi_i$ denotes the projection map defined by
$$
\pi_i : X^N \to X^{N_i} \,\, ; \,\, \left(x_1 , ..., x_{N_i} , x_{N_i +1} ,... , x_{N} \right) \mapsto \left(x_1 , ..., x_{N_i} \right).
$$

We fix Hermitian metrics $h_i $ on $L_i$ whose curvature forms (multiplied by $\sqrt{-1}$) are the {\kah} forms $\omega_i$ and satisfy $\sum_i^m \omega_i = {\rm Ric}\, dV$.
For $\phi_i \in \mathcal{H}(L_i)$, we can define the Hermitian metric $H(i,k,\phi_i)$ on the line bundle $\left(  \otimes_{j=1}^{N_i} p_j^*  L_i \right) \otimes \left( \otimes_{j=N_i + 1}^{N} p_j^* \mathscr{O}_X \right)$ by
$$
H(i, k, \phi_i) (u,u) \left(x_1 , ..., x_{N_i} , x_{N_i +1} ,... , x_{N} \right) := |u|^2 \left(  \prod_{j=1}^{N_i} h_i e^{-k\phi_i}    (x_j) \right) e^{-k (\phi_i (x_{N_i +1}) + ... + \phi_i (x_{N}))},
$$
where $u$ denotes a local holomorphic frame of $\left(  \otimes_{j=1}^{N_i} p_j^*  L_i \right) \otimes \left( \otimes_{j=N_i + 1}^{N} p_j^* \mathscr{O}_X \right)$.

We define the partition function in the coupled setting as follows.
\begin{definition}
Fix $p>1$.
For $f \in L^p(dV)$ such that $f\geq 0$ and $\int_X f dV =1$, we set $dV_f := f dV $.
For $\boldsymbol{\gamma} = (\gamma_1 ,..., \gamma_m) \in \mathbb{R}^m$, we define the coupled partition function by
$$
\mathcal{Z}_{N,f}^{\rm cp} (- \boldsymbol{\gamma}) := \int_{X^N} \prod_i  \left| \pi_{i}^* \det S_i^{(k)} \right|^{-2\gamma_i /k}_{H(i,k,0 ) }  dV_f^{\otimes N}  \in \mathbb{R}_{>0} \cup\{+\infty\}.
$$
\end{definition}

We can also define the coupled free energy and coupled Ding functional with respect to the (singular) measure $dV_f$.
For simplicity, we use the same symbols $F_{-\gamma}^{\rm cp}$ and $D_{-\gamma}^{\rm cp}$, i.e.,
$$
F_{-\gamma}^{\rm cp}(\mu) := {\rm Ent}_{dV_f} (\mu) - \gamma E^{\rm cp}(\mu) , \,\,\,  D_{-\gamma}^{\rm cp} ((\phi_i)_i) := -\frac{1}{\gamma} \log \int_X e^{-\gamma \sum_i \phi_i } dV_f  - \sum_i \mathcal{E}_{i} (\phi_i) .
$$

We prove the following coupled version of Berman's result \cite[Theorem 1.2]{BerQP}.
\begin{thm}
\label{partition and free energy}
Fix $\gamma_0 > 0$.
We set $\boldsymbol{\gamma} = (\gamma_1,...,\gamma_m)$, where 
\begin{equation}
\label{gamma}
\gamma_i := \frac{\gamma_0}{1-\epsilon_i} \frac{N}{N_1}
\end{equation}
for sufficiently small $\epsilon_i>0 ,\,\,i=1 ,2, ....,m$.
Then, there exists a positive constant $C_0 >  0$ such that for any sufficiently large integer $k \in \mathbb{Z}_{>0}$ and small $\tau >0$, we have
\begin{eqnarray*}
&& \hspace{-30pt}-\frac{1}{ N  (1+\tau)} \log \mathcal{Z}^{\rm cp}_{N,f} (-(1+\tau)\boldsymbol{\gamma}) \\ && \leq  F_{-(1 -k^{-1}C_0)\gamma_0}^{\rm cp}\left( \mu \right)  + \frac{k\tau-(1+\tau)M}{k \gamma_0 (1+\tau)} \log \int_{X} f^{k\tau / (k\tau-M \gamma(1+\tau))}  dV + O\left(\frac{ \log k}{k} \right).
\end{eqnarray*}
Here, we set $M := \sum_i \gamma_i$.


\end{thm}
Before proving this theorem, we introduce the notion of uniform Gibbs stability in the coupled setting.

\begin{definition}
The coupled microscopic stability threshold for $(X, \{ L_i \}_i,f)$ at level $k$ is defined by
$$
\gamma_{k,f}^{\rm cp}\left( X, \{ L_i \}_{i=1}^{m} \right) := \sup \left\{ \gamma \geq 0 \, \middle| \,  \mathcal{Z}_{N,f}^{\rm cp}(-(\gamma,..., \gamma)) < + \infty \right\} .
$$
We define the asymptotic threshold by
$$
\gamma_f^{\rm cp} \left( X, \{ L_i \}_{i=1}^{m} \right) := \liminf_{k \to \infty} \gamma_{k,f}^{\rm cp} \left( X, \{ L_i \}_{i=1}^{m} \right) .
$$
We say that $(X, \{L_i\}_i)$ is uniformly coupled Gibbs stable if $$ \gamma^{\rm cp} \left( X, \{ L_i \}_{i=1}^{m} \right) := \gamma_1^{\rm cp} \left( X, \{ L_i \}_{i=1}^{m} \right) >1.$$
\end{definition}

The microscopic stability threshold and uniform Gibbs stability were defined by Berman \cite{Ber18,Berinv} in order to study {\kah}--Einstein metrics on Fano manifolds.
Fujita \cite{Fj16} and Fujita--Odaka \cite{FO} proved that if a Fano manifold is uniformly Gibbs stable then it is (uniformly) K-stable, so it admits a {\kah}--Einstein metric.
After their works, Berman proved that a uniformly Gibbs stable Fano manifold has a {\kah}--Einstein metric by showing the coercivity of the free energy functional.
As a coupled version of \cite[Theorem 2.5]{FO} and \cite[Corollary 3.4]{BerQP}, we obtain the following by the coercivity of the coupled free energy (equivalently, the coupled Ding functional) and Kewei Zhang's result \cite{KZ}. 

\begin{cor}

\label{ext of BFO2}
The inequality $ \delta^{\rm cp} \left( X, \{ L_i \}_{i=1}^{m} \right) \geq \gamma^{\rm cp} \left( X, \{ L_i \}_{i=1}^{m} \right)$ holds.
In particular, if $(X , \{L_i \}_i)$ is uniformly coupled Gibbs stable, then there exists a unique coupled {\kah}--Einstein metric.
\end{cor}
\begin{proof}
For any $\gamma \in (0, \gamma^{\rm cp} \left( X, \{ L_i \}_{i=1}^{m} \right))$, we take $\gamma_0 < \gamma$.
Then, we can find sufficiently small $\tau, \epsilon_i$ and large $k$ such that $(1+\tau)\gamma_i < \gamma $ for all $i$.
Here, $\gamma_i$ is defined by (\ref{gamma}).
Theorem \ref{partition and free energy} implies that the coupled free energy $F_{-(1- k^{-1}C_0)\gamma_0}^{\rm cp}$ is bounded from below.
So, by Theorem \ref{coercivity}, the coupled Ding functional  $D_{-(1- k^{-1}C_0)\gamma_0 + d}^{\rm cp}$ is coercive for sufficiently small $d>0$.
The coupled delta-invariant $ \delta^{\rm cp} \left( X, \{ L_i \}_{i=1}^{m} \right)$ coincides with the coercivity threshold of the coupled Ding functional, so we obtain $(1-k^{-1} C_0) \gamma_0 - d <   \delta^{\rm cp} \left( X, \{ L_i \}_{i=1}^{m} \right)$.
As $k\to \infty$ and $d \to 0$, we have $\gamma_0 \leq \delta^{\rm cp} \left( X, \{ L_i \}_{i=1}^{m} \right)$, so we obtain $\gamma \leq \delta^{\rm cp} \left( X, \{ L_i \}_{i=1}^{m} \right)$.
Thus, we have proved the inequality $ \gamma^{\rm cp} \left( X, \{ L_i \}_{i=1}^{m} \right) \leq  \delta^{\rm cp} \left( X, \{ L_i \}_{i=1}^{m} \right)$.
\end{proof}

We now prove Theorem \ref{partition and free energy},
following the refinement in \cite{Ao} of Berman's argument \cite{BerQP}.

\medskip
\hspace{-17pt}{\it Proof of Theorem \ref{partition and free energy}.}
For a fixed finite energy measure $\mu$, we solve $\mu = {\rm MA}_{\omega_i} (\phi_i) $ for each $i$.
It suffices to prove the inequality in the case when $\mu$ is a volume form, i.e., all $\phi_i$ are smooth, by Remark \ref{reduction to sm}.

Firstly, we assume that for $\phi_i \in \mathcal{H}(L_i)$, the measure $e^{-\sum_{i=1}^{m} \gamma_i \phi_i} dV_f$ is a probability measure.
We consider the following normalization for $\phi_i \in \mathcal{H} (L_i)$:
$$
\int_X e^{-k \phi_i} dV =1 ,\,\,\,  (1 \leq i \leq m-1).
$$

For $\tau > 0$, the {\hol} inequality implies that
{\small
\begin{eqnarray*}
1&=& \int_{X^N}  (e^{- \sum_{i=1}^{m}\gamma_i \phi_i} dV_f)^{\otimes N} \\
&=& \int_{X^N}  \left( \prod_i  \left|\pi_i^* \det S_i^{(k)} \right|^{-2\gamma_i /k}_{H (i,k,\phi_i  ) } \prod^N e^{-\tau\sum_i \gamma_i  \phi_i/(1 + \tau)}  \right)\\
&& \hspace{50pt}\times  \left( \prod_i \left| \pi_i^* \det S_i^{(k)} \right|^{2\gamma_i /k}_{H(i,k,\phi_i  )  }  \prod^N e^{\tau  \sum_i \gamma_i\phi_i/(1 + \tau)} \right) (e^{-\sum_i \gamma_i \phi_i} dV_f)^{\otimes N}  \\
&\leq& \left( \int_{X^N} \prod_i \left| \pi_i^* \det S_i^{(k)} \right|^{-2\gamma_i (1+\tau)/k}_{H (i,k,\phi_i  ) }  \prod^N e^{-\tau \sum_i \gamma_i \phi_i}  (e^{-  \sum_i  \gamma_i \phi_i} dV_f)^{\otimes N} \right)^{1 /(1+\tau)} \\
&& \hspace{50pt}\times \left( \int_{X^N} \prod_i  \left| \pi_i^* \det S_i^{(k)} \right|^{2\gamma_i (1 + \tau)/k\tau}_{H(i,k,\phi_i  )  } \prod^N e^{ \sum_i\gamma_i  \phi_i}  (e^{- \sum_i \gamma_i \phi_i} dV_f)^{\otimes N} \right)^{\tau /(1 + \tau)}\\
&=& \left( \int_{X^N} \prod_i \left| \pi_i^* \det S_i^{(k)} \right|^{-2\gamma_i (1+\tau)/k}_{H(i,k,\phi_i  )  }  (e^{-(1+\tau) \sum_i \gamma_i \phi_i} dV_f)^{\otimes N} \right)^{1 /(1+\tau)} \\
&& \hspace{50pt}\times \left( \int_{X^N} \prod_i \left| \pi_i^*\det S_i^{(k)} \right|^{2\gamma_i (1 + \tau)/k\tau}_{H(i,k,\phi_i  ) } dV_f^{\otimes N} \right)^{\tau /(1 + \tau)}\\
&=& \left(\mathcal{Z}_{N,f}(- (1+\tau)\boldsymbol{\gamma})^{\rm cp}\right)^{1/(1+\tau)} 
\times \left( \int_{X^N} \prod_i \left| \pi_i^*\det S_i^{(k)} \right|^{2\gamma_i (1 + \tau)/k\tau}_{H(i,k,\phi_i  )  } dV_f^{\otimes N} \right)^{\tau /(1 + \tau)}.
\end{eqnarray*}
}

We write $F(x_1, ... ,x_N) := \prod_j  f (x_j)$.
Note that the assumptions that $f \in L^p (dV)$ and $\int_X dV_f = \int_X f dV = 1$ imply that $F \in L^p (dV^{\otimes N})$ and  $\int_{X^N } F (dV)^{\otimes N} =1$.
We set $M := \sum_j \gamma_j>0$.

We fix a sufficiently large integer $k$ so that $\frac{k\tau}{(1+\tau)M}>1$ and $1<\frac{k\tau }{ k\tau- (1+\tau)M} < p$.
The {\hol} inequality implies that
\begin{eqnarray*}
&&\hspace{-10pt} \int_{X^N} \prod_i \left| \pi_i^* \det S_i^{(k)} \right|^{2\gamma_i (1 + \tau)/k\tau}_{H(i,k, \phi_i) } dV_f^{\otimes N}\\
&&= \int_{X^N} \prod_i \left| \pi_i^* \det S_i^{(k)} \right|^{2\gamma_i (1 + \tau)/k\tau}_{H(i,k, \phi_i)} F dV^{\otimes N}\\
&&\leq \left( \int_{X^N} \prod_i \left| \pi_i^* \det S_i^{(k)} \right|^{2\gamma_i/M}_{H(i,k,\phi_i) } dV^{\otimes N} \right)^{ (1 + \tau) M /k\tau} \left( \int_{X^N} F^{k\tau / (k\tau- (1+\tau)M)} dV^{\otimes N} \right)^{(k\tau-(1+\tau)M)/k\tau} \\
&&= \left\Vert \prod_i \left| \pi_i^* \det S_i^{(k)} \right|^{\gamma_i/M}_{H(i,k,\phi_i) }  \right\Vert^{2(1 + \tau) M/k\tau}_{L^{2}(dV^{\otimes N})} \left( \int_{X^N} F^{k\tau / (k\tau- (1+\tau)M)} dV^{\otimes N} \right)^{(k\tau-(1+\tau)M)/k\tau} .
\end{eqnarray*}

Thus, we obtain the following.
\begin{eqnarray*}
&&\hspace{-20pt} -\frac{1}{(1+\tau) } \log \mathcal{Z}_{N, f}^{\rm cp}(-( 1+\tau)\boldsymbol{\gamma} ) \\
&&\leq \frac{M}{k} \log \left\Vert \prod_i \left| \pi_i^* \det S_i^{(k)} \right|^{\gamma_i/M}_{H_{i,k,\phi_i} }  \right\Vert^{2}_{L^{2}(dV^{\otimes N})} + \frac{N(k\tau-(1+\tau)M)}{ k (1+\tau)} \log \int_{X} f^{k\tau / (k\tau-(1+\tau)M)}  dV.
\end{eqnarray*}
By the definition of the Hermitian metrics $H(i,k,\phi_i)$, we have
$$
\prod_{i=1}^{m} \left| \pi_i^* \det S_i^{(k)} \right|^{2\gamma_i }_{H(i,k,\phi_i) }  = \prod_{i=1}^{m} \left| \pi_i^* \det S_i^{(k)} \right|^{2\gamma_i}_{H(i,k,0) } \prod_{i=1}^{m} e^{-k \sum_{j=1}^{N} \gamma_i \phi_i (x_j)}.
$$

We remove the normalization assumption 
on the measure $e^{-\sum_{i=1}^{m}\gamma_i \phi_i} dV_f$
by replacing $\phi_m$ with
$$
 \phi_m + \frac{1}{\gamma_m}  \log \int_X e^{- \sum_{i=1}^{m} \gamma_i \phi_i} dV_f.
$$
Then we have
\begin{eqnarray*}
&&\hspace{-20pt}-\frac{1}{(1+\tau)} \log \mathcal{Z}_{N,f}^{\rm cp} (- (1+\tau)\boldsymbol{\gamma})\\
&&\hspace{-10pt}\leq \frac{M}{k} \log \left\Vert \prod_i \left| \pi_i^* \det S_i^{(k)} \right|^{\gamma_i/M}_{H(i,k,\phi_i) }  \right\Vert^{2}_{L^{2}(dV^{\otimes N})} - N \log \int_X e^{-\sum_{i=1}^{m} \gamma_i \phi_i} dV_f \\
&& \hspace{70pt} + \frac{N(k\tau-M(1+\tau))}{ k (1+\tau)} \log \int_{X} f^{k\tau / (k\tau-(1+\tau)M)}  dV. 
\end{eqnarray*}
Since $ \frac{\gamma_1}{M} + \ ... + \frac{\gamma_m}{M} = 1$, the {\hol} inequality yields that
\begin{eqnarray*}
\int_{X^N} \prod_{i=1}^{m} \left| \pi_i^* \det S_i^{(k)} \right|^{2\gamma_i /M}_{H(i,k,\phi_i)} dV^{\otimes N}
&\leq& \prod_{i=1}^{m} \left( \int_{X^N} \left| \pi_i^* \det S_i^{(k)} \right|^{2}_{H(i,k,\phi_i)} dV^{\otimes N}  \right)^{\gamma_i /M}.
\end{eqnarray*}
Fubini's theorem and the normalization
$$
\int_X e^{-k \phi_i} dV =1 ,\,\,\, \, i=1, 2,...,m-1,
$$
gives
\begin{eqnarray*}
\int_{X^N} \left| \pi_i^* \det S_i^{(k)} \right|^{2 }_{H(i,k,\phi_i)} dV^{\otimes N}
&=& \int_{X^{N^i}}\left| \det S_i^{(k)} \right|^{2 }_{(h_i e^{-\phi_i})^k} dV^{\otimes N^i} .
\end{eqnarray*}
Then we have
\begin{eqnarray}
\label{PART}
\nonumber&&\hspace{-20pt}-\frac{1}{N  (1+\tau)} \log \mathcal{Z}^{\rm cp}_{N,f} (-(1+\tau)\boldsymbol{\gamma})\\
\nonumber&&\hspace{-10pt}\leq \sum_{i=1}^{m} \frac{\gamma_i}{Nk} \log \left\Vert \det S_i^{(k)} \right\Vert^{2}_{L^{2}((h_i e^{-\phi_i} )^k , dV^{\otimes N_i})} -  \log \int_X e^{- \sum_{i=1}^{m} \gamma_i \phi_i} dV_f \\
\nonumber&& \hspace{70pt} + \frac{k\tau-(1+\tau)M}{k (1+\tau)} \log \int_{X} f^{k\tau / (k\tau-M \gamma(1+\tau))}  dV \\
\nonumber&&\hspace{-10pt}\leq  - \sum_{i=1}^{m} \frac{N_i}{N} \gamma_i \mathcal{E}_{i,k}(\phi_i) -   \log \int_X e^{- \sum_{i=1}^{m} \gamma_i \phi_i} dV_f \\
&& \hspace{70pt} + \frac{k\tau-(1+\tau)M}{k (1+\tau)} \log \int_{X} f^{k\tau / (k\tau-(1+\tau)M)}  dV + \sum_i \frac{\gamma_i \log ( N_i!)}{Nk}.
\end{eqnarray}
Here, the symbol $\mathcal{E}_{i,k}$ denotes the quantized Aubin--Mabuchi energy of the Hermitian matrix $(s^i_\alpha, s_\beta^i)_{L^2 ((h_ie^{-\phi_i})^k,dV)}$:
$$
\mathcal{E}_{i,k} (\phi_i) = - \frac{1}{kN_i} \log \det (s^i_\alpha, s_\beta^i)_{L^2 ((h_ie^{-\phi_i})^k,dV)}
$$
where $(s_\alpha^i) _{\alpha =1}^{N_i}\subset H^0 (X, L_i^k)$ denotes a fixed (orthonormal) basis (see \cite{RTZ}) and we have used \cite[Lemma 5.3]{BB10} (see also the proof of \cite[Proposition 2.3]{BerQP}).
By using the weak partial $C^0$-estimate for each $i$ (see \cite{Ao,BerQP,KZ}), we can find positive constants $C_i>0$ such that we have
\begin{eqnarray*}
-\frac{1}{1-\epsilon_i} \mathcal{E}_{i,k} \left( \phi_{i}^{(\epsilon_i )} \right) \leq - \mathcal{E}_{i} (\phi) + \frac{C_i}{k} \left(- \mathcal{E}_{i} (\phi_i) + \sup_X \phi_i \right) +O(k^{-1}),
\end{eqnarray*}
for any sufficiently small $\epsilon_i >0$ and $\phi_i$.
Here $\phi_i^{(\epsilon_i)} := (1-\epsilon_i) \phi_i $.
The implicit constant in the $O(k^{-1})$ term is independent of $\phi_i$.

Recall that we set
$$
\gamma_i : = \frac{\gamma_0}{1 - \epsilon_i} \frac{N}{N_1}
$$
for the fixed $\gamma_0 \in \mathbb{R}$. 
For simplicity, we write $\gamma_i^{(\epsilon_i)} = (1-\epsilon_i) \gamma_i$.
By replacing $\phi_i$ with $\phi_i^{(\epsilon_i)}$ in the inequality (\ref{PART}), we can compute as follows:
\begin{eqnarray*}
&&\hspace{-20pt}-\frac{1}{N  (1+\tau)} \log \mathcal{Z}_{N,f}^{\rm cp} (- (1+\tau)\boldsymbol{\gamma})\\
&&\hspace{-10pt}\leq  - \sum_{i=1}^{m} \frac{N_i}{N}\gamma_i \mathcal{E}_{i,k}\left( \phi_i^{(\epsilon_i)} \right) -   \log \int_X e^{- \sum_{i=1}^{m} \gamma_i \phi_i^{(\epsilon_i)}} dV_f \\
&& \hspace{70pt} + \frac{k\tau-(1+\tau)M}{k (1+\tau)} \log \int_{X} f^{k\tau / (k\tau-(1+\tau)M)}  dV+ \sum_i \frac{\gamma_i \log ( N_i!)}{Nk} \\
&&\hspace{-10pt}\leq  - \sum_{i=1}^{m} \frac{N_i}{N}  \gamma_i^{(\epsilon_i)} \frac{1}{1-\epsilon_i}\mathcal{E}_{i,k}\left( \phi_i^{(\epsilon_i)} \right) -   \log \int_X e^{- \sum_{i=1}^{m} \gamma_i^{(\epsilon_i)}\phi_i }dV_f \\
&& \hspace{70pt} + \frac{k\tau-(1+\tau)M}{k (1+\tau)} \log \int_{X} f^{k\tau / (k\tau-(1+\tau)M)}  dV + \sum_i \frac{\gamma_i \log ( N_i!)}{Nk}\\
&&\hspace{-10pt}\leq  -  \sum_{i=1}^{m}\frac{N_i}{N}  \gamma_i^{(\epsilon_i)}  \mathcal{E}_{i}\left( \phi_i \right) -   \log \int_X e^{- \sum_{i=1}^{m} \gamma_i^{(\epsilon_i)}\phi_i }dV_f + k^{-1} \sum_{i=1}^{m} \frac{N_i}{N}  \gamma_i^{(\epsilon_i)}  C_i \left( - \mathcal{E}_i (\phi_i) + \sup_{X} \phi_i  \right)\\
&& \hspace{70pt} + \frac{k\tau-(1+\tau)M}{k (1+\tau)} \log \int_{X} f^{k\tau / (k\tau-(1+\tau)M)}  dV + \sum_i \frac{\gamma_i \log ( N_i!)}{Nk} +O(k^{-1}).
\end{eqnarray*}
For the first and second terms above, we also have
\begin{eqnarray*}
&&\hspace{-20pt} -  \sum_{i=1}^{m}\frac{N_i}{N}  \gamma_i^{(\epsilon_i)}  \mathcal{E}_{i}\left( \phi_i \right) -   \log \int_X e^{- \sum_{i=1}^{m} \gamma_i^{(\epsilon_i)}\phi_i }dV_f \\
&&\hspace{-0pt}= -  \sum_{i=1}^{m}\frac{N_1}{N}  \gamma_i^{(\epsilon_i)}  \mathcal{E}_{i}\left( \phi_i \right) -   \log \int_X e^{- \sum_{i=1}^{m} \gamma_i^{(\epsilon_i)}\phi_i }dV_f -  \sum_{i=1}^{m}\frac{N_i- N_1}{N}  \gamma_i^{(\epsilon_i)}  \mathcal{E}_{i}\left( \phi_i \right)  \\
&&\hspace{-0pt}\leq -  \sum_{i=1}^{m}\frac{N_1}{N}  \gamma_i^{(\epsilon_i)}  \mathcal{E}_{i}\left( \phi_i \right) -   \log \int_X e^{- \sum_{i=1}^{m} \frac{N_1}{N} \gamma_i^{(\epsilon_i)}\phi_i }dV_f \\
&& \hspace{20pt}+ \sum_{i=1}^{m}\gamma_i^{(\epsilon_i)} \frac{N_i- N_1}{N}  \left( \sup \phi_i-   \mathcal{E}_{i}\left( \phi_i  \right)\right)  \\
&&\hspace{-0pt}\leq  \gamma_0  \left( - \sum_{i=1}^{m}  \mathcal{E}_{i}\left( \phi_i \right) -  \frac{1}{\gamma_0} \log \int_X e^{- \gamma_0 \sum_{i=1}^{m} \phi_i }dV_f \right) + \sum_{i=1}^{m}\gamma_i^{(\epsilon_i)} \frac{N_i- N_1}{N}  E_i \left( {\rm MA}_{\omega_{i}}( \phi_i)  \right) + O(k^{-1})\\
\end{eqnarray*}
Here, $E_i$ denotes the pluricomplex energy with respect to $L_i$ and we have used the fact that $- \mathcal{E}_i(\phi_i) + \sup_X \phi_i \leq nE_i({\rm MA}_{\omega_i}(\phi_i)) + c$ for some uniform constant $c>0$ depending only on $(X,\omega_i)$ (see \cite[Lemma 2.1]{BerQP} and  \cite[Proposition 2.7]{GZ} ).
By Hypothesis \ref{volume hypothesis}, we can find a uniform constant $C_0$ such that
\begin{equation}
\label{purturb}
0 \leq (N_i - N_1)/N  \leq C_0 k^{-1}    
\end{equation}
for  any $i$.
Thus, we obtain
{\small
\begin{eqnarray*}
&&\hspace{-40pt}-\frac{1}{\gamma_0N  (1+\tau)} \log \mathcal{Z}_{N,f}^{\rm cp} (- (1+\tau)\boldsymbol{\gamma})\\
&&\hspace{-20pt}\leq  \mathcal{D}_{-\gamma_0}^{\rm cp}\left( ( \phi_i)_i \right)   + k^{-1} \sum_{i=1}^{m} C_0  E_i \left( {\rm MA}_{\omega_{i}}( \phi_i)  \right) \\
&& \hspace{20pt}+ \frac{k\tau-(1+\tau)M}{k \gamma_0 (1+\tau)} \log \int_{X} f^{k\tau / (k\tau-(1+\tau)M)}  dV + \sum_i \frac{\gamma_i \log ( N_i!)}{Nk} + O(k^{-1}). \\
\end{eqnarray*}}
Then, we can compute as follows.
\begin{eqnarray*}
&&\hspace{-20pt}-\frac{1}{ N  (1+\tau)} \log \mathcal{Z}_{N,f}^{\rm cp} (- (1+\tau)\boldsymbol{\gamma})\\
&&\hspace{-10pt}\leq  \gamma_0 \mathcal{D}_{-\gamma_0}^{\rm cp}\left(( \phi_i )_i \right)   + k^{-1}C_0 \gamma_0 \sum_{i=1}^{m} E_i \left( {\rm MA}_{\omega_{i}}( \phi_i)\right)  \\
&&\hspace{30pt}+ \frac{k\tau-(1+\tau)M}{k \gamma_0 (1+\tau)} \log \int_{X} f^{k\tau / (k\tau-(1+\tau)M)}  dV+ \sum_i \frac{\gamma_i \log ( N_i!)}{Nk}+O(k^{-1}) \\
&&\hspace{-10pt}\leq  F_{-\gamma_0}^{\rm cp}\left( \mu \right)   + k^{-1}C_0 \gamma_0 \sum_{i=1}^{m} E_i \left( {\rm MA}_{\omega_{i}}( \phi_i)\right)  + \frac{k\tau-(1+\tau)M}{k \gamma_0 (1+\tau)} \log \int_{X} f^{k\tau / (k\tau-(1+\tau)M)}  dV+O\left( \frac{\log k}{k} \right)\\
&&\hspace{-10pt}\leq  F_{-(1 -k^{-1}C_0)\gamma_0}^{\rm cp}\left( \mu \right) + \frac{k\tau-(1+\tau)M}{k \gamma_0 (1+\tau)} \log \int_{X} f^{k\tau / (k\tau-(1+\tau)M)}  dV + O\left( \frac{\log k}{k} \right)\\
\end{eqnarray*}

Here, we have used the equations $\mu = {\rm MA}_{\omega_i} (\phi_i)$ for any $i$, Proposition \ref{F geq G} and Stirling's formula.\sq

\medskip


\end{document}